\documentclass[a4paper,hidelinks,onefignum,onetabnum]{siamart251216}
\usepackage{geometry}
\usepackage{amsfonts}
\usepackage{graphicx}
\usepackage{epstopdf}
\usepackage{algorithmic}
\usepackage{booktabs}
\usepackage{multirow}
\usepackage{threeparttable}
\ifpdf
  \DeclareGraphicsExtensions{.eps,.pdf,.png,.jpg}
\else
  \DeclareGraphicsExtensions{.eps}
\fi

\newsiamremark{Hypothesis}{Hypothesis}
\crefname{Hypothesis}{Hypothesis}{Hypotheses}
\newsiamthm{claim}{Claim}
\newsiamremark{fact}{Fact}
\crefname{fact}{Fact}{Facts}

\numberwithin{equation}{section}

\newtheorem{Lemma}{Lemma}[section]
\newtheorem{Theorem}{Theorem}[section]
\newtheorem{Remark}{Remark}[section]
\newtheorem{Proposition}{Proposition}[section]

\headers{Optimal Error Estimates with  
Polynomial Dependence on $\epsilon^{-1}$ for CHNS system}{Haijun Gao, Xi Li, Xiaoli Li}

\title{Optimal Error Estimates with  
Polynomial Dependence on $\epsilon^{-1}$
for Mixed Finite Element Methods of the Cahn--Hilliard--Navier--Stokes System \thanks{Submitted to the editors DATE.
\funding{This work is supported by the National Natural Science Foundation of China (Grant Nos. 12671476, 12601824, 12271302, 12131014) and Shandong Provincial Natural Science Foundation for Outstanding Youth Scholar (Grant No. ZR2024JQ030), the Project SDCX-ZG-202603014 Funded by Postdoctoral Innovation Program of Shandong Province, Natural Science Foundation of Sichuan Province (Grant No. 2025ZNSFSC0070) }}}

\author{Haijun Gao\thanks{School of Mathematics, Shandong University, Jinan, Shandong, 250100, China. (\email{gaohaijun@sdu.edu.cn}).}
\and Xi Li\thanks{School of Mathematical Sciences \& Geomathematics Key Laboratory of Sichuan Province, Chengdu University of Technology, Chengdu, Sichuan 610059, China. (\email{lixi@cdut.edu.cn}).}
\and Xiaoli Li\thanks{Corresponding author. School of Mathematics, Shandong University, Jinan, Shandong, 250100, China. (\email{xiaolimath@sdu.edu.cn}).}
}

\ifpdf
\hypersetup{
	pdftitle={chns},
	pdfauthor={haijun gao}
}
\fi

\begin{document}
	
	\maketitle
	
	\begin{abstract}
		In this paper, we analyze a first-order convex-splitting mixed finite element scheme for
		the Cahn--Hilliard--Navier--Stokes system, tracking explicitly how the error constants depend on
		the interface width $\epsilon$. Due to the strong nonlinear coupling inherent in the CHNS system and the singular dependence of the Cahn--Hilliard equation on $\epsilon$, conventional energy-based analyzes often lead to severe dependence of the error constants on $\epsilon^{-1}$, and may even result in exponential growth. Such dependence significantly limits the applicability of the resulting error estimates in the small-interfacial-thickness regime. To overcome this difficulty, we develop a robust parameter error analysis framework with respect to interfacial thickness $\epsilon$. In particular, by introducing $H^{-1}$ approximation properties for the Ritz projection and the Ritz and Stokes quasi-projection operators, together with spectral-type estimates for the discrete perturbed Cahn--Hilliard operator, we effectively control the singular terms arising from the interfacial thickness parameter in the error analysis. Compared with the classical results, which address only the Allen–Cahn or Cahn–Hilliard equation, we adopt a three-term discrete Gagliardo–Nirenberg inequality analysis technique and achieve a substantial reduction in the degree of polynomial dependence for the coupled system. Numerical experiments are presented to verify the effectiveness of the theoretical results and to demonstrate the robustness of the method for different  parameters of interfacial thickness.
	\end{abstract}
	
	\begin{keywords}
		Cahn--Hilliard--Navier--Stokes system, finite element methods, quasi-projection, Cahn--Hilliard operator, optimal error estimates, polynomial dependence
	\end{keywords}
	
	\begin{MSCcodes}
		65M12, 65N30, 65M60, 35K55
	\end{MSCcodes}
		
	\section{Introduction}
	We consider the following Cahn–Hilliard–Navier–Stokes (CHNS) model \cite{2006_FengXiaobing_Fully_discrete_finite_element_approximations_of_the_Navier_Stokes_Cahn_Hilliard_diffuse_interface_model_for_two_phase_fluid_flows} in a convex bounded polygonal or polyhedral domain $\Omega\subset\mathbb{R}^d,~d=2,3$, for $t\in (0,T]$:
	\begin{subequations}
		\label{eq_CHNS}
		\begin{flalign}
			\label{eq_CHNS_phi_con}
			\phi_t+\nabla\cdot(\phi\boldsymbol{u})-M\varDelta\mu=0,\\
			\label{eq_CHNS_mu_con}
			\mu+\epsilon\varDelta\phi-\frac{1}{\epsilon}F'(\phi)=0,\\
			\label{eq_CHNS_ns_con}
			\boldsymbol{u}_t+(\boldsymbol{u}\cdot\nabla)\boldsymbol{u}-\nu\varDelta\boldsymbol{u}+\nabla p-\gamma\mu\nabla\phi=\boldsymbol{f},\\
			\label{eq_CHNS_incompressible_u_con}
			\nabla\cdot\boldsymbol{u}=0,
		\end{flalign}
	\end{subequations}
	where the no-flux/no-flow boundary conditions and initial conditions are:
	\begin{equation}
		\frac{\partial\phi}{\partial n}=\frac{\partial\mu}{\partial n}= 0,\quad\text{and}\quad\boldsymbol{u}=0 \quad\text{on}\quad \partial\Omega\times(0,T],
	\end{equation}
	and
	\begin{equation}\label{eq_CHNS_initial_boundary_conditions}
		\phi(x,0)=\phi^0,\quad\text{and}\quad\boldsymbol{u}(x,0)=\boldsymbol{u}^0 \quad\text{in}~\Omega.
	\end{equation}
	Here,
	$\phi$ is the phase-field variable, and $\mu$ is the chemical potential, the variables $\boldsymbol{u}$ and $p$ denote the velocity and pressure of the fluid, respectively. The positive physical parameters in the model are given by mobility $M$, viscosity $\nu$ (the inverse of Reynolds number), 
	$\gamma = \frac{1}{{We}^*}$, where ${We}^*$ denotes the modified Weber number measuring the ratio of kinetic to surface energies, 
	and the dimensionless parameter $\epsilon$ characterizing the interfacial width between the two phases \cite{2012_GarckeThermodynamically_consistent_frame_indifferent_diffuse_interface_models_for_incompressible_two_phase_flows_with_different_densities}. The nonlinear potential $F(\phi)$ is typically chosen as a double-well potential of the form: for $t\in (0,T]$,
	\begin{equation}
		F(\phi)=\frac{1}{4}(\phi^2-1)^2.
	\end{equation}
	For the CHNS system \eqref{eq_CHNS}-\eqref{eq_CHNS_initial_boundary_conditions}, the energy is defined by 
	\begin{equation}\label{eq_energy_functional}
		\mathcal{J}_{\epsilon}(\phi,\boldsymbol{u})=\int_{\Omega}\bigg\{\frac{\gamma}{\epsilon}|F(\phi)|+\frac{\gamma\epsilon}{2}|\nabla\phi|^2+\frac{1}{2}|\boldsymbol{u}|^2\bigg\}d\boldsymbol{x}.
	\end{equation}
	It is straightforward to observe that the CHNS system \eqref{eq_CHNS}-\eqref{eq_CHNS_initial_boundary_conditions} can be viewed as a mass-conserving gradient flow associated with the energy functional \eqref{eq_energy_functional}. Consequently, for any time $t \in [0,T]$, the system satisfies the following energy dissipation law:
	\begin{equation}
		\mathcal{J}_{\epsilon}(\phi(t),\boldsymbol{u}(t))+\int_{0}^{t}\left(\gamma M\|\nabla\mu(s)\|^2+\nu\|\nabla\boldsymbol{u}(s)\|^2\right)ds=\mathcal{J}_{\epsilon}(\phi^0,\boldsymbol{u}^0),
	\end{equation}
	and the mass conservation $\left(\phi(t),1\right)=\left(\phi^0,1\right)$ is maintained. Over the past many years, a large number of numerical methods have been proposed for the CHNS system to preserve its above physical properties, such as finite difference methods \cite{2001_Daisuke_A_stable_and_conservative_finite_difference_scheme_for_the_Cahn_Hilliard_equation,2024_WangCheng_A_second_order_numerical_scheme_of_the_Cahn_Hilliard_Navier_Stokes_system_with_Flory_Huggins_potential,2021_WangCheng_A_second_order_accurate_scalar_auxiliary_variable_SAV_numerical_method_for_the_square_phase_field_crystal_equation}, mixed finite element methods \cite{2006_FengXiaobing_Fully_discrete_finite_element_approximations_of_the_Navier_Stokes_Cahn_Hilliard_diffuse_interface_model_for_two_phase_fluid_flows,2008_Kay_David_and_Welford_Richard_Finite_element_approximation_of_a_Cahn_Hilliard_Navier_Stokes_system,2015_HanDaozhi_A_second_order_in_time_uniquely_solvable_unconditionally_stable_numerical_scheme_for_Cahn_Hilliard_Navier_Stokes_equation,2023_CaiWentao_Optimal_L2_error_estimates_of_unconditionally_stable_finite_element_schemes_for_the_Cahn_Hilliard_Navier_Stokes_system,2026_GaoHaijun_Optimal_L2_error_estimates_of_a_decoupled_FEM_for_the_Cahn_Hilliard_Navier_Stokes_equations}, discontinuous Galerkin (DG) method \cite{2020_LiuChen_A_priori_Error_Analysis_of_a_Discontinuous_Galerkin_Method_for_Cahn_Hilliard_Navier_Stokes_Equations,2023_LiuChen_Convergence_of_a_decoupled_splitting_scheme_for_the_Cahn_Hilliard_Navier_Stokes_system}, hybridizable discontinuous Galerkin (HDG) method \cite{2026_YeChanglun_An_unconditionally_stable_hybrid_discontinuous_Galerkin_scheme_for_a_Cahn_Hilliard_phase_field_model_of_two_phase_incompressible_flow_with_variable_densities,2020_FuGuosheng_A_divergence_free_HDG_scheme_for_the_Cahn_Hilliard_phase_field_model_for_two_phase_incompressible_flow} and Fourier spectral methods \cite{2009_ShenJie_An_efficient_moving_mesh_spectral_method_for_the_phase_field_model_of_two_phase_flows}.
	
	The numerical treatment of the CHNS system builds upon the fundamental theories of the Navier--Stokes (NS) and Cahn--Hilliard (CH) equations. In recent years, a key challenge has been the efficient approximation of nonlinear terms in the CH equation, aiming to simultaneously achieve computational efficiency and energy stability in discrete schemes. A variety of numerical strategies have been developed to address these nonlinearities, including convex splitting \cite{1998_Eyre_David_J_Unconditionally_gradient_stable_time_marching_the_Cahn_Hilliard_equation}, stabilized semi-implicit methods \cite{2010_ShenJie_Energy_stable_schemes_for_Cahn_Hilliard_phase_field_model_of_two_phase_incompressible_flows}, the invariant energy quadratization approach (IEQ) \cite{2017_YangXiaofeng_Numerical_approximations_for_the_molecular_beam_epitaxial_growth_model_based_on_the_invariant_energy_quadratization_method}, and the scalar auxiliary variable (SAV) method \cite{2018_Shenjie_Xujie_SAV}.
	For the NS equations, the primary difficulties lie in handling the nonlinear convective term and the incompressibility constraint coupling velocity and pressure. A widely used technique for decoupling velocity and pressure is the projection method \cite{2006_Guermond_ShenJie_An_overview_of_projection_methods_for_incompressible_flows}, originally introduced by Chorin and Temam \cite{1968_Chorin_Alexandre_Numerical_solution_of_the_Navier_Stokes_equations}. The nonlinear terms in the NS equations can be treated either implicitly or through explicit SAV-based techniques \cite{2022JieShen_LiXiaoliMSAVCHNStwo_phase_incompressible_flows}. Within the CHNS system, substantial research efforts have focused on ensuring both energy stability and computational efficiency \cite{2021_Yang_Xiaofeng_Decoupled_linear_and_unconditionally_energy_stable_fully_discrete_finite_element_numerical_scheme_for_a_two_phase_ferrohydrodynamics_model,2018_GaoYali_Decoupled_linear_and_energy_stable_finite_element_method_for_the_Cahn_Hilliard_Navier_Stokes_Darcy_phase_field_model,2020_YangXiaofeng_Error_Analysis_of_a_Decoupled_Linear_Stabilization_Scheme_for_the_Cahn_Hilliard_Model_of_Two_Phase_Incompressible_Flows,2015_HanDaozhi_A_second_order_in_time_uniquely_solvable_unconditionally_stable_numerical_scheme_for_Cahn_Hilliard_Navier_Stokes_equation,2018_CaiYongyong_Error_estimates_for_a_fully_discretized_scheme_to_a_Cahn_Hilliard_phase_field_model_for_two_phase_incompressible_flows}.
	In addition to spatial discretization of finite elements \cite{2015_Diegel_Analysis_of_a_mixed_finite_element_method_for_a_Cahn_Hilliard_Darcy_Stokes_system,2008_Kay_David_and_Welford_Richard_Finite_element_approximation_of_a_Cahn_Hilliard_Navier_Stokes_system}, numerous studies have employed DG methods for the CHNS system \cite{2023_LiuChen_Convergence_of_a_decoupled_splitting_scheme_for_the_Cahn_Hilliard_Navier_Stokes_system}. For example, Liu et al. \cite{2020_LiuChen_A_priori_Error_Analysis_of_a_Discontinuous_Galerkin_Method_for_Cahn_Hilliard_Navier_Stokes_Equations} proposed and analyzed an interior penalty DG method for the CHNS model. Subsequently, they \cite{2023_LiuChen_Convergence_of_a_decoupled_splitting_scheme_for_the_Cahn_Hilliard_Navier_Stokes_system} developed a decoupled and energy-stable DG scheme for the CHNS equations. Despite the fact that the continuous finite element method is relatively simple and straightforward to implement given the complexity of the CHNS system, the DG method offers local mass conservation at the cost of increased computational effort.
	
	More recently, Xu et al. \cite{2020_YangXiaofeng_Error_Analysis_of_a_Decoupled_Linear_Stabilization_Scheme_for_the_Cahn_Hilliard_Model_of_Two_Phase_Incompressible_Flows} introduced a first-order, linear, fully decoupled and energy-stable numerical scheme for the CHNS model. By incorporating an intermediate velocity variable in the discretization of the CH equation, they clarified the coupling between velocity and phase function; however, this also complicated error analysis. Their study focused on establishing the mathematical theory for the numerical solution in the $L^\infty$-norm.
	Based on this work, Wen et al. \cite{2022_Wenjuan_Semi_implicit_unconditionally_energy_stable_stabilized_finite_element_method_based_on_multiscale_enrichment_for_the_Cahn_Hilliard_Navier_Stokes_phase_field_model} developed a novel fully discrete stabilized finite element method using the lowest equal-order elements, providing optimal $H^1$ error estimates for the phase field variable, chemical potential, and velocity. Cai et al. \cite{2023_CaiWentao_Optimal_L2_error_estimates_of_unconditionally_stable_finite_element_schemes_for_the_Cahn_Hilliard_Navier_Stokes_system} presented two quasi-projection operators, which led to the optimal $L^2$ error estimates for the convex splitting finite element method applied to the CHNS system.
	Later, Yang et al. \cite{2024_ChenYaoyao_SAVFEM} proposed a semi-implicit Euler scheme based on the SAV approach and analyzed the optimal error estimates of the MINI element approximation for the CHNS system. 
	Although substantial progress has been made in the numerical analysis of the CHNS system, the dependence of error constants on the inverse interfacial thickness parameter $\epsilon^{-1}$ remains a major challenge. Due to the singular dependence of the Cahn--Hilliard operator on $\epsilon$ and the strong nonlinear coupling between the phase and fluid variables, conventional energy arguments may lead to exponentially growing error constants in $\epsilon^{-1}$. For the CH equation, spectral-type estimates have been successfully developed to reduce this dependence to polynomial order, but such techniques have not yet been fully extended to the coupled CHNS system within the mixed finite element framework.
	
	In this work, we analyze the optimal error estimates for the classical convex-splitting mixed finite element discretization of the CHNS system. By employing a discrete spectral estimate for the perturbed Cahn--Hilliard operator, we further show that the error constants depend on the inverse interface-width parameter $\epsilon^{-1}$ only polynomially. The main difficulty in the optimal error analysis of the mixed finite element method for the CHNS system lies in the strong coupling among the physical variables, through which error contamination terms arising from the coupling terms may destroy the optimal error estimates. To overcome this difficulty, inspired by the Ritz quasi-projection and Stokes quasi-projection operators \cite{2023_CaiWentao_Optimal_L2_error_estimates_of_unconditionally_stable_finite_element_schemes_for_the_Cahn_Hilliard_Navier_Stokes_system}, we decompose the total error into projection and discrete components, and exploit the quasi-projection structure to eliminate the leading contamination terms arising from the coupling. Combined with the $H^{-1}$ approximation properties of the Ritz and Stokes quasi-projection operators and the discrete spectral estimate for the perturbed Cahn--Hilliard operator, this approach effectively controls the contamination of the velocity error induced by the phase-field and chemical-potential errors through the nonlinear coupling terms, thereby yielding an optimal stability estimate for the coupled error equations. Meanwhile, we theoretically determine the admissible range of the perturbation magnitude in the perturbed Cahn--Hilliard operator, which ensures the validity of the corresponding discrete spectral estimate.
	The conventional optimal error analysis encounters two major difficulties: i) gradient terms of the phase-field error may appear on the right-hand side of the error estimates, which prevents the error estimates from being optimal under the available regularity assumptions. ii) a direct application of Grönwall's inequality to coupled error equations may absorb high powers of $\epsilon^{-1}$ into the stability coefficient, resulting in an exponential dependence of the error constants on $\epsilon^{-1}$. To overcome these difficulties, we combine a quasi-projection technique with a discrete spectral estimate, allowing us to derive optimal error estimates while avoiding exponential dependence of the error constants on the interface-width parameter $\epsilon$. Consequently, we obtain optimal $L^2$ error estimates whose constants depend polynomially on $\epsilon^{-1}$. 
	Furthermore, by choosing the \(L^p\) norms in a new optimal three-term Gagliardo--Nirenberg interpolation inequality rather than employing the classical Gagliardo--Nirenberg inequality, we substantially reduce the degree of polynomial dependence on \(\epsilon\). Compared with the classical analytical results in the literature \cite{2015_FengXiaobing_Analysis_of_symmetric_interior_penalty_discontinuous_Galerkin_methods_for_the_Allen_Cahn_equation_and_the_mean_curvature_flow,2023_ChenGang_On_the_superconvergence_of_a_hybridizable_discontinuous_Galerkin_method_for_the_Cahn_Hilliard_equation}, it is worth noting that those results address the Allen--Cahn or Cahn--Hilliard equation individually and do not provide an analysis for coupled problems, whereas our results pertain to a phase-field model coupled with multiple physical fields and can be extended to the optimal error analysis of other coupled systems.
	
	The rest of this paper is organized as follows. In section \ref{section_Notations}, we present the basic mathematical notation used throughout this paper and introduce several projection operators,  including the Ritz quasi-projection and the Stokes quasi-projection operator and their approximation properties, the spectral estimate for the Cahn-Hilliard operator.  In section \ref{section_numerical_scheme}, we introduce the fully discrete numerical scheme and establish its stability. The section \ref{section_error_estimates} is devoted to the optimal $L^2$ error estimates for the fully discrete scheme, which does not exhibit exponential dependence on $\epsilon^{-1}$. The section \ref{section_numerical_results} provides numerical experiments to validate the theoretical analysis. The conclusion and discussion are given in section \ref{section_conclusions}.
	
	\section{Preliminaries}\label{section_Notations}	
	In this section, we mainly introduce the notation and mixed formulations required for the finite element method applied to the CHNS equations \eqref{eq_CHNS}.
	\subsection{Notations and mixed formulation}
	To introduce the fully discrete mixed finite element methods for the CHNS system, we first give some notation. For $k \geq 0$ and $1 \leq p \leq \infty$, we employ standard Sobolev spaces $W^{k,p}(\Omega)$, 
	with the usual notation $L^p(\Omega)=W^{0,p}(\Omega)$ and $H^k(\Omega)=W^{k,2}(\Omega)$. 
	We further introduce the following scalar function spaces:
	\begin{equation}
		L_0^2(\Omega)=\left\{ v\in L^2(\Omega): \int_{\Omega} v\,dx=0 \right\}, \quad 
		H_0^1(\Omega)=\left\{ v\in H^1(\Omega): v|_{\partial\Omega}=0 \right\},
	\end{equation}
	and define the associated vector-valued spaces by
	\begin{equation}
		\begin{aligned}
			&\boldsymbol{L}^p(\Omega)=[L^p(\Omega)]^d, \quad 
			\boldsymbol{W}^{k,p}(\Omega)=[W^{k,p}(\Omega)]^d, \quad 
			\boldsymbol{H}^k(\Omega)=[H^k(\Omega)]^d,\\
			&\boldsymbol{H}_0^1(\Omega)=[H_0^1(\Omega)]^d, \quad 
			\mathring{H}^1(\Omega)=H^1(\Omega)\cap L_0^2(\Omega).
		\end{aligned}
	\end{equation}
	The dual spaces of $H_0^1(\Omega)$, $\mathring{H}^1(\Omega)$, and $H^1(\Omega)$ are denoted by 
	$H_0^{-1}(\Omega)$, $\mathring{H}^{-1}(\Omega)$, and ${H}^{-1}(\Omega)$, respectively. 
	The norms $\|\cdot\|_{W^{k,p}}$, $\|\cdot\|_k$, and $\|\cdot\|_{L^p}$ represent the standard Sobolev norms, 
	while $\|\cdot\|$ denotes the $L^2$-norm. 
	
	Let $\mathcal{T}_h$ denote a regular and quasi-uniform partition of $\Omega$ into the triangles or tetrahedra $\{\mathcal{K}_j\}, ~j=1,2,\cdots,K$ with the mesh size $h=\max\limits_{0\leq j \leq K}\{\dim\mathcal{K}_j\}$.
	We introduce the finite element spaces as follows: 
	\begin{equation}
		\label{eq_defined_finite_element_space}
		\begin{aligned}
			&S_h^r=\{v_h\in C(\bar{\Omega}):v_h|_{\mathcal{K}_j}\in P_r(\mathcal{K}_j),\forall \mathcal{K}_j\in \mathcal{T}_h\},~
			\mathring{S}_h^r=S_h^r\cap L_0^2(\Omega),\\
			&\boldsymbol{X}_h^{r+1}=\{\boldsymbol{v}_h\in\boldsymbol{H}_0^1(\Omega):\boldsymbol{v}_h|_{\mathcal{K}_j}\in\boldsymbol{P}_{r+1}(\mathcal{K}_j),\forall \mathcal{K}_j\in \mathcal{T}_h\},\\
			&\boldsymbol{X}_h^{1b}=(S_h^1\oplus B_{d+1})^d\cap \boldsymbol{H}_0^1(\Omega),~d=2,3,
		\end{aligned}
	\end{equation}
	where $P_r\left(\mathcal{K}_j\right)$, and $\boldsymbol{P}_r\left(\mathcal{K}_j\right) $ are the spaces of polynomials of total degree $r$ on $\mathcal{K}_j$, $B_3$ and $B_4$ are the spaces of cubic bubbles and quartic bubbles, respectively. 
	We note that the finite element spaces $\boldsymbol{X}_h^{r+1}\times \mathring{S}_h^r$ and  $\boldsymbol{X}_h^{1b}\times \mathring{S}_h^1$ are the generalized Taylor-Hood element and MINI element which satisfy the inf-sup condition \cite{1986_Girault_Vivette_Finite_element_methods_for_Navier_Stokes_equations,2008_Brenner_Susanne_C_The_mathematical_theory_of_finite_element_methods}:
	\begin{equation}\label{eq_inf_sup_condition}
		c\left\|q_h\right\| \leq  \sup_{\boldsymbol{v}_h\neq 0} 
		\frac{\left(q_h, \nabla \cdot \boldsymbol{v}_h\right)}{\left\|\nabla\boldsymbol{v}_h\right\|},~ q_h\in \mathring{S}_h^r,
	\end{equation}
	where $c>0$ is  a positive constant that does not depend on the mesh size $h$. 
	
	In the CHNS model \eqref{eq_CHNS}--\eqref{eq_CHNS_initial_boundary_conditions}, 
	we set the parameters $M$, $\gamma$, and $\nu$ to 1. 
	The extension to the general case with $M,\gamma > 0$ can be carried out 
	without essential difficulty in the analysis, and for the error analysis, we take $\boldsymbol{f} \equiv 0$. 
	From the CHNS system \eqref{eq_CHNS},  it is readily observed that the exact solution $\left(\phi,\mu,\boldsymbol{u},p\right)$ the variational formulation satisfies below:
	\begin{flalign}
		\label{eq_CHNS_variational_form_phi}
		\left(\phi_t,w\right)+b\left(\phi,\boldsymbol{u}, w\right)+\left(\nabla\mu,\nabla w\right)=0,\\
		\label{eq_CHNS_variational_form_mu}
		\left(\mu,\varphi\right)-\epsilon\left(\nabla\phi,\nabla\varphi\right)-\epsilon^{-1}\left(F'(\phi),\varphi\right)=0,\\
		\label{eq_CHNS_variational_form_ns}
		\left(\boldsymbol{u}_t,\boldsymbol{v}\right)+\boldsymbol{B}\left(\boldsymbol{u},\boldsymbol{u},\boldsymbol{v}\right)+\nu\left(\nabla\boldsymbol{u},\nabla\boldsymbol{v}\right)-\left(p,\nabla\cdot \boldsymbol{v}\right)-\left(\mu\nabla\phi,\boldsymbol{v}\right)=0,\\
		\label{eq_CHNS_variational_form_inc}
		\left(\nabla\cdot\boldsymbol{u},q\right)=0,
	\end{flalign}
	for any test functions $\left(w,\varphi,\boldsymbol{v},q\right)\in \left(H^1(\Omega),H^1(\Omega),\boldsymbol{H}_0^1(\Omega),L_0^2(\Omega)\right)$, where
	\begin{flalign}
		&b(\phi,\boldsymbol{u},w):=(\nabla\phi\cdot\boldsymbol{u},w),\\
		&\boldsymbol{B}(\boldsymbol{u},\boldsymbol{v},\boldsymbol{w}):=\frac{1}{2}[(\boldsymbol{u}\cdot\nabla\boldsymbol{v},\boldsymbol{w})-(\boldsymbol{u}\cdot\nabla\boldsymbol{w},\boldsymbol{v})].
	\end{flalign}
	Here, setting $w = 1$ in \eqref{eq_CHNS_variational_form_phi}, the mass conservation $\left(\phi(t),1\right)=\left(\phi^0,1\right)$ is maintained. 
	Firstly, we assume that there exist positive $\epsilon-$independent constants $m_0$ and the nonnegative constants $\sigma_i$ for $i = 1, 2, 3, 4$ such that
	\begin{flalign}
		&m_0:=\frac{1}{|\Omega|}\int_{\Omega}\phi^{0}dx\in (-1,1),\\
		&\mathcal{J}_{\epsilon}(\phi^0,\boldsymbol{u}^0):=\frac{\epsilon}{2}\|\nabla\phi^0\|^2+\frac{1}{\epsilon}\|F(\phi^0)\|_{L^1}+\frac{1}{2}\|\boldsymbol{u}^0\|^2\leq C\epsilon^{-2\sigma_1},\\
		&\|\mu^0\|_{H^l}:=\|-\epsilon\varDelta\phi^0+\frac{1}{\epsilon}f(\phi^0)\|_{H^l}\leq C\epsilon^{-\sigma_{l+2}},~~l=0,1,2,
	\end{flalign} 
	where $f(\phi)$ is defined by $f(\phi)=F'(\phi)=\phi^3-\phi$.
	According to \cite{2006_FengXiaobing_Fully_discrete_finite_element_approximations_of_the_Navier_Stokes_Cahn_Hilliard_diffuse_interface_model_for_two_phase_fluid_flows}, it holds the following proposition:
	\begin{Proposition} \label{Proposition_initial_condition}
		Suppose that the initial values $\phi^0$ and $\boldsymbol{u}^0$ satisfy $|\phi^0|\leq 1$ and $\mathcal{J}_{\epsilon}(\phi^0,\boldsymbol{u}^0)\leq C\epsilon^{-2\sigma_1}<\infty$, i.e., the initial energy is bounded, then the regular solutions $(\phi,\boldsymbol{u},p)$ of the system \eqref{eq_CHNS} satisfy the estimates as follows: for all $T\in (0,\infty)$,
		\begin{flalign}
			\label{eq_mass_conversation0001}
			\int_{\Omega}\phi(x,t)dx=\int_{\Omega}\phi(x,0)dx, ~\forall t\in(0,T),\\
			\label{eq_ennergy_stable00001}
			\operatorname*{ess\,sup}_{t \in [0, T]} \left\{ \| \boldsymbol{u}(t) \|^2 + \epsilon \| \nabla \phi(t) \|^2 +  \epsilon^{-1} (F(\phi(t)), 1) \right\} \leq C\epsilon^{-2\sigma_1},\\
			\label{eq_nabla_u_mu_bounded0001}
			\int_0^T  \|\nabla \boldsymbol{u}(t)\|^2 \, dt + \int_0^T   \|\nabla \mu(t)\|^2 dt\leq C\epsilon^{-2\sigma_1},\\
			\label{eq_phi_t_bounded0001}
			\int_0^T\left(  \|\phi_t(t) + \boldsymbol{u}(t) \cdot \nabla \phi(t)\|_{H^{-1}}^2+\|\nabla\mu(t)\|^2\right)dt \leq C\epsilon^{-2\sigma_1},\\
			\label{eq_phi_t_H_neg1bounded0001}
			\int_{0}^{T}\|\phi_{t}(t)\|_{H^{-1}}^2dt\leq C\epsilon^{-2\sigma_1}.
		\end{flalign}
		In addition, there holds
		\begin{flalign}
			&\operatorname*{ess\,sup}_{t \in [0, T]} \left\{ \begin{aligned} &\|\phi_t\|_{H^{-1}}^2 \\ &\|\nabla \mu\|^2 \end{aligned} \right\} + \epsilon \int_0^\infty \|\nabla \phi_t\|^2 \,\mathrm{d}s \le C\epsilon^{-\max\{2\sigma_1+3,2\sigma_3\}},\\
			\label{eq_Delta_phi_boundedness_continuous}
			&\operatorname*{ess\,sup}_{t \in [0, T]} \|\varDelta \phi\| \le C\epsilon^{-\max\{\sigma_1+5/2,\sigma_3+1\}},\\
			&\operatorname*{ess\,sup}_{t \in [0, T]} \|\nabla \varDelta \phi\| \le C\epsilon^{-\max\{\sigma_1+5/2,\sigma_3+1\}},\\
			&\left\{ \begin{aligned} &\int_0^\infty \|\phi_t\|^2 \,\mathrm{d}s \\ &\int_0^\infty \|\varDelta \mu\|^2 \,\mathrm{d}s \end{aligned} \right\} + \operatorname*{ess\,sup}_{t \in [0, T]} \epsilon\|\varDelta \phi\|^2 \le C\epsilon^{-\max\{2\sigma_1+7/2,2\sigma_2+1,2\sigma_3+1/2\}},\\
			&\operatorname*{ess\,sup}_{t \in [0, T]} (\|\phi_t\|^2 + \|\varDelta \mu\|^2) + \epsilon \int_0^\infty \|\varDelta \phi_t\|^2 \,\mathrm{d}s \le C\epsilon^{-\max\{2\sigma_1+13/2,2\sigma_2+4,2\sigma_3+7/2,2\sigma_4\}}.
		\end{flalign}
		Moreover, according to the GA2 \cite{2015_FengXiaobing_Analysis_of_symmetric_interior_penalty_discontinuous_Galerkin_methods_for_the_Allen_Cahn_equation_and_the_mean_curvature_flow} suppose that there exists $\sigma_5>0$ such that
		\begin{equation}
			\lim_{s\to 0}\|\nabla\phi_t(s)\|\leq C\epsilon^{-\sigma_5}.
		\end{equation}
		Then the solution of \eqref{eq_CHNS} also satisfies the following estimates: for $d = 2, 3$,
		\begin{flalign}
			&\operatorname*{ess\,sup}_{t \in [0, T]} \|\nabla\phi_t\|^2 + \epsilon \int_0^\infty \|\nabla \varDelta \phi_t\|^2 \,\mathrm{d}s \le C\Big(\epsilon^{-\max\{2\sigma_1+7,2\sigma_3+4\}} && \nonumber \\
			&\qquad + \epsilon^{-\frac{2}{6-d}\max\{2\sigma_1+5,2\sigma_3+2\}-\max\{2\sigma_1+13/2,2\sigma_2+4,2\sigma_3+7/2\}} + \epsilon^{-2\sigma_5}\Big) \equiv C\rho_0(\epsilon, d), && \\
			&\int_0^\infty \|\phi_{tt}\|_{H^{-1}}^2 \,\mathrm{d}s \le C\epsilon\rho_0(\epsilon, d) \equiv C\rho_1(\epsilon, d), && 
		\end{flalign}
		and
		\begin{flalign}
			\int_{0}^T\|\phi_t(s)+\boldsymbol{u}(s)\cdot\nabla\phi(s)\|^2+\|\varDelta\mu(s)\|^2 ds\leq C\rho_0(\epsilon,d).
		\end{flalign}
		Moreover, when $d=2$, there are also the following additional estimates:
		\begin{flalign}
			\label{eq_nablau_Deltau_boundedness}
			\operatorname*{ess\,sup}_{t \in [0, T]}\|\nabla\boldsymbol{u}(t)\|^2+\int_{0}^{T}\|\varDelta\boldsymbol{u}(t)\|^2dt\leq C\epsilon^{-2\sigma_1},\\
			\label{eq_ut_nablap_boundedness}
			\int_{0}^{T}[\|\boldsymbol{u}_t(s)\|^2+\|\nabla p(s)\|^2]ds\leq C\epsilon^{-2\sigma_1}.
		\end{flalign}
	\end{Proposition}
	As the proof of Proposition \ref{Proposition_initial_condition} is given in considerable detail in \cite{2006_FengXiaobing_Fully_discrete_finite_element_approximations_of_the_Navier_Stokes_Cahn_Hilliard_diffuse_interface_model_for_two_phase_fluid_flows,2007_FengXiaobing_Analysis_of_finite_element_approximations_of_a_phase_field_model_for_two_phase_fluids}, it is omitted here.
	\begin{Remark}
		(a). When $d = 3$, the estimates \eqref{eq_nablau_Deltau_boundedness} and \eqref{eq_ut_nablap_boundedness} hold only locally in time, which is consistent with the case of the Navier–Stokes equations \cite{2001_TemamRoger_Navier_Stokes_equations_Theory_and_numerical_analysis}. Furthermore, according to \cite{2006_FengXiaobing_Fully_discrete_finite_element_approximations_of_the_Navier_Stokes_Cahn_Hilliard_diffuse_interface_model_for_two_phase_fluid_flows}, the difficulty in extending these estimates to all times stems precisely from the nonlinear term $(\boldsymbol{u}\cdot\nabla)\boldsymbol{u}$, rather than from the nonlinear coupling terms $\boldsymbol{u}\nabla\phi$ and $\mu\nabla\phi$.
		
		(b). Unlike the Navier–Stokes-Allen–Cahn model for two-phase fluids \cite{2007_FengXiaobing_Analysis_of_finite_element_approximations_of_a_phase_field_model_for_two_phase_fluids}, the higher-order norm estimates for solutions of the CHNS model depend on $\frac{1}{\epsilon}$ only polynomially, rather than exponentially \cite{2006_FengXiaobing_Fully_discrete_finite_element_approximations_of_the_Navier_Stokes_Cahn_Hilliard_diffuse_interface_model_for_two_phase_fluid_flows}. This important fact provides the possibility of deriving \emph{a priori} error estimates for numerical solutions of the CHNS model, such that these error estimates also depend polynomially on $\frac{1}{\epsilon}$. 
		For the Cahn–Hilliard equation, $\frac{1}{\epsilon}$-polynomial error estimates have been successfully established in the literature \cite{2004_FengXiaobing_Analysis_of_a_fully_discrete_finite_element_method_for_the_phase_field_model_and_approximation_of_its_sharp_interface_limits,2004_FengXiaobing_Error_analysis_of_a_mixed_finite_element_method_for_the_Cahn_Hilliard_equation,2023_ChenGang_On_the_superconvergence_of_a_hybridizable_discontinuous_Galerkin_method_for_the_Cahn_Hilliard_equation}. In the case of the CHNS coupled system, although \cite{2006_FengXiaobing_Fully_discrete_finite_element_approximations_of_the_Navier_Stokes_Cahn_Hilliard_diffuse_interface_model_for_two_phase_fluid_flows} has touched upon the idea of polynomial-dependent error analysis, the specific error estimation procedure provided therein does not achieve the desired polynomial dependence on $\frac{1}{\epsilon}$. Hence, the $\frac{1}{\epsilon}$-polynomial error analysis for the CHNS system remains unresolved, and to the best of our knowledge, no related results have been reported in the literature. 
	\end{Remark}
	
	Next, we recall and prove some useful inequalities necessary for the error analysis, introduce the quasi-projection, and study its approximation properties.  Throughout, the positive constant $C$ is independent of mesh size $h$, the time step $\tau$, interface width $\epsilon$ and can change from line to line.
	First, for all $v_h\in \mathring{S}_h^r$ and $\xi_h\in \mathring{S}_h^r$, we define the discrete Laplacian operator $-\varDelta_h :\mathring{S}_h^r\to \mathring{S}_h^r$ by 
	\begin{equation}
		\label{eq_discrete_Laplacian_operator}
		\left(-\varDelta_hv_h,\xi_h\right)=\left(\nabla v_h,\nabla \xi_h\right).
	\end{equation}
	Recall that for all $v\in H^1(\Omega)\cap L_0^2(\Omega)$ and $w\in H^1(\Omega)$, the  inverse Laplacian  operator $(-\varDelta)^{-1}:H^{-1}(\Omega)\cap L_0^2(\Omega)\to H^1(\Omega)\cap L_0^2(\Omega)$ is such that 
	\begin{equation}
		\left(\nabla(-\varDelta)^{-1} v, \nabla w\right)=\left(v,w\right), ~\text{condition on}~\frac{\partial(-\varDelta)^{-1} v}{\partial n}=0.
	\end{equation}
	The discrete inverse Laplacian operator $-(\varDelta_h)^{-1}:L_0^2(\Omega)\to \mathring{S}_h^r$ is defined by
	\begin{equation}\label{eq_discrete_Laplacian_operator0002}
		\left(\nabla(-\varDelta_h)^{-1}v,\nabla \xi_h\right)=\left(v,\xi_h\right),~\forall \xi_h\in S_h^r,
	\end{equation}
	and the induced discrete $H^{-1}$-norm is given by
	\begin{equation}\label{eq_Hnev1_definition}
		\|v_h\|_{-1,h}:=\sqrt{(\nabla(-\varDelta_h)^{-1}v_h,\nabla(-\varDelta_h)^{-1}v_h)}=\|\nabla(-\varDelta_h)^{-1}v_h\|=\sup_{{\xi}_h\neq 0}\frac{(v_h,\xi_h)}{\|\nabla\xi_h\|},
	\end{equation}
	which is  symmetric and positive definite on $\mathring{S}_h^r$. 
	\begin{Lemma}[Lemma 2.6 in \cite{2015_Diegel_Analysis_of_a_mixed_finite_element_method_for_a_Cahn_Hilliard_Darcy_Stokes_system}]\label{lemma_H1Hnev1}
		For all \( \chi \in S_h^r \) and all \( \zeta \in \mathring{S}_h^r \),		
		\begin{equation}\label{eq_Hnev1_H1_boundedness00001}
			|(\zeta, \chi)| \leq \|\zeta\|_{-1,h} \|\nabla \chi\|.
		\end{equation}
		The following Poincaré-type estimate holds:		
		\begin{equation}\label{eq_embedding_inequality00001}
			\|\zeta\|_{-1,h} \leq C \|\zeta\|,\quad \forall \ \zeta \in \mathring{S}_h^r,
		\end{equation}
		for some \( C > 0 \) that is independent of \( h \). Finally, if \( \mathcal{T}_h \) is globally quasi-uniform, then the following inverse estimate holds:
		\begin{equation}
			\|\zeta\|  \leq C h^{-1} \|\zeta\|_{-1,h}, \quad \forall \ \zeta \in \mathring{S}_h^r,
		\end{equation}
		for some \( C > 0 \) that is independent of \( h \).
	\end{Lemma}
	If $v_h$ is a constant function, we define $(-\varDelta_h)^{\frac{1}{2}} v_h := 0$ and $(-\varDelta_h)^{-\frac{1}{2}} v_h := 0$. 
	\begin{Lemma}[\cite{2023_CaiWentao_Optimal_L2_error_estimates_of_unconditionally_stable_finite_element_schemes_for_the_Cahn_Hilliard_Navier_Stokes_system}]\label{lemma_laplcian_boundedness}
		For the operators $(-\varDelta_h)^{\frac{1}{2}}:\mathring{S}_h^r\to\mathring{S}_h^r$ and $(-\varDelta_h)^{-\frac{1}{2}}:\mathring{S}_h^r\to\mathring{S}_h^r$, we assume that the $\varepsilon$ is an arbitrary positive constant independent of $h$, for $v_h\in \mathring{S}_h^r$, and  the estimates hold as follows:
		\begin{flalign}
			\label{eq_discrete_Laplacian_bounded00001}
			&C^{-1}\|v_h\|_{H^1}\leq \|(-\varDelta_h)^{\frac{1}{2}} v_h\|\leq C\|v_h\|_{H^1},\\
			\label{eq_discrete_Laplacian_bounded00002}
			&C^{-1}\|v_h\|_{-1,h}\leq \|(-\varDelta_h)^{-\frac{1}{2}} v_h\|\leq C\|v_h\|_{-1,h},
		\end{flalign}
		and
		\begin{equation}
			\label{eq_discrete_Laplacian_bounded00003}
			\|v_h\|\leq \left(\varepsilon\|\nabla v_h\|^2+C\varepsilon^{-1}\|(-\varDelta_h)^{-\frac{1}{2}}v_h\|^2\right)^{\frac{1}{2}}.
		\end{equation}
	\end{Lemma}
	By the above Lemma \ref{lemma_laplcian_boundedness}, if $\boldsymbol{v}_h\in \boldsymbol{X}_h^{r+1}$, the inequalities \eqref{eq_discrete_Laplacian_bounded00001}-\eqref{eq_discrete_Laplacian_bounded00003} still hold.
	For $v_h\in S_h^r$, by subtracting its mean value, we obtain from \eqref{eq_discrete_Laplacian_bounded00003} that 
	\begin{equation}
		\label{eq_vh_boundedness000111}
		\begin{aligned}
			\|v_h\|\leq &~\|v_h-\frac{1}{|\Omega|}(v_h,1)\|+\|\frac{1}{|\Omega|}(v_h,1)\|\\
			\leq&~\varepsilon^{1/2}\|\nabla v_h\|+C\varepsilon^{-1/2}\|(-\varDelta_h)^{-\frac{1}{2}}\tilde{v}_h\|+|(v_h,1)|,
		\end{aligned}
	\end{equation}
	where $\tilde{v}_h:= v_h-\frac{1}{|\Omega|}(v_h,1)$. The several fundamental inequalities are presented as follows: 
	\begin{align}
		\label{eq_Poincare_Wirtinger_inequality}
		&\|v\| \leq C \left( \|\nabla v\|^2 + |(v,1)|^2 \right)^{\frac{1}{2}}, \\
		\label{eq_embedding_inequality}
		&\|v\|_{L^\ell} \leq C \|v\|_{H^1},  \\
		\label{eq_inverse_inequality}
		&\|v_h\|_{W^{m,s}} \leq C h^{n-m+\frac{d}{s}-\frac{d}{q}} \|v_h\|_{W^{n,q}}, 
	\end{align}
	for \( v \in H^1(\Omega) \), \( v_h \in S_h^r, X_h^{r+1}/X_h^{1b} \), \( 1 \leq \ell \leq 6 \), \( 0 \leq n \leq m \leq 1 \), and \( 1 \leq q \leq s \leq \infty \).	

	\subsection{Projection operators}\label{subsec_prjection_opreators}
	Let $R_h: H^1(\Omega)\to S_h^r$ be a classic Ritz projection defined by \cite{1973_Wheeler_Mary_Fanett_A_priori_L2_error_estimates_for_Galerkin_approximations_to_parabolic_partial_differential_equations},
	\begin{equation}
		\left(\nabla(\phi-R_h\phi),\nabla w_h\right)=0,
	\end{equation}
	for all $w_h\in S_h^r$ with $\int_{\Omega}(\phi-R_h\phi)dx=0$. For $s\in [2,\infty) $ and $n=1,2,\cdots,N$ and from the finite element theory \cite{2008_Brenner_Susanne_C_The_mathematical_theory_of_finite_element_methods}, it holds that
	\begin{flalign}
		\label{eq_Ritz_projections_boundedness0001}
		&\|\phi-R_h\phi\|_{L^s}+h\|\phi-R_h\phi\|_{W^{1,s}}\leq Ch^{r+1}\|\phi\|_{W^{r+1,s}},\\		
		\label{eq_Ritz_projections_boundedness0003}
		&\|d_{\tau}(\phi^n-R_h\phi^n)\|+h\|d_{\tau}(\phi^n-R_h\phi^n)\|_{H^1}\leq Ch^{r+1}\|d_{\tau}\phi^n\|_{H^{r+1}}.
	\end{flalign}
	According to Lemma 2 in \cite{1972_ToddDupont_Some_L2_Error_Estimates_for_Parabolic_Galerkin_Methods}, we have used the following approximation properties in $H^{-1}$ of the Ritz projection:
	\begin{flalign}
		\label{eq_Ritz_projections_boundedness0002}
		&\|\phi-R_h\phi\|_{-1,h}\leq C\mathcal{E}_h\|\phi\|_{H^{r+1}},\\
		\label{eq_Ritz_projections_boundedness0004}
		&\|d_{\tau}(\phi^n-R_h\phi^n)\|_{-1,h}\leq C\mathcal{E}_h\|\phi_{t}\|_{H^{r+1}},
	\end{flalign}
	where $\mathcal{E}_h$ is defined by
	\begin{equation}
		\mathcal{E}_h = 
		\begin{cases} 
			h^{r+2}, & \text{for } r \geq 2, \\ 
			h^{2}, & \text{for } r = 1.
		\end{cases}
	\end{equation}
	Now we recall the standard $L^2$ projections $\Pi_h: L^2(\Omega)\to S_h^r$ and $\boldsymbol{I}_h:\boldsymbol{L}^2(\Omega)\to \boldsymbol{X}_h^{1b}/\boldsymbol{X}_h^{r+1}$ by
	\begin{flalign}
		&\left(v-\Pi_hv,w_h\right)=0,\\
		&\left(\boldsymbol{v}-\boldsymbol{I}_h\boldsymbol{v},\boldsymbol{w}_h\right)=0,
	\end{flalign} 
	and the approximation results are classical as  follows \cite{2008_Brenner_Susanne_C_The_mathematical_theory_of_finite_element_methods}:
	\begin{flalign}
		\label{eq_L2_projections_boundedness0001}
		&\|v-\Pi_hv\|+h\|\nabla(v-\Pi_hv)\|\leq Ch^{r+1}\|v\|_{H^{r+1}},\\
		\label{eq_L2_projections_boundedness0002}
		&\|\boldsymbol{v}-\boldsymbol{I}_h\boldsymbol{v}\|+h\|\nabla(\boldsymbol{v}-\boldsymbol{I}_h\boldsymbol{v})\|\leq Ch^{2}\|\boldsymbol{v}\|_{H^{2}}, ~\text{if} ~\boldsymbol{I}_h\boldsymbol{v}\in \boldsymbol{X}_h^{1b},\\
		\label{eq_L2_projections_boundedness0003}
		&\|\boldsymbol{v}-\boldsymbol{I}_h\boldsymbol{v}\|+h\|\nabla(\boldsymbol{v}-\boldsymbol{I}_h\boldsymbol{v})\|\leq Ch^{r+2}\|\boldsymbol{v}\|_{H^{r+2}}, ~\text{if} ~\boldsymbol{I}_h\boldsymbol{v}\in \boldsymbol{X}_h^{r+1}.
	\end{flalign}
	According to \cite{2004_FengXiaobing_Analysis_of_a_fully_discrete_finite_element_method_for_the_phase_field_model_and_approximation_of_its_sharp_interface_limits,2016_FengXiaobing_Analysis_of_mixed_interior_penalty_discontinuous_Galerkin_methods_for_the_Cahn_Hilliard_equation_and_the_Hele_Shaw_flow,2023_ChenGang_On_the_superconvergence_of_a_hybridizable_discontinuous_Galerkin_method_for_the_Cahn_Hilliard_equation}, the $L^{\infty}$-norm estimates for $\phi-R_h\phi$ and $f'(R_h\phi)$ have the following result:
	\begin{equation}
		\label{eq_phiRh_phi_Linfty_bounded}
		\begin{aligned}
			\|\phi-R_h\phi\|_{L^{\infty}}\leq&~ \|\phi-\Pi_h\phi\|_{L^{\infty}}+\|\Pi_h\phi-R_h\phi\|_{L^{\infty}}\\
			\leq&~ Ch^{-\frac{d}{2}}\|\Pi_h\phi-R_h\phi\|+Ch^{2-\frac{d}{2}}\|\phi\|_{H^2}
			\leq Ch^{2-\frac{d}{2}}\|\phi\|_{H^2},
		\end{aligned}
	\end{equation}
	and from the estimate for the continuous solution, \(\|\phi\|_{H^2}\le C\epsilon^{-\alpha'_3}\) with \(\alpha'_3:=\max\{\sigma_1+\frac{5}{2},\sigma_3+1\}\) (see \eqref{eq_Delta_phi_boundedness_continuous}), we obtain the condition \(h^{2-\frac{d}{2}}\le C\epsilon^{\alpha'_3}\), i.e., according to the continuous estimate, \(h\le C\epsilon^{2\alpha'_3/(4-d)}=\epsilon^{\max\{2\sigma_1+5,\,2\sigma_3+2\}/(4-d)}\).
	\begin{equation}
		\label{eq_fphiRh_phi_Linfty_bounded}
		\begin{aligned}
			\|f'(R_h\phi)\|_{L^{\infty}}\leq&~ \|f'(\phi)\|_{L^{\infty}}+\|f'(\phi)-f'(R_h\phi)\|_{L^{\infty}}\\
			\leq &~C\|f'(\phi)\|_{L^{\infty}}+C\|\phi-R_h\phi\|_{L^{\infty}}\leq C,
		\end{aligned}
	\end{equation}
	provided $h\leq \epsilon^{\max\{2\sigma_1+5,\,2\sigma_3+2\}/(4-d)}$. From (2.9) in \cite{2005_FengXiaobing_Numerical_analysis_of_the_Cahn_Hilliard_equation_and_approximation_of_the_Hele_Shaw_problem}, one has
	\begin{equation}\label{eq_boundedness_Linfty_phi_Rhphi}
		\|\phi\|_{L^{\infty}}\leq C, ~\|R_h\phi\|_{L^{\infty}}\leq C.
	\end{equation}
	\subsection{Quasi-projection operators}
	\label{section_projections}
	In this subsection, we will focus on introducing the Ritz quasi-projection operator \cite{2023_CaiWentao_Optimal_L2_error_estimates_of_unconditionally_stable_finite_element_schemes_for_the_Cahn_Hilliard_Navier_Stokes_system} and Stokes quasi-projection operator and its optimal approximation properties.
	Based on the classic Ritz projection, we define a Ritz quasi-projection $\mathcal{G}_h:H^1(\Omega)\to S_h^r$ by
	\begin{equation}\label{eq_Ritz_quasi_projection}
		\left(\nabla(\mu-\mathcal{G}_h\mu),\nabla w_h\right)+\left(\boldsymbol{u}\cdot\nabla(\phi-R_h\phi),w_h\right)=0,
	\end{equation}
	for all $w_h\in S_h^r$ with $\int_{\Omega}(\mu-\mathcal{G}_h\mu)dx=0$. The corresponding estimates are presented in the
	following lemma.
	\begin{Lemma}[Lemma 3.1 in \cite{2023_CaiWentao_Optimal_L2_error_estimates_of_unconditionally_stable_finite_element_schemes_for_the_Cahn_Hilliard_Navier_Stokes_system}]
		For the Ritz quasi-projection \eqref{eq_Ritz_quasi_projection}, it holds that
		\begin{flalign}
			&\|\mu-\mathcal{G}_h\mu\|+h\|\nabla(\mu-\mathcal{G}_h\mu)\|\leq Ch^{r+1}\left(\|\boldsymbol{u}\|_{L^{\infty}}\|\phi\|_{H^{r+1}}+\|\mu\|_{H^{r+1}}\right),\\
			\label{eq_Hnev_Gh}
			&\|\mu-\mathcal{G}_h\mu\|_{-1,h}\leq C\mathcal{E}_h(\|\phi\|_{H^{r+1}}\|\boldsymbol{u}\|_{W^{1,4}}+\|\mu\|_{H^{r+1}}).
		\end{flalign}
	\end{Lemma}
	For the classical Stokes projection \cite{1986_Girault_Vivette_Finite_element_methods_for_Navier_Stokes_equations} and Stokes quasi-projection in \cite{2023_CaiWentao_Optimal_L2_error_estimates_of_unconditionally_stable_finite_element_schemes_for_the_Cahn_Hilliard_Navier_Stokes_system}, we introduce the  Stokes quasi-projection $\left(\boldsymbol{P}_h,Q_h\right):\boldsymbol{H}_0^1(\Omega)\times L_0^2(\Omega)\to \boldsymbol{X}_h^{1b}/\boldsymbol{X}_h^{r+1}\times\mathring{S}_h^r$ defined by
	\begin{flalign}
		\label{eq_Stokes_quasiprojection_definitions0001}
		\left(\nabla(\boldsymbol{u}-\boldsymbol{P}_h\boldsymbol{u}),\nabla\boldsymbol{v}_h\right)-\left(p-Q_hp,\nabla\cdot\boldsymbol{v}_h\right)&=\left((\mu-\mathcal{G}_h\mu)\nabla R_h\phi+\mu\nabla(\phi-R_h\phi),\boldsymbol{v}_h\right),\\
		\label{eq_Stokes_quasiprojection_definitions0002}
		\left(\nabla\cdot(\boldsymbol{u}-\boldsymbol{P}_h\boldsymbol{u}),q_h\right)&=0.
	\end{flalign}
	Its approximation properties are as follows:
	\begin{Lemma}\label{lemma_stokes_projection001}
		For the Stokes quasi-projection defined in \eqref{eq_Stokes_quasiprojection_definitions0001}-\eqref{eq_Stokes_quasiprojection_definitions0002} and $\boldsymbol{P}_h\boldsymbol{u}\in \boldsymbol{X}_h^{r+1}$, it holds that 
		\begin{equation}
			\label{eq_u_boundedness_error_estimate001}
			\|\boldsymbol{u }- \boldsymbol{P}_h \boldsymbol{u}\| \leq C\mathcal{E}_h\bigl(\|\boldsymbol{u}\|_{H^{r+2}} + \|p\|_{H^{r+1}} + \|\phi\|_{H^{r+1}} \|\mu\|_{H^2}+ \|\phi\|_{H^{2}} \|\mu\|_{H^{r+1}}\bigr),
		\end{equation}
		and
		\begin{equation}
			\label{eq_u_boundedness_error_estimate002}
			\begin{aligned}
				\|\nabla(\boldsymbol{u }- \boldsymbol{P}_h \boldsymbol{u})\|  &+ \|p - Q_h p\| \\
				\leq &~	Ch^{r+1}\bigl(\|\boldsymbol{u}\|_{H^{r+2}} + \|p\|_{H^{r+1}} + \|\phi\|_{H^{r+1}} \|\mu\|_{W^{1,4}}+ \|\phi\|_{W^{1,4}} \|\mu\|_{H^{r+1}}\bigr),
			\end{aligned}
		\end{equation}
		and for $n=1,2,\cdots,N$,
		\begin{equation}
			\label{eq_deltau_boundedness_error_estimate001}
			\begin{aligned}
				&\|d_\tau (\boldsymbol{u}^n - \boldsymbol{P}_h \boldsymbol{u}^n)\|\\
				&~\leq 	C\mathcal{E}_h\left(\|d_{\tau}\boldsymbol{u}^n\|_{H^{r+2}}+\|d_{\tau}p^n\|_{H^{r+1}}+\|d_{\tau}\mu^n\|_{H^{r+1}}\|\phi^n\|_{H^{2}}+\|\mu^{n-1}\|_{H^{r+1}}\|d_{\tau}\phi^n\|_{H^2}\right.\\
				&\qquad\left.+\|d_{\tau}\phi^n\|_{H^{r+1}}\|\mu^n\|_{H^2}+\|\phi^{n-1}\|_{H^{r+1}}\|d_{\tau}\mu^n\|_{H^2}\right), \\ 
			\end{aligned}
		\end{equation}
	\end{Lemma}

	\begin{Remark}
		For $\boldsymbol{P}_h\boldsymbol{u}\in \boldsymbol{X}_h^{1b}$, the following approximation holds:
		\begin{equation}
			\label{eq_u_boundedness_error_estimate001_MINI}
			\|\boldsymbol{u }- \boldsymbol{P}_h \boldsymbol{u}\| \leq Ch^{2}\bigl(\|\boldsymbol{u}\|_{H^{2}} + \|p\|_{H^{r+1}} + \|\phi\|_{H^{r+1}} \|\mu\|_{H^2}+ \|\phi\|_{H^{2}} \|\mu\|_{H^{r+1}}\bigr),
		\end{equation}
		and
		\begin{equation}
			\label{eq_u_boundedness_error_estimate002_MINI}
			\begin{aligned}
				\|\nabla(\boldsymbol{u }- \boldsymbol{P}_h \boldsymbol{u})\|  &+ \|p - Q_h p\| \\
				\leq &~	Ch \bigl(\|\boldsymbol{u}\|_{H^{2}} + \|p\|_{H^{2}} + \|\phi\|_{H^{r+1}} \|\mu\|_{W^{1,4}}+ \|\phi\|_{W^{1,4}} \|\mu\|_{H^{r+1}}\bigr),
			\end{aligned}
		\end{equation}
		and for $n=1,2,\cdots,N$,
		\begin{equation}
			\label{eq_deltau_boundedness_error_estimate001_MINI}
			\begin{aligned}
				&\|d_\tau (\boldsymbol{u}^n - \boldsymbol{P}_h \boldsymbol{u}^n)\|\\
				&~\leq 	Ch^{2}\left(\|d_{\tau}\boldsymbol{u}^n\|_{H^{2}}+\|d_{\tau}p^n\|_{H^{r+1}}+\|d_{\tau}\mu^n\|_{H^{r+1}}\|\phi^n\|_{H^{2}}+\|\mu^{n-1}\|_{H^{r+1}}\|d_{\tau}\phi^n\|_{H^2}\right.\\
				&\qquad\left.+\|d_{\tau}\phi^n\|_{H^{r+1}}\|\mu^n\|_{H^2}+\|\phi^{n-1}\|_{H^{r+1}}\|d_{\tau}\mu^n\|_{H^2}\right).
			\end{aligned}
		\end{equation}
		In this case, the optimal $L^2$ error order for the velocity is $2$, and the analytical framework and methodology adopted in this paper remain applicable. Therefore,  the MINI finite element subspace yields an optimal $L^2$-error estimate for the velocity \cite{2023_CaiWentao_Optimal_L2_error_estimates_of_unconditionally_stable_finite_element_schemes_for_the_Cahn_Hilliard_Navier_Stokes_system,2026_GaoHaijun_Optimal_L2_error_estimates_of_a_decoupled_FEM_for_the_Cahn_Hilliard_Navier_Stokes_equations}. Next, the subsequent discussion will focus mainly on finite element spaces $S_h^r\times S_h^r\times\boldsymbol{X}_h^{r+1}\times \mathring{S}_h^r$. 
	\end{Remark}
	From \eqref{eq_u_boundedness_error_estimate001} and \eqref{eq_u_boundedness_error_estimate002}, we can obtain the boundedness
	\begin{equation}
		\label{eq_Phu_boundedness00001}
		\|\boldsymbol{P}_h\boldsymbol{u}\|_{L^{\infty}}+\|\boldsymbol{P}_h\boldsymbol{u}\|_{W^{1,3}}\leq C\left(\|\boldsymbol{u}\|_{H^2}+\|p\|_{H^2}+\|\phi\|_{H^2}\|\mu\|_{H^2}\right).
	\end{equation}
	The proof of the approximation property of the Stokes quasi-projection presented above is similar to that of Lemma 4 in \cite{2026_GaoHaijun_Optimal_L2_error_estimates_of_a_decoupled_FEM_for_the_Cahn_Hilliard_Navier_Stokes_equations}, and is therefore omitted here for brevity.
	\begin{Remark}
		The presence of the term $(\mu-\mathcal{G}_h\mu)\nabla\phi$ on the right-hand side of \eqref{eq_Stokes_quasiprojection_definitions0001} can be explained as follows. In the error estimation procedure, terms of the form $\left((\mu-\mathcal{G}_h\mu)\nabla R_h\phi, e_{\boldsymbol{u}}\right)$ naturally arise. If the term $(\mu-\mathcal{G}_h\mu)\nabla R_h\phi$ is not properly eliminated, the undesired terms appear in the resulting estimates as follows:
		\begin{equation}
			\begin{aligned}
				\left((\mu-\mathcal{G}_h\mu)\nabla R_h\phi, e_{\boldsymbol{u}}\right)\leq &~C\|\mu-\mathcal{G}_h\mu\|_{-1,h}\|\nabla R_h\phi e_{\boldsymbol{u}}\|_{H^1}\\
				\leq &~C\|\mu-\mathcal{G}_h\mu\|_{-1,h}\left(\|\nabla R_h\phi\|_{L^3}\|e_{\boldsymbol{u}}\|_{L^6}+\|e_{\boldsymbol{u}}\|_{L^{\infty}}\|R_h\phi\|_{H^2}\right),
			\end{aligned}
		\end{equation}
		which prevents us from obtaining optimal convergence results. In particular, we are unable to obtain an effective control of the $L^{\infty}$-norm of $e_{\boldsymbol{u}}$ and $\nabla R_h\phi$, for instance, by bounding it in terms of the $H^1$-norm of $e_{\boldsymbol{u}}$. Similar difficulties arise in several coupling terms $\boldsymbol{u}\cdot\nabla\phi$ and $\mu\nabla\phi$ throughout the analysis. Consequently, we employ the quasi-projection to properly handle these coupling terms that hinder optimal error convergence.
	\end{Remark}
	
	\subsection{The spectral estimate of the Cahn-Hilliard operator}In this subsection, we introduce spectral estimation methods for the discrete Cahn-Hilliard operator. Based on Proposition 1 in \cite{2004_FengXiaobing_Error_analysis_of_a_mixed_finite_element_method_for_the_Cahn_Hilliard_equation} concerning the spectral estimation of the continuous Cahn–Hilliard operator, we present a spectral estimation lemma for the discrete Cahn-Hilliard operator.
	\begin{Lemma}\label{lemma_spectral_estimates_boundedness}
		Assume that the spectral estimate of the continuous Cahn-Hilliard operator holds. Then for any fixed $\beta\in (0,1)$, the spectral estimate of a discrete Cahn-Hilliard operator satisfies for all $\psi_h\in \mathring{S} _h^r$,
		\begin{equation}\label{eq_spectral_estimate_CHNS}
			\epsilon\|\nabla\psi_h\|^2+\frac{1-\beta}{\epsilon}\left(f'(R_h\phi)\psi_h,\psi_h\right)\geq -(2C_0+1)\|\psi_h\|_{-1,h}^2,
		\end{equation}
		where $h\leq C\min\{\beta^{\frac{1}{2}}\epsilon^{\frac{11}{4}},\beta^{\frac{3}{4}}\epsilon^{\frac{7}{4}},\epsilon^{\frac{2\alpha'_3+6}{4-d}}\}$ and $C_0$ is a constant in the spectral estimation of the continuous Cahn–Hilliard operator.
	\end{Lemma}
	For the spectral estimate of the Cahn--Hilliard operator, proofs for its DG and HDG formulations have been provided in \cite{2016_FengXiaobing_Analysis_of_mixed_interior_penalty_discontinuous_Galerkin_methods_for_the_Cahn_Hilliard_equation_and_the_Hele_Shaw_flow,2023_ChenGang_On_the_superconvergence_of_a_hybridizable_discontinuous_Galerkin_method_for_the_Cahn_Hilliard_equation}, respectively. Since the proof for the continuous finite element formulation considered here follows essentially the same arguments as those in \cite{2016_FengXiaobing_Analysis_of_mixed_interior_penalty_discontinuous_Galerkin_methods_for_the_Cahn_Hilliard_equation_and_the_Hele_Shaw_flow,2023_ChenGang_On_the_superconvergence_of_a_hybridizable_discontinuous_Galerkin_method_for_the_Cahn_Hilliard_equation}, we omit the details for brevity.

	We next give the Gr\"{o}nwall inequality and provide some explanation.
	The classical Gr\"{o}nwall lemma provides a fundamental estimate for functions satisfying a first-order linear differential inequality. This lemma serves as a primary technical tool in the error analysis of continuous-in-time semi-discrete schemes for deriving error estimates of numerical solutions to initial-boundary value problems of partial differential equations. Consequently, its discrete counterpart constitutes an indispensable theoretical foundation in the error analysis of fully discrete schemes. However, for a wide class of nonlinear partial differential equations, the classical Gr\"{o}nwall lemma is no longer applicable due to the presence of nonlinear effects, and thus one should resort to its nonlinear generalizations. For nonlinearities of power type, Feng et al. \cite{2005_FengXiaobing_Numerical_analysis_of_the_Cahn_Hilliard_equation_and_approximation_of_the_Hele_Shaw_problem} established a corresponding generalized Gr\"{o}nwall lemma. In the lemma presented below, we provide a discrete version of their result, for which a similar proof can be found in \cite{2001_Feng_Prohl_Numerical_analysis_of_the_Cahn_Hilliard_equation_Part_II}. This discrete lemma will play a crucial role in the theoretical error analysis of the next section.
	\begin{Lemma}[Nonlinear Gr\"{o}nwall's inequality \cite{2015_FengXiaobing_Analysis_of_symmetric_interior_penalty_discontinuous_Galerkin_methods_for_the_Allen_Cahn_equation_and_the_mean_curvature_flow}]
		\label{lemma_nonlinear_gronwall_inequalities}
		Let \(\{S_m\}_{m \ge 1}\) be a positive nondecreasing sequence and \(\{b_m\}_{m \ge 1}\) and \(\{\tau_m\}_{m \ge 1}\) be non-negative sequences, and \(p > 1\) be a constant. If
		\begin{equation}
			S_{m+1} - S_m \le b_m S_m + \tau_m S_m^p, \quad \text{for } m \ge 1, 
		\end{equation}
		and for  $m \ge 2$,
		\begin{equation}
			S_1^{1-p} + (1-p) \sum_{n=1}^{m-1} \tau_n a_{n+1}^{1-p} > 0,
		\end{equation}
		then
		\begin{equation}
			S_m \le \frac{1}{a_m} \left\{ S_1^{1-p} + (1-p) \sum_{n=1}^{m-1} \tau_n a_{n+1}^{1-p} \right\}^{1/(1-p)},
		\end{equation}
		where $a_m$ is defined by \[a_m := \prod_{n=1}^{m-1} \frac{1}{1 + b_n}.\]
	\end{Lemma}
	\section{Numerical scheme and stability}\label{section_numerical_scheme}
	In this section, we introduce the fully discrete scheme of the mixed finite element method and discuss the stability of numerical scheme.
	\subsection{Numerical scheme}
	For a fixed integer $N$, let $\{t_n=n\tau\}_{n=0}^N$ be a uniform partition of the time interval $[0,T]$ with $\tau=T/N$. We denote $g^n=g(x,t_n)$, and for any sequence $\{g^n\}_{n=1}^N$, define 
	$$
	d_{\tau}g^n=\frac{g^n-g^{n-1}}{\tau}.
	$$
	With the above notation, the fully discrete mixed finite element scheme based on the first order convex-splitting approach can be constructed utilizing the backward Euler method:
	find $\left(\phi_h^n,\mu_h^n,\boldsymbol{u}_h^n,p_h^n\right)\in S_h^r\times S_h^r \times \boldsymbol{X}_h^{r+1}\times \mathring{S}_h^r$
	such that
	\begin{flalign}
		\label{eq_CHNS_variational_form_phi_discretization}
		\left(d_{\tau}\phi_h^n,w_h\right)+b\left( \phi_h^{n-1},\boldsymbol{u}_h^n,w_h\right)+\left(\nabla\mu_h^n,\nabla w_h\right)&~=0,\\
		\label{eq_CHNS_variational_form_mu_discretization}
		\left(\mu_h^n,\varphi_h\right)-\epsilon\left(\nabla\phi_h^n,\nabla\varphi_h\right)-\epsilon^{-1}\left(F'(\phi_h^{n},\phi_h^{n-1}),\varphi_h\right)&~=0,\\
		\label{eq_CHNS_variational_form_ns_discretization}
		\left(d_{\tau}\boldsymbol{u}_h^n,\boldsymbol{v}_h\right)+\boldsymbol{B}\left(\boldsymbol{u}_h^{n-1},\boldsymbol{u}_h^n,\boldsymbol{v}_h\right)+\nu\left(\nabla\boldsymbol{u}_h^n,\nabla\boldsymbol{v}_h\right)\notag\\
		-\left( p_h^n,\nabla\cdot\boldsymbol{v}_h\right)
		-b\left(\phi_h^{n-1},\boldsymbol{v}_h,\mu_h^n\right)&~=0,\\
		\label{eq_CHNS_variational_form_incompressible_u_discretization}
		\left(\nabla\cdot\boldsymbol{u}_h^n,q_h\right)&~=0.
	\end{flalign}
	
	Here,  for any test functions $\left(w_h,\varphi_h,\boldsymbol{v}_h,q_h\right)\in S_h^r\times S_h^{r}\times \boldsymbol{X}_h^{r+1}\times \mathring{S}_h^r $ and $n=1,2,\cdots,N$, the nonlinear potential energy function $F'(\phi_h^{n},\phi_h^{n-1})= (\phi_h^{n})^3-\phi_h^{n-1}$. We choose the initial value $\phi_h^0=R_h\phi^0$, 
	$\boldsymbol{u}_h^0=\boldsymbol{I}_h\boldsymbol{u}^0$. Here, $R_h$ and $\boldsymbol{I}_h$  denote the Ritz projection, the $L^2$ projection defined in subsection \ref{subsec_prjection_opreators}, respectively. 
	
	\subsection{Stability}
	In this subsection, we establish the stability of the fully discrete scheme \eqref{eq_CHNS_variational_form_phi_discretization}-\eqref{eq_CHNS_variational_form_incompressible_u_discretization}, including initial value assumptions and the main conclusions of this paper. First, it is easy to establish the unconditional energy stability for the fully discrete scheme \eqref{eq_CHNS_variational_form_phi_discretization}-\eqref{eq_CHNS_variational_form_incompressible_u_discretization}, where the discrete energy is defined as follows:
	\begin{equation}\label{eq_discrete_energy}
		\mathcal{J}_{\epsilon}(\phi_h^n,\boldsymbol{u}_h^n)=\frac{\epsilon}{2}\|\nabla\phi_h^n\|^2+\frac{1}{4\epsilon}\|(\phi_h^n)^2-1\|^2+\frac{1}{2}\|\boldsymbol{u}_h^n\|^2.
	\end{equation}
	According to the unconditional discrete energy stability, we have the following stability estimates for the fully discrete scheme.
	\begin{Lemma}\label{lemma_initial_data_conditional}
		The solution $\left(\phi_h^n,\mu_h^n,\boldsymbol{u}_h^n,p_h^n\right)$ of the fully discrete scheme \eqref{eq_CHNS_variational_form_phi_discretization}-\eqref{eq_CHNS_variational_form_incompressible_u_discretization} for $n=1,2,\cdots,N$ satisfies 
		\begin{flalign}
			\label{eq_mass_conservation}
			&\left(\phi_h^n,1\right)=\left(\phi^0,1\right),\\
			&\max_{0\leq n\leq N}\bigg\{\epsilon\|\nabla\phi_h^n\|^2+\epsilon^{-1}\|F(\phi_h^n)\|_{L^1}+\|\boldsymbol{u}_h^n\|^2\bigg\}\notag\\
			&\qquad\qquad+\tau\sum_{n=1}^{N}\left(\|\nabla\mu_h^n\|^2+\epsilon \tau\|\nabla d_{\tau}\phi_h^n\|^2+\|\nabla\boldsymbol{u}_h^n\|^2\right)\leq C\epsilon^{-2\sigma_1},\\
			&\max_{0\leq n\leq N}\|\nabla\mu_h^n\|^2+\tau\sum_{n=1}^{N}\left(\tau\|d_{\tau}\nabla\mu_h^n\|^2+\epsilon\|\nabla d_{\tau}\phi_h^n\|^2\right)\leq C\epsilon^{-\max\{2\sigma_1+3,2\sigma_3}\}.
		\end{flalign}
	\end{Lemma}
	Since the fully discrete scheme \eqref{eq_CHNS_variational_form_phi_discretization}-\eqref{eq_CHNS_variational_form_incompressible_u_discretization} is an unconditional energy stability and according to Lemma 3.1 in \cite{2005_FengXiaobing_Numerical_analysis_of_the_Cahn_Hilliard_equation_and_approximation_of_the_Hele_Shaw_problem}, the above lemma can be easily derived. Therefore, the detailed proof of this Lemma \ref{lemma_initial_data_conditional} is omitted here.
	We assume that the initial energy $\mathcal{J}_{\epsilon}(\phi_h^0,\boldsymbol{u}_h^0)$ is uniformly bounded in terms of $\epsilon$ as follows:
	\begin{equation}\label{eq_initial_energy_bounded}
		\mathcal{J}_{\epsilon}(\phi_h^0,\boldsymbol{u}_h^0)=\frac{1}{4\epsilon}\|(\phi_h^0)^2-1\|^2+\frac{\epsilon}{2}\|\nabla\phi_h^0\|^2+\frac{1}{2}\|\boldsymbol{u}_h^0\|^2\leq C\epsilon^{-2\sigma_1},
	\end{equation}
	and that $\|\phi_h^0\|\leq C\epsilon^{-(2\sigma_1+1)}$. 
	This assumption is reasonable from the standpoint of the sharp interface limit of the Cahn-Hilliard free energy functional \cite{1987_Elliott_Numerical_studies_of_the_Cahn_Hilliard_equation_for_phase_separation}. 
	By an elementary fixed point argument and the energy method, the unconditional unique solvability of the fully discrete scheme \eqref{eq_CHNS_variational_form_phi_discretization}-\eqref{eq_CHNS_variational_form_incompressible_u_discretization} is established in \cite{2023_CaiWentao_Optimal_L2_error_estimates_of_unconditionally_stable_finite_element_schemes_for_the_Cahn_Hilliard_Navier_Stokes_system,2010_ShenJie_Energy_stable_schemes_for_Cahn_Hilliard_phase_field_model_of_two_phase_incompressible_flows}. The basic energy stability estimates are provided in the following lemma.
	\begin{Lemma}\label{lemma_stability_bounds}
		For any $\tau,h>0$ and $m=1,2,\cdots,N$, the stability bound holds as follows:
		\begin{equation}
			\frac{1}{4\epsilon}\|(\phi_h^m)^2-1\|^2+\frac{\epsilon}{2}\|\nabla\phi_h^m\|^2+\frac{1}{2}\|\boldsymbol{u}_h^m\|^2+\tau\sum_{n=1}^{m}(\|\nabla\mu_h^n\|^2+\|\nabla\boldsymbol{u}_h^n\|^2)\leq C\epsilon^{-2\sigma_1}.
		\end{equation}
	\end{Lemma}
	From Lemmas \ref{lemma_initial_data_conditional} and \ref{lemma_stability_bounds}, we derive the boundedness of a discrete solution, as shown in the following Lemmas \ref{lemma_varDelta_phi_bounded} and \ref{lemma_phi_Linf_boundedness}.
	\begin{Lemma}\label{lemma_varDelta_phi_bounded}
		For all $m=1,2,\cdots,N$, it holds
		\begin{equation}\label{eq_L4_boundedness0001}
			\|\phi_h^m\|_{L^4}^4+\tau\sum_{n=1}^{m}\|\nabla\mu_h^n\|^2\leq C(1+\epsilon^{-(2\sigma_1+1)}),~\|\nabla\phi_h^n\|^2\leq C\epsilon^{-(2\sigma_1+1)},
		\end{equation}
		and for all $n=1,2,\cdots,m\leq N$,
		\begin{equation}\label{eq_varDeltaphimu}
			\|\varDelta_h\phi_h^n\|\leq C\epsilon^{-\max\{3\sigma_1+7/2,\sigma_3+1\}}=C\epsilon^{-\alpha_3},
		\end{equation}
		where $C$ depends on the time $T$ and the initial condition.
	\end{Lemma}
	\begin{proof}
		According to Lemmas \ref{lemma_initial_data_conditional} and \ref{lemma_stability_bounds}, we obtain the result of inequalities \eqref{eq_L4_boundedness0001}. Next, taking $\varphi_h=1$ in \eqref{eq_CHNS_variational_form_mu_discretization} yields 
		\begin{equation}\label{eq_takinginner_one00001}
			|\left(\mu_h^n,1\right)|=\frac{1}{\epsilon}\bigg|\left(f(\phi_h^n,\phi_h^{n-1}),1\right)\bigg|\leq\frac{1}{\epsilon}\left(\|\phi_h^n\|_{L^3}^3+\|\phi_h^{n-1}\|_{L^1}\right)\leq C\epsilon^{-(3\sigma_1+5/2)}.
		\end{equation}
		From \eqref{eq_Poincare_Wirtinger_inequality}, we have
		\begin{equation}\label{eq_mu_bounded00001}
			\|\mu_h^n\|^2\leq C\left(\|\nabla\mu_h^n\|^2+|(\mu_h^n,1)|^2\right)\leq C\epsilon^{-\max\{6\sigma_1+5,2\sigma_3\}}.
		\end{equation}
		Taking $\varphi_h=\varDelta_h\phi_h^n$ in \eqref{eq_CHNS_variational_form_mu_discretization}, we have
		\begin{equation}
			\begin{aligned}
				\epsilon\|\varDelta_h\phi_h^n\|^2=&~\epsilon^{-1}\left(f(\phi_h^n,\phi_h^{n-1}),\varDelta_h\phi_h^n\right)-\left(\mu_h^n,\varDelta_h\phi_h^n\right).
			\end{aligned}
		\end{equation}
		Using the Cauchy-Schwarz inequality, we get
		\begin{equation}
			\epsilon\|\varDelta_h\phi_h^n\|^2
			\leq \epsilon^{-1}\left(\|\phi_h^n\|_{L^6}^3+\|\phi_h^{n-1}\|\right)\|\varDelta_h\phi_h^n\|+\|\mu_h^n\|\|\varDelta_h\phi_h^n\|
		\end{equation}
		Next we use the Sobolev inequality \eqref{eq_embedding_inequality} and the Lemma \ref{lemma_stability_bounds}, and obtain
		\begin{equation}
			\|\phi_h^n\|_{L^6}\leq C\left(\|\phi_h^n\|_{L^4}+\|\nabla\phi_h^n\|\right)\leq C\epsilon^{-(\sigma_1+1/2)}.
		\end{equation}
		Using Young's inequality and \eqref{eq_L4_boundedness0001}, we obtain the first desired result \eqref{eq_varDeltaphimu},
		which completes the proof of this lemma.
	\end{proof}
	Using the Lemma \ref{lemma_varDelta_phi_bounded} and applying the discrete Gagliardo-Nirenberg inequality 
	\begin{equation}\label{eq_Gagliardo_Nirenberg_inequality}
		\|\phi_h\|_{L^{\infty}}\leq C\left(\|\varDelta_h\phi_h\|^{\frac{d}{2(6-d)}}\|\phi_h\|_{L^6}^{\frac{3(4-d)}{2(6-d)}}+\|\phi_h\|_{L^6}\right),
	\end{equation}
	we deduce the $L^{\infty}$-norm estimate of $\phi_h^n$ as follows:
	\begin{Lemma}
		\label{lemma_phi_Linf_boundedness}
		For all $n=1,2,\cdots,N$, it holds 
		\begin{equation}\label{eq_phi_Linf_boundedness000121203}
			\tau\sum_{n=1}^{N}\left(\|\phi_h^n\|_{L^{\infty}}^2+\|\nabla\phi_h^n\|_{L^6}^2\right)\leq C\epsilon^{-\max\{6\sigma_1+7,2\sigma_3+2\}}.
		\end{equation}
	\end{Lemma}
	\begin{Remark}
		If no amplification treatment is applied to \eqref{eq_phi_Linf_boundedness000121203}, the exponent given by the Gagliardo-Nirenberg inequality is 
		\begin{equation}
			\|\phi_h^n\|_{L^{\infty}}\leq C\left(\epsilon^{-\frac{d}{2(6-d)}\max\{3\sigma_1+7/2,\sigma_3+1\}}(1+\epsilon^{\frac{-2\sigma_1+1}{4}})^{\frac{3(4-d)}{2(6-d)}}+\epsilon^{-(\sigma_1+1/2)}\right).
		\end{equation}
		In the cases \(d=2,3\), these are approximately \(\epsilon^{-\frac{1}{4}\alpha_3-\frac{3}{8}\sigma_1-\frac{3}{16}}\) (for \(d=2\)) and \(\epsilon^{-\frac{1}{2}\alpha_3-\frac{1}{4}\sigma_1-\frac{1}{8}}\) (for \(d=3\)), respectively, which are slightly better than \(\epsilon^{-\alpha_3}\) and $\alpha_3:=\max\{3\sigma_1+7/2,\sigma_3+1\}$. To maintain the consistency of scaling with Lemma \ref{lemma_phi_Linf_boundedness}, the concise form \eqref{eq_phi_Linf_boundedness000121203} is therefore adopted in this paper.
	\end{Remark}

	\begin{Lemma}\label{lemma_L4_bounded000011}
		For any $\tau, h>0$ and $m=1,2,\cdots,N$, it holds
		\begin{equation}
			\|\mu_h^m\|^2+\sum_{n=1}^{m}\|\mu_h^n-\mu_h^{n-1}\|^2+\epsilon\tau\sum_{n=1}^{m}\|d_{\tau}\phi_h^n\|^2\leq C\epsilon^{-\alpha_4},
		\end{equation}
		where we defined by $\alpha_4:=\max\{6\sigma_1+8,4\sigma_1+2\sigma_3+5,8\sigma_1+6,2\sigma_2\}$.
	\end{Lemma}
	\begin{proof}
		From \eqref{eq_Poincare_Wirtinger_inequality} and \eqref{eq_takinginner_one00001}, we have
		\begin{equation}\label{eq_tau_mu_bounded00001}
			\tau\sum_{n=1}^{N}\|\mu_h^n\|^2\leq C\tau\sum_{n=1}^{N}\left(\|\nabla\mu_h^n\|^2+|(\mu_h^n,1)|^2\right)\leq C\epsilon^{-\max\{6\sigma_1+5,2\sigma_3\}}.
		\end{equation}
		Applying $d_{\tau}$ to \eqref{eq_CHNS_variational_form_mu_discretization}, we give
		\begin{equation}
			\left(d_{\tau}\mu_h^n,\varphi_h\right)=\epsilon\left(d_{\tau}\nabla\phi_h^n,\nabla\varphi_h\right)+\epsilon^{-1}\left(d_{\tau}f(\phi_h^{n},\phi_h^{n-1}),\varphi_h\right).
		\end{equation}
		Taking $\left(w_h,\varphi_h\right)=(d_{\tau}\phi_h^n,\mu_h^n)$ in \eqref{eq_CHNS_variational_form_phi_discretization} and \eqref{eq_CHNS_variational_form_mu_discretization}, and combining the two equations, one obtains
		\begin{equation}
			\left(d_{\tau}\mu_h^n,\mu_h^n\right)+\epsilon\|d_{\tau}\phi_h^n\|^2=\epsilon^{-1}\left(d_{\tau}f(\phi_h^n,\phi_h^{n-1}),\mu_h^n\right)-\epsilon b\left(\phi_h^{n-1},\boldsymbol{u}_h^n, d_{\tau}\phi_h^n\right),
		\end{equation}
		and we have the following result 
		\begin{equation}
			\label{eq_combining_results_00001}
			\begin{aligned}
				\|\mu_h^m\|^2+&\sum_{n=1}^{m}\|\mu_h^n-\mu_h^{n-1}\|^2+2\tau\epsilon\sum_{n=1}^{m}\|d_{\tau}\phi_h^n\|^2\\
				=&~2\tau\sum_{n=1}^{m}\left(\epsilon^{-1}\left(d_{\tau}f(\phi_h^n,\phi_h^{n-1}),\mu_h^n\right)-\epsilon\left(\nabla\phi_h^{n-1}\cdot\boldsymbol{u}_h^n, d_{\tau}\phi_h^n\right)\right)+\|\mu_h^0\|^2,
			\end{aligned}
		\end{equation}
		where $\mu_h^0$	defined by \eqref{eq_CHNS_variational_form_mu_discretization} with $n = 0$ and $\phi_h^{-1}:=\phi_h^0$.
		We will now prove the estimate on the right-hand side of \eqref{eq_combining_results_00001}, obtaining
		\begin{equation}
			\begin{aligned}
				\frac{\tau}{\epsilon}\sum_{n=1}^{m}\left(d_{\tau}f(\phi_h^n,\phi_h^{n-1}),\mu_h^n\right)\leq& \frac{C\tau}{\epsilon}\sum_{n=1}^{m}\|d_{\tau}\phi_h^n\|\|\mu_h^n\|_{L^6}\left(\|\phi_h^n\|_{L^6}^2+\|\phi_h^{n-1}\|_{L^6}^2+1\right)\\
				\leq&\frac{\epsilon\tau}{4}\sum_{n=1}^{m}\|d_{\tau}\phi_h^n\|^2+C\tau\epsilon^{-(4\sigma_1+5)}\sum_{n=1}^{m}\|\nabla\mu_h^n\|^2,
			\end{aligned}
		\end{equation}
		and
		\begin{equation}
			\begin{aligned}
				\epsilon\tau\sum_{n=1}^{m}\left(\nabla\phi_h^{n-1}\cdot\boldsymbol{u}_h^n, d_{\tau}\phi_h^n\right)\leq& C\epsilon\tau\sum_{n=1}^{m}\|\nabla\phi_h^{n-1}\|_{L^6}\|\boldsymbol{u}_h^n\|_{L^{3}}\|d_{\tau}\phi_h^n\|\\
				\leq&\frac{\epsilon\tau}{4}\sum_{n=1}^{m}\|d_{\tau}\phi_h^n\|^2+C\tau\epsilon\sum_{n=1}^{m}\|\varDelta_h\phi_h^{n-1}\|^2\|\nabla\boldsymbol{u}_h^n\|^2\\
				\leq &~\frac{\epsilon\tau}{4}\sum_{n=1}^{m}\|d_{\tau}\phi_h^n\|^2+C\tau\epsilon^{-\max\{8\sigma_1+6,2\sigma_1+2\sigma_3+1\}}.
			\end{aligned}
		\end{equation}
		Using the Sobolev embedding inequality, we have 
		\begin{equation}
			\begin{aligned}
				\|\mu_h^m\|^2+&\sum_{n=1}^{m}\|\mu_h^n-\mu_h^{n-1}\|^2+\tau\epsilon\sum_{n=1}^{m}\|d_{\tau}\phi_h^n\|^2\\
				\leq &~C\tau\sum_{n=1}^{m}\left(\epsilon^{-(4\sigma_1+5+\max\{2\sigma_1+3,2\sigma_3\})}+\epsilon^{-\max\{8\sigma_1+6,2\sigma_1+2\sigma_3+1\}}\right)+\|\mu_h^0\|^2\\
				\leq &C\epsilon^{-\max\{6\sigma_1+8,4\sigma_1+2\sigma_3+5,8\sigma_1+6,2\sigma_2\}},
			\end{aligned}
		\end{equation}
		where we use the estimate \eqref{eq_tau_mu_bounded00001} and the bounds in Lemma \ref{lemma_stability_bounds}, thereby completing the proof of the lemma.
	\end{proof}
	\begin{Lemma}\label{lemma_dtau_phi_Linf_bounded}
		For any $\tau, h>0$ and $m=1,2,\cdots,N$, it holds
		\begin{equation}
			\tau\sum_{n=1}^m\|\varDelta_h\mu_h^n\|\leq C\epsilon^{-\max\{4\sigma_1+7/2,\sigma_1+\sigma_3+1\}},
		\end{equation}
		and 
		\begin{equation}
			\tau\sum_{n=1}^m\|d_{\tau}\phi_h^{n}\|_{-1,h}\leq C\epsilon^{-\max\{\sigma_1+3/2,\sigma_3,4\sigma_1+7/2,\sigma_1+\sigma_3+1\}}.
		\end{equation}
	\end{Lemma}
	\begin{proof}   
		Taking $w_h=\varDelta_h\mu_h^n$ in \eqref{eq_CHNS_variational_form_phi_discretization} yields 
		\begin{equation}
			\begin{aligned}
				\|\varDelta_h\mu_h^n\|^2=&~\left(d_{\tau}\phi_h^n,\varDelta_h\mu_h^n\right)+\left(\boldsymbol{u}_h^n\nabla\phi_h^{n-1},\varDelta_h\mu_h^n\right)\\
				\leq &~\| d_{\tau}\phi_h^n\|\|\varDelta\mu_h^n\|+\|\nabla\boldsymbol{u}_h^n\|\|\varDelta_h\phi_h^{n-1}\|\|\varDelta_h\mu_h^n\|.
			\end{aligned}
		\end{equation}
		Thus, according to the above inequality, it is easy to see that 
		\begin{equation}\label{eq_Linf_muh}
			\tau\sum_{n=1}^m\|\varDelta_h\mu_h^n\|\leq  C\epsilon^{-\max\{4\sigma_1+7/2,\sigma_1+\sigma_3+1\}}.
		\end{equation}
		Then, we take $w_h=(-\varDelta_h)^{-1}d_{\tau}\phi_h^n$ in \eqref{eq_CHNS_variational_form_phi_discretization} to yield
		\begin{equation}
			\begin{aligned}
				\|d_{\tau}\phi_h^n\|_{-1,h}^2=&~\left(\mu_h^n,d_{\tau}\phi_h^n\right)+(\phi_h^{n-1}\boldsymbol{u}_h^n,\nabla(-\varDelta_h)^{-1}d_{\tau}\phi_h^n)\\
				\leq &~\|\nabla\mu_h^n\|\|d_{\tau}\phi_h^n\|_{-1,h}+C\|\phi_h^{n-1}\|_{L^{\infty}}\|\boldsymbol{u}_h^n\|\|d_{\tau}\phi_h^n\|_{-1,h}.
			\end{aligned}
		\end{equation}
		Thus, it is easy to arrive at the conclusion in Lemma \ref{lemma_dtau_phi_Linf_bounded}.
	\end{proof}
	
	
	\section{Error estimates}\label{section_error_estimates}
	In this section, we establish the optimal $L^2$ error estimates of the fully discrete
	scheme \eqref{eq_CHNS_variational_form_phi_discretization}-\eqref{eq_CHNS_variational_form_incompressible_u_discretization}  that depends only on polynomials $\epsilon^{-1}$ to solve the CHNS equations. We first perform the error estimate for the phase-field variable in the negative norm and then under the \( L^2 \) norm, thus establishing the main convergence result. Throughout the analysis, we assume that the initial data and the solution of equation \eqref{eq_CHNS} are smooth enough. In this paper, we assume that the solution to the CHNS model \eqref{eq_CHNS}-\eqref{eq_CHNS_initial_boundary_conditions} exists and satisfies the following regularities (for $r \geq 1$):
	\begin{equation}\label{eq_regularity_condition}
		\begin{aligned}
			\phi &\in H^2(0,T;L^2(\Omega))
			\cap H^1(0,T;H^{r+1}(\Omega))
			\cap C([0,T];W^{2,4}(\Omega)),
			\mu \in H^1(0,T;H^{r+1}(\Omega)),\\
			\boldsymbol{u} &\in H^2(0,T;L^2(\Omega))
			\cap H^1(0,T;\boldsymbol{H}^{r+2}(\Omega)),
			p \in L^2(0,T;H^{r+1}(\Omega)\cap L_0^2(\Omega)).
		\end{aligned}
	\end{equation}
	We present our main result in the following theorem. The proof will be given in
	subsection \ref{subsection_error_analysis}.
	\subsection{The main result}
	In this subsection, we give the main convergence result of the fully discrete scheme \eqref{eq_CHNS_variational_form_phi_discretization}-\eqref{eq_CHNS_variational_form_incompressible_u_discretization} in the following theorem.
	\begin{Theorem}\label{theorem_main_result}
		We assume that the system
		\eqref{eq_CHNS}-\eqref{eq_CHNS_initial_boundary_conditions} admits a unique
		solution $(\phi,\mu,\boldsymbol{u},p)$ with the regularity specified in
		\eqref{eq_regularity_condition}, which is sufficient to guarantee the
		required best-approximation estimates. Let
		$(\phi_h^n,\mu_h^n,\boldsymbol{u}_h^n,p_h^n)$ denote the corresponding discrete
		solution of scheme
		\eqref{eq_CHNS_variational_form_phi_discretization}--
		\eqref{eq_CHNS_variational_form_incompressible_u_discretization}.
		Then, provided
		\begin{equation}
			\begin{aligned}
				&h\leq C\min\biggl\{\epsilon^{\frac{17}{4}},~\epsilon^{\frac{2\alpha_3+2}{r+1}},~\epsilon^{\frac{2\alpha_3}{4-d}},~\left(\frac{\epsilon^{\alpha_5+5d+16}}{C(T)}\right)^{\frac{1}{2k_r}}\biggr\},\\
				&\tau\leq C\epsilon^{5d+16},
			\end{aligned}
		\end{equation}
		where the $\alpha_5$ is defined by 
		\begin{equation}\label{eq_alpha5}
			\alpha_5:=\max\{8\sigma_1+10,2\sigma_2+4,2\sigma_3+4,2\sigma_5\},
		\end{equation} 
		$k_r=2$ if $r=1$, $k_r=r+2$ if $r\geq 2$, 
		under these assumptions, we derive the following optimal error estimates for polynomials of degree $r\geq 1$:
		\begin{align}
			&\max_{1\leq n\leq N}\|\phi^n-\phi_h^n\|^2+\tau\sum_{n=1}^{N}\|\mu^n-\mu_h^n\|^2\leq \frac{1}{\epsilon^{\alpha_5+1}}C(T)(\tau^2+h^{2(r+1)}),\\
			&\max_{1\leq n\leq  N}\|\boldsymbol{u}^n-\boldsymbol{u}_h^n\|^2\leq  \frac{1}{\epsilon^{\alpha_5}}C(T)(\tau^2+\mathcal{E}_h^2),\\
			&\tau\sum_{n=1}^{N}\|\nabla(\boldsymbol{u}^n-\boldsymbol{u}_h^n)\|^2\leq  \frac{1}{\epsilon^{\alpha_5}}C(T)(\tau^2+h^{2(r+1)}),
		\end{align}
		where $C(T):= Ce^{CT}$, and $C > 0$ is a constant independent of $h$, $\tau$, and $\epsilon$. 
	\end{Theorem}
	\subsection{Error analysis}\label{subsection_error_analysis}
	Next, we will prove the optimal convergence theory described above, which requires several steps and lemmas.
	To simplify, we define the following notation: 
	\begin{equation}
		e_{\phi}^n=R_h\phi^n-\phi_h^n,~e_{\mu}^n=\mathcal{G}_h\mu^n-\mu_h^n,~e_{\boldsymbol{u}}^n=\boldsymbol{P}_h\boldsymbol{u}^n-\boldsymbol{u}_h^n,~e_p^n=Q_hp^n-p_h^n,
	\end{equation}
	for $n=0,1,\cdots,N$. Subtracting \eqref{eq_CHNS_variational_form_phi}-\eqref{eq_CHNS_variational_form_inc} from \eqref{eq_CHNS_variational_form_phi_discretization}-\eqref{eq_CHNS_variational_form_incompressible_u_discretization} and using the Ritz and Stokes quasi-projection defined in section \ref{section_projections}, we get the following error equations for $\left(e_{\phi}^n,e_{\mu}^n,e_{\boldsymbol{u}}^n,e_p^n\right)$,
	\begin{flalign}
		\label{eq_error_equations_phi}
		&\left(d_{\tau}e_{\phi}^n,w_h\right)+\left(\nabla e_{\mu}^n,\nabla w_h\right)=\left(d_{\tau}(R_h\phi^n-\phi^n),w_h\right)+\left(\boldsymbol{u}_h^n\nabla\phi_h^{n-1}-\boldsymbol{u}^n\nabla R_h\phi^n, w_h\right)+\left(R_1^n,w_h\right),\\
		\label{eq_error_equations_mu}
		&\left(e_\mu^n,\varphi_h\right)-\epsilon\left(\nabla e_{\phi}^n,\nabla\varphi_h\right)=\left(\mathcal{G}_h\mu^n -\mu^n,\varphi_h\right)
		+\epsilon^{-1}\left(f(\phi^{n})-f(\phi_h^n,\phi_h^{n-1}),\varphi_h\right),\\
		\label{eq_error_estimates_u}
		&\left(d_{\tau}e_{\boldsymbol{u}}^n,\boldsymbol{v}_h\right)+\left(\nabla e_{\boldsymbol{u}}^n,\nabla\boldsymbol{v}_h\right)-\left(e_p^n,\nabla\cdot\boldsymbol{v}_h\right)
		=\left(R_2^n,\boldsymbol{v}_h\right)+\left(d_{\tau}(\boldsymbol{P}_h\boldsymbol{u}^n-\boldsymbol{u}^n),\boldsymbol{v}_h\right)\\
		&\qquad+\boldsymbol{B}(\boldsymbol{u}_h^{n-1},\boldsymbol{u}_h^n,\boldsymbol{v}_h)-
		\boldsymbol{B}\left(\boldsymbol{u}^{n-1},\boldsymbol{u}^n,\boldsymbol{v}_h\right)+\left(\mathcal{G}_h\mu^n\nabla R_h\phi^{n}-\mu_h^n\nabla\phi_h^{n-1},\boldsymbol{v}_h\right),\notag\\
		\label{eq_error_estimates_incompressible_condition}
		&\left(\nabla\cdot e_{\boldsymbol{u}}^n,q_h\right)=0,
	\end{flalign}
	for any $\left(w_h,\varphi_h,\boldsymbol{v}_h,q_h\right)\in S_h^r\times S_h^r\times \boldsymbol{X}_h^{r+1}\times \mathring{S}_h^r$ and $n=1,2,\cdots,N$. In the error equations \eqref{eq_error_equations_phi}-\eqref{eq_error_estimates_incompressible_condition}, the truncation errors $R_1^n$ and $R_2^n$ denote as follows:
	\begin{equation}
		\label{eq_truncation_errors}
		\begin{aligned}
			&R_1^n=d_{\tau}\phi^n-\phi_t(t^n),\\
			&R_2^n=d_{\tau}\boldsymbol{u}^n-\boldsymbol{u}_t(t^n)-((\boldsymbol{u}^n\cdot\nabla)\boldsymbol{u}^n-(\boldsymbol{u}^{n-1}\cdot\nabla)\boldsymbol{u}^n).
		\end{aligned}
	\end{equation}
	From the truncation errors \eqref{eq_truncation_errors}, we have the following lemma.
	\begin{Lemma}\label{lemma_R1_Hneg_L2_boundedness0001}
		For any $h,\tau>0$ and $1\leq n\leq m\leq N$, there exists the following error estimate:
		\begin{equation}
			\tau\sum_{n=1}^{m}\left(\|R_1^n\|_{-1,h}^2+\|R_2^n\|^2\right)\leq C\tau^2(1+\rho_1(\epsilon,d)+\epsilon^{-2\sigma_1}).
		\end{equation}
	\end{Lemma}
	\begin{proof}
		It follows from solution estimates that
		\begin{equation}
			\tau\sum_{n=1}^{m}\|R_1^n\|_{-1,h}^2\leq \frac{1}{\tau}\sum_{n=1}^{m}\left(\int_{t^{n-1}}^{t^n}(s-t^{n-1})ds\right)\left(\int_{t^{n-1}}^{t^n}\|\phi_{tt}(s)\|_{-1,h}^2ds\right)\leq C\rho_1(\epsilon,d)\tau^2.
		\end{equation}
		Using \eqref{eq_ennergy_stable00001} and \eqref{eq_nabla_u_mu_bounded0001}, we next estimate the term $R_2^n$ as follows:
		\begin{equation}
			\begin{aligned}
				\tau\sum_{n=1}^{m}\|R_2^n\|^2\leq&~ \tau\sum_{n=1}^{m}\left(\int_{t^{n-1}}^{t^n}(s-t^{n-1})\boldsymbol{u}_{tt}(s)ds+\|\boldsymbol{u}^n-\boldsymbol{u}^{n-1}\|_{L^4}\|\nabla\boldsymbol{u}^n\|_{L^4}\right)^2
				\leq C(1+\epsilon^{-2\sigma_1})\tau^2.
			\end{aligned}
		\end{equation}
		Combining the above two inequalities, we obtain the conclusion of this lemma.
	\end{proof}
	Next, we need to obtain optimal error estimates for the bounds $\|e_{\boldsymbol{u}}^n\|$ and $\|e_{\phi}^n\|_{-1,h}$. Since the right-hand side of the coupling term $\mu\nabla\phi$ in the overall error estimate contains an estimate for $\|\nabla e_{\mu}^n\|$, we present the following estimate lemma for $\|\nabla e_{\mu}^n\|$.
	\begin{Lemma}
		\label{lemma_nabla_emu}
		For any $h,\tau>0$ and $1\leq n\leq m\leq N$, it holds 
		\begin{equation}
			\begin{aligned}
				\tau\sum_{n=1}^{m}\|\nabla e_{\mu}^n\|^2\leq&~C\tau\sum_{n=1}^{m}\left(\|d_{\tau}e_{\phi}^n\|_{-1,h}^2+\|\nabla e_{\mu}^n\|^2+\|e_{\boldsymbol{u}}^n\|^2+\|e_{\phi}^n\|^2\right)\\ 
				&~+C(1+\rho_0(\epsilon,d)+\rho_1(\epsilon,d))\tau^2+C\epsilon^{-\alpha_3}\mathcal{E}_h^2.
			\end{aligned}
		\end{equation}
	\end{Lemma}
	\begin{proof}
		Taking $w_h=e_{\mu}^n$ in \eqref{eq_error_equations_phi} and summing over $n$ from $1$ to $m$, we obtain
		\begin{equation}
			\label{eq_taking_wh_emu_summ00001}
			\tau\sum_{n=1}^{m}\|\nabla e_{\mu}^n\|^2=\tau\sum_{n=1}^{m}\left(\left(d_{\tau}(R_h\phi^n-\phi^n)-d_{\tau}e_{\phi}^n+R_1^n,e_{\mu}^n\right)+\left(\boldsymbol{u}_h^n\nabla\phi_h^{n-1}-\boldsymbol{u}^n\nabla R_h\phi^n, e_{\mu}^n\right)\right).
		\end{equation}
		Now, we estimate each term on the right-hand side of \eqref{eq_taking_wh_emu_summ00001} one by one, and obtain the following results.
		Using the properties of the approximation \eqref{eq_Ritz_projections_boundedness0002}, we have
		\begin{equation}
			\label{eq_taking_wh_emu_summ00001_part01111}
			\begin{aligned}
				&\tau\sum_{n=1}^{m}\left(d_{\tau}(R_h\phi^n-\phi^n)-d_{\tau}e_{\phi}^n+R_1^n,e_{\mu}^n\right)\\
				\leq&~ C\tau\sum_{n=1}^{m}\left(\|d_{\tau}(R_h\phi^n-\phi^n)\|_{-1,h}+\|d_{\tau}e_{\phi}^n\|_{-1,h}+\|R_1^n\|_{-1,h}\right)\|\nabla e_{\mu}^n\|\\
				\leq &~C\tau\sum_{n=1}^{m}\left(\epsilon^{-\alpha_3}\mathcal{E}_h^2+\|d_{\tau}e_{\phi}^n\|_{-1,h}^2+(1+\rho_1(\epsilon,d))\tau^2+\|\nabla e_{\mu}^n\|^2\right),
			\end{aligned}
		\end{equation}
		and
		\begin{equation}
			\label{eq_taking_wh_emu_summ00001_part02222}
			\begin{aligned}
				&~\tau\sum_{n=1}^{m}\left(\boldsymbol{u}_h^n\nabla\phi_h^{n-1}-\boldsymbol{u}^n\nabla R_h\phi^n,e_{\mu}^n\right)\\
				=&~\tau\sum_{n=1}^{m}\left((R_h\phi^n-R_h\phi^{n-1})\boldsymbol{u}^n+(R_h\phi^{n-1}-\phi_h^{n-1})\boldsymbol{u}^n+\phi_h^{n-1}(\boldsymbol{u}^n-\boldsymbol{u}_h^n),\nabla e_{\mu }^n\right)\\
				\leq&~ C\tau\sum_{n=1}^{m}\left((\|\boldsymbol{u}^n\|_{L^{3}}\|R_h\phi^n-R_h\phi^{n-1}\|_{L^6}+\|\phi_h^{n-1}\|_{L^{\infty}}\|\boldsymbol{u}^n-\boldsymbol{u}_h^n\|)\|\nabla e_{\mu }^n\|\right)\\
				&~+C\tau\sum_{n=1}^{m}\|e_{\phi}^{n-1}\|\|\boldsymbol{u}^n\|_{L^{\infty}}\|\nabla e_{\mu }^n\|\\
				\leq&~C\tau\sum_{n=1}^{m}\left(\rho_0(\epsilon,d)\tau^2+\epsilon^{-\alpha_3}\mathcal{E}_h^2+\|e_{\boldsymbol{u}}^n\|^2+\|e_{\phi}^n\|^2+\|\nabla e_{\mu}^n\|^2\right).
			\end{aligned}
		\end{equation}
		Combining the above inequalities \eqref{eq_taking_wh_emu_summ00001_part01111} and \eqref{eq_taking_wh_emu_summ00001_part02222} with \eqref{eq_taking_wh_emu_summ00001}, we obtain the result of this lemma.
		This completes the proof.
	\end{proof}
	\begin{Lemma}\label{lemma_eu_ephiHnev1000001}
		Under some conditions such as the discrete spectral estimate in Lemma \ref{lemma_spectral_estimates_boundedness}, provided
		\begin{equation}
			\begin{aligned}
				&h\leq C\min\biggl\{\epsilon^{\frac{17}{4}},~\epsilon^{\frac{2\alpha_3+2}{r+1}},~\epsilon^{\frac{2\alpha_3}{4-d}},~\left(\frac{\epsilon^{\alpha_5+5d+16}}{C(T)}\right)^{\frac{1}{2k_r}}\biggr\},\\
				&\tau\leq C\epsilon^{5d+16},
			\end{aligned}
		\end{equation} 
		then, for any $h,\tau>0$ and $1\leq n\leq m\leq N$, it holds
		\begin{equation}
			\label{eq_eu_ephiHnev1000001_boundedness}
			\begin{aligned}
				\|e_{\phi}^m\|_{-1,h}^2&+\|e_{\boldsymbol{u}}^m\|^2+\sum_{n=1}^{m}\|e_{\phi}^n-e_{\phi}^{n-1}\|_{-1,h}^2+\tau\sum_{n=1}^{m}\|\nabla e_{\boldsymbol{u}}^n\|^2\\
				&+\tau\sum_{n=1}^{m}\left(\epsilon^{-1}\|e_{\phi}^n\|_{L^4}^4+\frac{4\epsilon^4}{1-\epsilon^3}\|\nabla e_{\phi}^n\|^2\right)
				\leq C(T)\epsilon^{-\alpha_5}(\tau^2+\mathcal{E}_h^2),
			\end{aligned}
		\end{equation}
		with constant $C$ independent of mesh size $h$, time step $\tau$ and interface width $\epsilon$, and $C(T)$ is defined by $C(T)=  Ce^{CT}$ and $k_r=r+1$ if $r=1$, $k_r=r+2$ if $r\geq 2$ and where $\alpha_5$ is defined in \eqref{eq_alpha5}.
	\end{Lemma}
	\begin{proof}
		To obtain error estimates for $\|e_{\phi}^n\|_{-1,h}$ and $\|e_{\boldsymbol{u}}^n\|$, we will proceed in three steps. 
		
		\textbf{Step1.} First, we will calculate the error estimate for $\|e_{\phi}^n\|_{-1,h}$ using the error equations \eqref{eq_error_equations_phi} and \eqref{eq_error_equations_mu}.
		Taking $w_h=(-\varDelta_h)^{-1}e_{\phi}^n$ in \eqref{eq_error_equations_phi} and $\varphi_h=e_{\phi}^n$ in \eqref{eq_error_equations_mu}, adding the resulting equations, we get the following result:
		\begin{equation}
			\label{eq_Hneg_1norm_sum00001}
			\begin{aligned}
				\frac{1}{2\tau}\|e_{\phi}^n\|_{{H^{-1}}}^2-&~\frac{1}{2\tau}\|e_{\phi}^{n-1}\|_{-1,h}^2+\frac{\tau}{2}\|d_{\tau}e_{\phi}^n\|_{{H^{-1}}}^2+\epsilon\|\nabla e_{\phi}^n\|^2+\frac{1}{\epsilon}\left(f(\phi^n)-f(\phi_h^n,\phi_h^{n-1}),e_{\phi}^n\right)\\
				=&~\left((R_1^n,(-\varDelta_h)^{-1}e_{\phi}^n)+(d_{\tau}(R_h\phi^n-\phi^n),(-\varDelta_h)^{-1}e_{\phi}^n)\right)\\
				&~+\left(\left(\boldsymbol{u}_h^n\nabla\phi_h^{n-1}-\boldsymbol{u}^n\nabla R_h\phi^n,(-\varDelta_h)^{-1}e_{\phi}^n\right)+\left(\mu^n-\mathcal{G}_h\mu^n,e_{\phi}^n\right)\right)
				:=\sum_{i=1}^{4}T_{ai}.
			\end{aligned}
		\end{equation}
		Next, we estimate the terms $T_{a1}\sim T_{a4}$ on the right-hand side of \eqref{eq_Hneg_1norm_sum00001}. Using the Lemma  \ref{lemma_R1_Hneg_L2_boundedness0001}, the first term $T_{a1}$ is bounded by 
		\begin{equation}
			\label{eq_Ta1_boundedness}
			\begin{aligned}
				T_{a1}\leq &~C\|R_1^n\|_{-1,h}\|(-\varDelta_h)^{-1}e_{\phi}^n\|_{H^1}
				\leq C\left(\|R_1^n\|_{-1,h}^2+\|\nabla(-\varDelta_h)^{-1}e_{\phi}^n\|^2\right).
			\end{aligned}
		\end{equation}
		Using the properties of the approximation \eqref{eq_Ritz_projections_boundedness0002}, the second term $T_{a2}$ is estimated by 
		\begin{equation}
			\begin{aligned}
				T_{a2}\leq &~C\|d_{\tau}(R_h\phi^n-\phi^n)\|_{-1,h}\|(-\varDelta_h)^{-1}e_{\phi}^n\|_{H^1}\\
				\leq &~C\left(\|d_{\tau}(R_h\phi^n-\phi^n)\|_{-1,h}^2+\|\nabla(-\varDelta_h)^{-1}e_{\phi}^n\|^2\right)\\
				\leq &~C\mathcal{E}_h^2\|\phi\|_{H^{r+1}}^2+C_0\|\nabla(-\varDelta_h)^{-1}e_{\phi}^n\|^2.
			\end{aligned}
		\end{equation}
		The third term $T_{a3}$ is estimated by 
		\begin{equation}
			\begin{aligned}
				T_{a3}=&~\left((R_h\phi^n-R_h\phi^{n-1})\boldsymbol{u}^n+(R_h\phi^{n-1}-\phi_h^{n-1})\boldsymbol{u}^n+\phi_h^{n-1}(\boldsymbol{u}^n-\boldsymbol{u}_h^n),\nabla(-\varDelta_h)^{-1}e_{\phi}^n\right)\\
				\leq&~ C\left((\|\boldsymbol{u}^n\|_{L^{3}}\|R_h\phi^n-R_h\phi^{n-1}\|_{L^6}+\|\phi_h^{n-1}\|_{L^{\infty}}\|\boldsymbol{u}^n-\boldsymbol{u}_h^n\|)\|\nabla(-\varDelta_h)^{-1}e_{\phi}^n\|\right)\\
				&~+C\|e_{\phi}^{n-1}\|\|\boldsymbol{u}^n\|_{L^{\infty}}\|\nabla(-\varDelta_h)^{-1}e_{\phi}^n\|\\
				\leq &~C\left(\|\nabla\boldsymbol{u}^n\|\|\nabla(\phi^n-\phi^{n-1})\|+\|\phi_h^{n-1}\|_{L^{\infty}}(\|\boldsymbol{u}^n-\boldsymbol{P}_h\boldsymbol{u}^n\|+\|e_{\boldsymbol{u}}^n\|)\right)\|\nabla(-\varDelta_h)^{-1}e_{\phi}^n\|\\
				&~+C\|e_{\phi}^{n-1}\|\|\boldsymbol{u}^n\|_{L^{\infty}}\|\nabla(-\varDelta_h)^{-1} e_{\phi}^n\|\\\
				\leq&~C\left(\tau^2\|\nabla\phi_{t}\|^2+\epsilon^{-2\alpha_3}\mathcal{E}_h^2\|\boldsymbol{u}\|_{H^{r+2}}^2+\|e_{\boldsymbol{u}}^n\|^2\right)+C_0\|\nabla(-\varDelta_h)^{-1}e_{\phi}^n\|^2+\| e_{\phi}^{n-1}\|^2\\
				\leq&~C\left(\rho_0(\epsilon,d)\tau^2+\epsilon^{-2\alpha_3}\mathcal{E}_h^2+\|e_{\boldsymbol{u}}^n\|^2+\|e_{\phi}^{n-1}\|^2\right)+C_0\|\nabla(-\varDelta_h)^{-1}e_{\phi}^n\|^2,
			\end{aligned}
		\end{equation}
		where used equation \eqref{eq_CHNS_variational_form_incompressible_u_discretization}, discrete mass conservation \eqref{eq_mass_conservation}.
		Using the approximate of the Ritz quasi-projection \eqref{eq_Hnev_Gh}, the boundedness of $T_{a4}$ is estimated by 
		\begin{equation}
			\begin{aligned}
				T_{a4}\leq&~ C\|\mu^n-\mathcal{G}_h\mu^n\|_{-1,h}\|e_{\phi}^n\|_{H^1}\\
				\leq&~ \frac{C(1-\beta)\mathcal{E}_h^2}{\beta\epsilon}(\|\phi\|_{H^{r+1}}^2\|\boldsymbol{u}\|_{W^{1,4}}^2+\|\mu\|_{H^{r+1}}^2)+\frac{\beta\epsilon}{1-\beta}\|\nabla e_{\phi}^n\|^2.
			\end{aligned}
		\end{equation}
		Combining the four terms above $T_{a1}\sim T_{a4}$, one gets the following result:
		\begin{equation}
			\label{eq_combining_above_four_terms_result}
			\begin{aligned}
				\frac{1}{2\tau}&\|e_{\phi}^n\|_{{H^{-1}}}^2-\frac{1}{2\tau}\|e_{\phi}^{n-1}\|_{-1,h}^2+\frac{\tau}{2}\|d_{\tau}e_{\phi}^n\|_{{H^{-1}}}^2+\epsilon\|\nabla e_{\phi}^n\|^2+\frac{1}{\epsilon}\left(f(\phi^n)-f(\phi_h^n,\phi_h^{n-1}),e_{\phi}^n\right)\\
				\leq &~C\left(\|R_1^{n}\|_{-1,h}^2+\frac{1-\beta}{\beta\epsilon}\mathcal{E}_h^2+\|e_{\boldsymbol{u}}^n\|^2+\| e_{\phi}^{n-1}\|^2\right)+C_0\|\nabla(-\varDelta_h)^{-1}e_{\phi}^n\|^2+\frac{\beta\epsilon}{1-\beta}\|\nabla e_{\phi}^n\|^2.
			\end{aligned}
		\end{equation}
		
		\textbf{Step 2:} To account for the fourth and fifth terms on the left-hand side of \eqref{eq_combining_above_four_terms_result}, we applied the  spectral estimation of the Cahn--Hilliard operator. Since $f(\phi_h^n,\phi_h^{n-1})=(\phi_h^n)^3-\phi_h^{n-1}$, the nonlinear term can be written as  follows:
		\begin{equation}
			\begin{aligned}
				f(\phi^n)-f(\phi_h^n,\phi_h^{n-1})=f(\phi^n)-f(R_h\phi^n)+f'(R_h\phi^n)e_{\phi}^n-3R_h\phi^n(e_{\phi}^n)^2+(e_{\phi}^n)^3-\tau d_{\tau}\phi_h^n.
			\end{aligned}
		\end{equation}
		Using the Cauchy-Schwarz inequality, the Sobolev embedding, the uniform estimate \eqref{eq_Hnev1_H1_boundedness00001} and the boundedness for \eqref{eq_boundedness_Linfty_phi_Rhphi}, we obtain
		\begin{equation}
			\label{eq_ep_geqphie_00001}
			\begin{aligned}
				\epsilon^{-1}\left(f(\phi^n)-f(\phi_h^n,\phi_h^{n-1}),e_{\phi}^n\right)\geq&~ \epsilon^{-1}\left(f'(R_h\phi^n)e_{\phi}^n,e_{\phi}^n\right)+\epsilon^{-1}\|e_{\phi}^n\|_{L^4}^4\\
				&-\frac{C(1-\beta)}{\beta\epsilon^3}\left(\|f(\phi^n)-f(R_h\phi^n)\|_{-1,h}^2+\tau^2\|d_{\tau}\phi_h^n\|_{-1,h}^2\right)\\
				&-\frac{\beta\epsilon}{1-\beta}\|\nabla e_{\phi}^n\|^2-C\epsilon^{-1}\|e_{\phi}^n\|_{L^3}^3.
			\end{aligned}
		\end{equation}
		Therefore, we obtain the estimate for the right-hand side of the equation \eqref{eq_ep_geqphie_00001} as follows:
		\begin{equation}
			\label{eq_epsilon6_f_boundedness}
			\begin{aligned}
				\frac{1-\beta}{\beta\epsilon^3}&\left(\|f(\phi^n)-f(R_h\phi^n)\|_{-1,h}^2+\tau^2\|d_{\tau}\phi_h^n\|_{-1,h}^2\right)\\
				\leq&~ \frac{1-\beta}{\beta\epsilon^3}\|\phi^n-R_h\phi^n\|_{-1,h}^2+\frac{1-\beta}{\beta\epsilon^3}\tau^2\|d_{\tau}\phi_h^n\|_{-1,h}^2\\
				\leq&~\frac{C(1-\beta)\mathcal{E}_h^2}{\beta\epsilon^3}\|\phi\|_{H^2}^2+\frac{C(1-\beta)}{\beta\epsilon^3}\tau^2\|d_{\tau}\phi_h^n\|_{-1,h}^2.
			\end{aligned}
		\end{equation}
		According to the above inequalities \eqref{eq_ep_geqphie_00001} and \eqref{eq_epsilon6_f_boundedness}, the error equation \eqref{eq_combining_above_four_terms_result} can be recast as 
		\begin{equation}
			\label{eq_error_equation_recasted0001}
			\begin{aligned}
				\frac{1}{2\tau}\|e_{\phi}^n\|_{-1,h}^2-&~\frac{1}{2\tau}\|e_{\phi}^{n-1}\|_{-1,h}^2+\frac{\tau}{2}\|d_{\tau}e_{\phi}^n\|_{{H^{-1}}}^2+\frac{1}{\epsilon}\|e_{\phi}^n\|_{L^4}^4\\
				&+\frac{1-7\beta}{1-\beta}\left(\epsilon\|\nabla e_{\phi}^n\|^2+\frac{1-\beta}{\epsilon}(f'(R_h\phi^n)e_{\phi}^n,e_{\phi}^n)\right)+\frac{4\beta\epsilon}{1-\beta}\|\nabla e_{\phi}^n\|^2\\
				\leq&~C\left(\frac{1-\beta}{\beta\epsilon^3}\tau^2\|d_{\tau}\phi_h^n\|_{-1,h}^2+\frac{1-\beta}{\beta\epsilon^3}\mathcal{E}_h^2\|\phi\|_{H^2}^2+\|e_{\boldsymbol{u}}^n\|^2+\|e_{\phi}^n\|^2+\|e_{\phi}^n\|_{-1,h}^2\right)\\
				&+\frac{C}{\epsilon}\|e_{\phi}^n\|_{L^3}^3+\frac{C\beta}{\epsilon}\left(f'(R_h\phi^n)e_{\phi}^n,e_{\phi}^n\right).
			\end{aligned}
		\end{equation}
		Using the inequality  \eqref{eq_Hnev1_H1_boundedness00001}, \eqref{eq_phiRh_phi_Linfty_bounded}, \eqref{eq_fphiRh_phi_Linfty_bounded} and inequality (2.5) of \cite{2005_FengXiaobing_Numerical_analysis_of_the_Cahn_Hilliard_equation_and_approximation_of_the_Hele_Shaw_problem}, we have
		\begin{equation}
			\frac{C\beta}{\epsilon}|\left(f'(R_h\phi^n)e_{\phi}^n,e_{\phi}^n\right)|\leq \frac{C\beta}{\epsilon}\|e_{\phi}^n\|^2\leq \frac{C(1-\beta)}{\epsilon^2}\|\nabla(-\varDelta_h)^{-1}e_{\phi}^n\|^2+\frac{\beta\epsilon}{1-\beta}\|\nabla e_{\phi}^n\|^2.
		\end{equation}
		The parameter $\beta=\epsilon^3$ seems to give the "best" constant according to \cite{2001_Feng_Prohl_Numerical_analysis_of_the_Cahn_Hilliard_equation_Part_II,2005_FengXiaobing_Numerical_analysis_of_the_Cahn_Hilliard_equation_and_approximation_of_the_Hele_Shaw_problem}.  Here, we provide a detailed explanation of the optimal results for the perturbation parameter $\beta$ in Remark \ref{remark_paraters_beta}.
		
		To control the second term on the right-hand side of \eqref{eq_error_equation_recasted0001}, 
		by using the H\"{o}lder, the discrete Gagliardo--Nirenberg interpolation inequality, the Young inequalities, and the inequality \eqref{eq_Hnev1_H1_boundedness00001}  and the fact $\|e_{\phi}^n\|\leq C\epsilon^{-2\sigma_1}$ when $d=2,3$, we obtain the following result:
		\begin{equation}\label{eq_L3_d2_boundedness}
			\begin{aligned}
				C\epsilon^{-1}\|e_{\phi}^n\|_{L^3}^3\leq &~C\epsilon^{-1}\|e_{\phi}^n\|_{L^4}^2\|e_{\phi}^n\|\\
				\leq &~C\epsilon^{-1}\|e_{\phi}^n\|_{L^4}^{\frac{2d}{2+d}}\|e_{\phi}^n\|^{\frac{6}{2+d}}\|\nabla e_{\phi}^n\|^{\frac{d}{2+d}}\\
				\leq &~C\epsilon^{-1}\|e_{\phi}^n\|_{L^4}^{\frac{2d}{2+d}}\|\nabla e_{\phi}^n\|^{\frac{3+d}{2+d}}\|e_{\phi}^n\|_{-1,h}^{\frac{3}{2+d}}\\
				\leq &~\frac{1}{2\epsilon}\|e_{\phi}^n\|_{L^4}^4+\frac{\epsilon^4}{1-\epsilon^3}\|\nabla e_{\phi}^n\|^2+C\epsilon^{-(5d+16)}\|e_{\phi}^n\|_{-1,h}^6.
			\end{aligned}
		\end{equation}
		Using the spectral estimate in Lemma \ref{lemma_spectral_estimates_boundedness}, and setting $\beta=\epsilon^3$ in \eqref{eq_spectral_estimate_CHNS}, we have
		\begin{equation}\label{eq_epsilon_ephi_boundedness}
			\epsilon\|\nabla e_{\phi}^n\|^2+\frac{1-\epsilon^3}{\epsilon}\left(f'(R_h\phi^n)e_{\phi}^n,e_{\phi}^n\right)\geq -2C_0\|\nabla(-\varDelta_h)^{-1}e_{\phi}^n\|^2-\|e_{\phi}^n\|_{-1,h}^2.
		\end{equation}
		Hence, the error estimate of $C_0\|\nabla(-\varDelta_h)^{-1}e_{\phi}^n\|^2$ is bounded by
		\begin{equation}
			\label{eq_C0_nablaDelta_ephi_boundedness}
			\begin{aligned}
				C_0\|\nabla(-\varDelta_h)^{-1}e_{\phi}^n\|^2\leq&~ C(\|\nabla(-\varDelta_h)^{-1}[e_{\phi}^n-\Pi_he_{\phi}^n]\|^2+\|\nabla (-\varDelta_h)^{-1}\Pi_he_{\phi}^n\|^2)\\
				\leq &~Ch^2\|e_{\phi}^n\|^2+C\|e_{\phi}^n\|_{-1,h}^2\\
				\leq &~\frac{\epsilon^4}{1-\epsilon^3}\|\nabla e_{\phi}^n\|^2+\frac{Ch^4}{\epsilon^4}\|e_{\phi}^n\|_{-1,h}^2+C\|e_{\phi}^n\|_{-1,h}^2\\
				\leq&~\frac{\epsilon^4}{1-\epsilon^3}\|\nabla e_{\phi}^n\|^2+C\|e_{\phi}^n\|_{-1,h}^2,
			\end{aligned}
		\end{equation}
		where we used the condition $h\leq \epsilon$. We get the following result from \eqref{eq_epsilon_ephi_boundedness},
		\begin{equation}
			\epsilon\|\nabla e_{\phi}^n\|^2+\frac{1-\epsilon^3}{\epsilon}\left(f'(R_h\phi^n)e_{\phi}^n,e_{\phi}^n\right)\geq-\frac{\epsilon^4}{1-\epsilon^3}\|\nabla e_{\phi}^n\|^2-C\|e_{\phi}^n\|_{-1,h}^2.
		\end{equation}
		Thus, by combining the above inequalities in \eqref{eq_error_equation_recasted0001}, multiplying both sides by $2\tau$ and summing up the results obtained from $n=1$ to $m$, we obtain 
		\begin{equation}\label{eq_combining_ede00001}
			\begin{aligned}
				&\|e_{\phi}^m\|_{-1,h}^2+2\sum_{n=1}^{m}\|e_{\phi}^n-e_{\phi}^{n-1}\|_{-1,h}^2+\frac{4\epsilon^4\tau}{1-\epsilon^3}\sum_{n=1}^{m}\|\nabla e_{\phi}^n\|^2+\frac{2\tau}{\epsilon}\sum_{n=1}^{m}\|e_{\phi}^n\|_{L^4}^4\\
				\leq &~C\tau\sum_{n=1}^{m}\left(\rho_0(\epsilon,d)\tau^2+\epsilon^{-2\alpha_3}\mathcal{E}_h^2+\|e_{\boldsymbol{u}}^n\|^2+\|e_{\phi}^n\|^2+\|e_{\phi}^n\|_{-1,h}^2+\epsilon^{-(5d+16)}\|e_{\phi}^n\|_{-1,h}^6\right)\\
				&+\|e_{\phi}^0\|_{-1,h}^2.
			\end{aligned}
		\end{equation}
		In view of embedding \eqref{eq_embedding_inequality00001}, 
		choosing $\tau\leq C\epsilon^{5d+16}$, it is easy to check that
		\begin{equation}
			C\tau\left(\epsilon^{-(5d+16)}\|e_{\phi}^n\|_{-1,h}^6+\|e_{\phi}^n\|_{-1,h}^2\right)\leq \frac{1}{2}\|e_{\phi}^n\|_{-1,h}^2.
		\end{equation}
		Under this restriction, we obtain the following result:
		\begin{equation}
			\begin{aligned}
				\|e_{\phi}^m\|_{-1,h}^2+&~2\sum_{n=1}^{m}\|e_{\phi}^n-e_{\phi}^{n-1}\|_{-1,h}^2+\frac{4\epsilon^4\tau}{1-\epsilon^3}\sum_{n=1}^{m}\|\nabla e_{\phi}^n\|^2+\frac{\tau}{\epsilon}\sum_{n=1}^{m}\|e_{\phi}^n\|_{L^4}^4\\
				\leq &~C\tau\sum_{n=1}^{m}\|e_{\boldsymbol{u}}^n\|^2+C\left((1+\rho_0(\epsilon,d)+\rho_1(\epsilon,d))\tau^2+(1+\epsilon^{-\alpha_3})\mathcal{E}_h^2\right)\\
				&~+C\tau\sum_{n=1}^{m-1}\left(\epsilon^{-(5d+16)}\|e_{\phi}^n\|_{-1,h}^6+\|e_{\phi}^n\|_{-1,h}^2\right)+2\|e_{\phi}^0\|_{-1,h}^2.
			\end{aligned}
		\end{equation}
		\textbf{Step 3.} We provide the error estimate for $e_{\boldsymbol{u}}^n$. Taking $\boldsymbol{v}_h=e_{\boldsymbol{u}}^n$ in \eqref{eq_error_estimates_u} and using the equation \eqref{eq_error_estimates_incompressible_condition} yields the following result:
		\begin{equation}
			\label{eq_Tb_boundedness00001}
			\begin{aligned}
				\frac{1}{2}d_{\tau}\|e_{\boldsymbol{u}}^n\|^2+\frac{1}{2\tau}\|e_{\boldsymbol{u}}^n-e_{\boldsymbol{u}}^{n-1}\|^2+\|\nabla e_{\boldsymbol{u}}^n\|^2=\left(R_2^n,e_{\boldsymbol{u}}^n\right)+\left(d_{\tau}(\boldsymbol{P}_h\boldsymbol{u}^n-\boldsymbol{u}^n),e_{\boldsymbol{u}}^n\right)&\\
				+\boldsymbol{B}\left(\boldsymbol{u}_h^{n-1},\boldsymbol{u}_h^n,e_{\boldsymbol{u}}^n\right)-\boldsymbol{B}\left(\boldsymbol{u}^{n-1},\boldsymbol{u}^n,e_{\boldsymbol{u}}^n\right)+\left(\mathcal{G}_h\mu^n\nabla R_h\phi^n-\mu_h^n\nabla\phi_h^{n-1},e_{\boldsymbol{u}}^n\right)&:=\sum_{i=1}^{4}T_{bi}.
			\end{aligned}
		\end{equation}
		We first estimate $T_{bi},~1\leq i\leq 4$, below. By the Lemma \ref{lemma_R1_Hneg_L2_boundedness0001}, it is easy to see that
		\begin{equation}\label{eq_Tb1_boundedness00001}
			T_{b1}\leq C\|R_2^n\|\|e_{\boldsymbol{u}}^n\|\leq C\left(\|R_2^n\|^2+\|e_{\boldsymbol{u}}^n\|^2\right),
		\end{equation}
		and 
		\begin{equation}
			T_{b2}\leq C\|d_{\tau}(\boldsymbol{P}_h\boldsymbol{u}^n-\boldsymbol{u}^n)\|\|e_{\boldsymbol{u}}^n\|\leq C(\mathcal{E}_h^2+\|e_{\boldsymbol{u}}^n\|^2).
		\end{equation}
		In the following, we rewrite $T_{b3}$ as follows:
		\begin{equation}
			\begin{aligned}
				T_{b3}=&-\boldsymbol{B}\left(\boldsymbol{u}^{n-1}-\boldsymbol{u}_h^{n-1},\boldsymbol{u}^n,e_{\boldsymbol{u}}^n\right)-\boldsymbol{B}\left(\boldsymbol{u}_h^{n-1},\boldsymbol{u}^n-\boldsymbol{u}_h^n,e_{\boldsymbol{u}}^n\right)\\
				=&-\boldsymbol{B}\left(\boldsymbol{u}^{n-1}-\boldsymbol{P}_h\boldsymbol{u}^{n-1},\boldsymbol{u}^n,e_{\boldsymbol{u}}^n\right)-\boldsymbol{B}\left(e_{\boldsymbol{u}}^{n-1},\boldsymbol{P}_h\boldsymbol{u}^n,e_{\boldsymbol{u}}^n\right)\\
				&-\boldsymbol{B}\left(\boldsymbol{P}_h\boldsymbol{u}^{n-1},\boldsymbol{u}^n-\boldsymbol{P}_h\boldsymbol{u}^n,e_{\boldsymbol{u}}^n\right)-\boldsymbol{B}\left(\boldsymbol{u}_h^{n-1},e_{\boldsymbol{u}}^n,e_{\boldsymbol{u}}^n\right):=\sum_{i=1}^{4}T_{b3,i}.
			\end{aligned}
		\end{equation}
		It follows from the skew-symmetry that $T_{b3,4}=0$, and the two terms $T_{b3,1}$ and $T_{b3,2}$ are estimated by
		\begin{equation}
			\begin{aligned}
				T_{b3,1} =&-\frac{1}{2} \left( (\boldsymbol{u}^{n-1} - \boldsymbol{P}_{h} \boldsymbol{u}^{n-1}) \cdot \nabla \boldsymbol{u}^n, e_{\boldsymbol{u}}^n \right) + \frac{1}{2} \left( (\boldsymbol{u}^{n-1} - \boldsymbol{P}_h \boldsymbol{u}^{n-1}) \cdot \nabla e_{\boldsymbol{u}}^n, \boldsymbol{u}^n \right)	\\
				\leq &~C \|\boldsymbol{u}^{n-1} - \boldsymbol{P}_h \boldsymbol{u}^{n-1}\| \|\nabla \boldsymbol{u}^n\|_{L^3} \|e_{\boldsymbol{u}}^n\|_{L^6}			
				+ C \|\boldsymbol{u}^{n-1} - \boldsymbol{P}_h\boldsymbol{u}^{n-1}\| \|\nabla e_{\boldsymbol{u}}^n\|\|\boldsymbol{u}^n\|_{L^\infty}\\
				\leq &~C \|\nabla e_{\boldsymbol{u}}^n\|^2 + C\epsilon^{-2\alpha_3}\mathcal{E}_h^2,
			\end{aligned}
		\end{equation}
		and
		\begin{equation}
			\begin{aligned}
				T_{b3,2}=&-\frac{1}{2} \left( e_{\boldsymbol{u}}^{n-1} \cdot \nabla \boldsymbol{P}_h\boldsymbol{u}^n, e_{\boldsymbol{u}}^n \right) + \frac{1}{2} \left( e_{\boldsymbol{u}}^{n-1} \cdot \nabla e_{\boldsymbol{u}}^n, \boldsymbol{P}_h\boldsymbol{u}^n \right)\\				
				\leq &~C \|e_{\boldsymbol{u}}^{n-1}\| \|\nabla \boldsymbol{P}_h\boldsymbol{u}^n\|_{L^3} \|e_{\boldsymbol{u}}^n\|_{L^6} + C \|e_{\boldsymbol{u}}^{n-1}\| \|\nabla e_{\boldsymbol{u}}^n\| \|\boldsymbol{P}_h\boldsymbol{u}^n\|_{L^\infty}\\
				\leq &~ C\|\nabla e_{\boldsymbol{u}}^n\|^2 + C  \|e_{\boldsymbol{u}}^{n-1}\|^2,
			\end{aligned}
		\end{equation}
		where we used the inequality \eqref{eq_u_boundedness_error_estimate001}, the projection error estimates \eqref{eq_u_boundedness_error_estimate002}, \eqref{eq_deltau_boundedness_error_estimate001}, and \eqref{eq_Phu_boundedness00001}. Using integration by parts, we further obtain the following result:
		\begin{equation}
			\label{eq_Tb33_boundedness00001}
			\begin{aligned}
				T_{b3,3}=&-\frac{1}{2} \left( \boldsymbol{P}_h \boldsymbol{u}^{n-1} \cdot \nabla (\boldsymbol{u}^n - \boldsymbol{P}_h \boldsymbol{u}^n), e_{\boldsymbol{u}}^n \right) + \frac{1}{2} \left( \boldsymbol{P}_h \boldsymbol{u}^{n-1} \cdot \nabla e_{\boldsymbol{u}}^n, \boldsymbol{u}^n - \boldsymbol{P}_h \boldsymbol{u}^n \right)\\				
				= &~\frac{1}{2} \left( (\nabla \cdot \boldsymbol{P}_h \boldsymbol{u}^{n-1})(\boldsymbol{u}^n - \boldsymbol{P}_h \boldsymbol{u}^n), e_{\boldsymbol{u}}^n \right) + \frac{1}{2} \left( \boldsymbol{u}^n - \boldsymbol{P}_h \boldsymbol{u}^n, \boldsymbol{P}_h \boldsymbol{u}^{n-1} \cdot \nabla e_{\boldsymbol{u}}^n \right)	\\			
				&+ \frac{1}{2} \left( \boldsymbol{P}_h \boldsymbol{u}^{n-1} \cdot \nabla e_{\boldsymbol{u}}^n, \boldsymbol{u}^n - \boldsymbol{P}_h \boldsymbol{u}^n \right)\\
				\leq&~ C \left\| \nabla \cdot \boldsymbol{P}_h \boldsymbol{u}^{n-1} \right\|_{L^3} \left\| \boldsymbol{u}^n - \boldsymbol{P}_h \boldsymbol{u}^n \right\| \left\| e_{\boldsymbol{u}}^n \right\|_{L^6}\\
				&+ C \left\| \boldsymbol{u}^n - \boldsymbol{P}_h \boldsymbol{u}^n \right\| \left\| \boldsymbol{P}_h \boldsymbol{u}^{n-1} \right\|_{L^\infty} \left\| \nabla e_{\boldsymbol{u}}^n \right\| \\
				&+ C \left\| \boldsymbol{P}_h \boldsymbol{u}^{n-1} \right\|_{L^\infty} \left\| \nabla e_{\boldsymbol{u}}^n \right\| \left\| \boldsymbol{u}^n - \boldsymbol{P}_h \boldsymbol{u}^n \right\|\\
				\leq&~ C\| \nabla e_{\boldsymbol{u}}^n\|^2 + C\epsilon^{-2\alpha_3}\mathcal{E}_h^2.
			\end{aligned}
		\end{equation}
		Next,  the term $T_{b4}$ can be estimated by 
		\begin{equation}
			\label{eq_Tb4_boundedness00001}
			\begin{aligned}
				T_{b4}=&~\left(e_{\mu}^n\nabla R_h\phi^n+\mu_h^n\nabla e_{\phi}^{n-1},e_{\boldsymbol{u}}^n\right)+\left(\mu_h^n\nabla(R_h\phi^n-R_h\phi^{n-1}),e_{\boldsymbol{u}}^n\right)\\
				\leq &~C(\|e_{\mu}^n\|_{L^4}\|\nabla R_h\phi^n\|_{L^4}\|e_{\boldsymbol{u}}^n\|+\| e_{\phi}^{n-1}\|(\|e_{\boldsymbol{u}}^n\|_{L^4}\|\nabla\mu_h^n\|_{L^4}+\|\mu_h^n\|_{L^{\infty}}\|\nabla e_{\boldsymbol{u}}^n\|))\\
				&~+C\|\mu_h^n\|_{L^{\infty}}\|\nabla(R_h\phi^n-R_h\phi^{n-1})\|\|e_{\boldsymbol{u}}^n\|\\
				\leq&~C\|\nabla e_{\mu}^n\|\|\nabla\phi^n\|\|\nabla e_{\boldsymbol{u}}^n\|+C\| e_{\phi}^n\|\|\mu_h^n\|_{H^{2}}\|\nabla e_{\boldsymbol{u}}^n\|\\
				&~+C\|\mu_h^n\|_{L^{\infty}} \|\nabla(R_h\phi^n-R_h\phi^{n-1})\|\| e_{\boldsymbol{u}}^n\|\\
				\leq &~C\left(\|\nabla e_{\mu}^n\|^2+\|e_{\boldsymbol{u}}^n\|^2+\|\nabla e_{\boldsymbol{u}}^n\|^2+\|\mu_h^n\|_{H^2}^2\| e_{\phi}^n\|^2+\tau^2\|\mu_h^n\|_{L^{\infty}}^2\|\nabla\phi_{t}\|^2\right).
			\end{aligned}
		\end{equation}
		Combining the previous inequalities \eqref{eq_Tb1_boundedness00001}-\eqref{eq_Tb4_boundedness00001}, we obtain the following result:
		\begin{equation}
			\label{eq_combining_results_boundedness00001}
			\begin{aligned}
				\frac{1}{2}d_{\tau}\|e_{\boldsymbol{u}}^n\|^2&+\frac{1}{2\tau}\|e_{\boldsymbol{u}}^n-e_{\boldsymbol{u}}^{n-1}\|^2+\|\nabla e_{\boldsymbol{u}}^n\|^2
				\leq C\left(\|R_2^n\|^2+\|e_{\boldsymbol{u}}^n\|^2+\|\nabla e_{\boldsymbol{u}}^n\|^2\right)\\
				&+C\left(\|\nabla e_{\mu}^n\|^2+\|\mu_h^n\|_{H^2}^2\| e_{\phi}^n\|^2+\tau^2\|\mu_h^n\|_{L^{\infty}}^2\|\nabla\phi_{t}\|^2+\epsilon^{-2\alpha_3}\mathcal{E}_h^2\right).
			\end{aligned}
		\end{equation}
		Applying the Agmon inequality $\|\mu_h^n\|_{L^{\infty}}\leq C\|\mu_h^n\|^{1-d/4}\|\varDelta_h\mu_h^n\|^{d/4}+C\|\mu_h^n\|$ and multiplying \eqref{eq_combining_results_boundedness00001} by $2\tau$ and taking the summation of $n$ from $1$ to $m$, we derive
		\begin{equation}
			\begin{aligned}
				\|e_{\boldsymbol{u}}^m\|^2&+\tau\sum_{n=1}^{m}\|\nabla e_{\boldsymbol{u}}^n\|^2\leq \|e_{\boldsymbol{u}}^0\|^2+C\tau\sum_{n=1}^{m}\left(\|R_2^n\|^2+\|e_{\boldsymbol{u}}^n\|^2+\|\nabla e_{\boldsymbol{u}}^n\|^2\right)\\
				&+C\tau\sum_{n=1}^{m}\left(\|\nabla e_{\mu}^n\|^2+\|\mu_h^n\|_{H^2}^2\|e_{\phi}^n\|^2+\tau^2(\|\mu_h^n\|^{\frac{4-d}{2}}\|\varDelta_h\mu_h^n\|^{\frac{d}{2}}+\|\mu_h^{n}\|^2)\|\nabla\phi_{t}\|^2+\epsilon^{-2\alpha_3}\mathcal{E}_h^2\right).
			\end{aligned}
		\end{equation}
		Then, combining \textbf{Step 1} to \textbf{Step 3} yields the following result
		\begin{equation}
			\label{eq_combiningStep1_3_boundedness}
			\begin{aligned}
				\|e_{\phi}^m\|_{-1,h}^2+\|e_{\boldsymbol{u}}^m\|^2&~+\tau^2\sum_{n=1}^{m}\|d_{\tau}e_{\phi}^n\|_{-1,h}^2+\tau\sum_{n=1}^{m}\|\nabla e_{\boldsymbol{u}}^n\|^2+\tau\sum_{n=1}^{m}\left(\epsilon^{-1}\|e_{\phi}^n\|_{L^4}^4+\frac{4\epsilon^4}{1-\epsilon^3}\|\nabla e_{\phi}^n\|^2\right)\\
				\leq &~C\tau\sum_{n=1}^{m}\left(\|e_{\boldsymbol{u}}^n\|^2+\|\nabla e_{\boldsymbol{u}}^n\|^2+\epsilon^{-2\alpha_3}\| e_{\phi}^n\|^2+\|\nabla e_{\mu}^n\|^2\right)\\
				&~+C\tau\sum_{n=1}^{m-1}\left(\epsilon^{-(5d+16)}\|e_{\phi}^n\|_{-1,h}^6+\|e_{\phi}^n\|_{-1,h}^2\right)+\|e_{\phi}^0\|_{-1,h}^2\\
				&~+C\left((1+\rho_0(\epsilon,d)+\rho_1(\epsilon,d))\tau^2+(1+\epsilon^{-2\alpha_3})\mathcal{E}_h^2\right).
			\end{aligned}
		\end{equation} 		
		Therefore, based on the definition \eqref{eq_Hnev1_definition}, inequality \eqref{eq_Hnev1_H1_boundedness00001} in the Lemma \ref{lemma_H1Hnev1}, inequality \eqref{eq_Linf_muh} in the Lemma \ref{lemma_dtau_phi_Linf_bounded} and inequality \eqref{eq_C0_nablaDelta_ephi_boundedness}, we obtain the estimation results for the term $\epsilon^{-2\alpha_3}\|e_{\phi}^n\|^2$ on the right-hand side of \eqref{eq_combiningStep1_3_boundedness} as follows:
		\begin{equation}
			\begin{aligned}
				\epsilon^{-2\alpha_3}\|e_{\phi}^n\|^2\leq &~\epsilon^{-2\alpha_3}\|\nabla e_{\phi}^n\|\|\nabla(-\varDelta_h)^{-1}e_{\phi}^n\|\\
				\leq &~\frac{\epsilon^4}{2(1-\epsilon^3)}\|\nabla e_{\phi}^n\|^2+\frac{1-\epsilon^3}{2\epsilon^{4\alpha_3+4}}\|\nabla (-\varDelta_h)^{-1}e_{\phi}^n\|^2\\
				\leq &~\frac{\epsilon^4}{1-\epsilon^3}\|\nabla e_{\phi}^n\|^2+\frac{1-\epsilon^3}{2\epsilon^{4\alpha_3+4}}\|\nabla (-\varDelta_h)^{-1}(e_{\phi}^n-\Pi_he_{\phi}^n+\Pi_he_{\phi}^n)\|^2\\
				\leq &~\frac{\epsilon^4}{1-\epsilon^3}\|\nabla e_{\phi}^n\|^2+\frac{1-\epsilon^3}{2\epsilon^{4\alpha_3+4}}\|\nabla (-\varDelta_h)^{-1}(e_{\phi}^n-\Pi_he_{\phi}^n)\|^2+2\|\nabla (-\varDelta_h)^{-1}\Pi_he_{\phi}^n\|^2\\
				\leq &~\frac{\epsilon^4}{1-\epsilon^3}\|\nabla e_{\phi}^n\|^2+\frac{C(1-\epsilon^3)h^{2(r+1)}}{2\epsilon^{4\alpha_3+4}}\|e_{\phi}^n\|_{-1,h}^6+C\|e_{\phi}^n\|_{-1,h}^2\\
				\leq &~\frac{\epsilon^4}{1-\epsilon^3}\|\nabla e_{\phi}^n\|^2+C\|e_{\phi}^n\|_{-1,h}^2,
			\end{aligned}
		\end{equation}
		provided by condition $h<\epsilon^{\frac{2\alpha_3+2}{r+1}}$. Here, we use the inequality $c\|a-b+b\|^2\leq 2c\|a-b\|^2+2\|b\|^2,~ 0<c\leq 1$ and address the error estimate for the term $\epsilon^{-2\alpha_3}\|e_{\phi}^n\|^2$ on the right-hand side of \eqref{eq_combiningStep1_3_boundedness}. For the term $\epsilon^{-2\alpha_3}\|e_{\phi}^n\|^2$, we only have the error estimates of the $H^{-1}$- and $H^1$-norms, but not the $L^2$-norm error estimate. Therefore, we employ the Lemma \ref{lemma_H1Hnev1} together with the relationship between the mesh size $h$ and the interface thickness $\epsilon$ to derive an estimate of the $L^2$-norm for $e_{\phi}^n$ in terms of its $H^{-1}$- and $H^1$-norms.
		By using the Lemma \ref{lemma_nabla_emu} and the nonlinear Gr\"{o}nwall's Lemma \ref{lemma_nonlinear_gronwall_inequalities}, we define the slack variable $d_m\geq 0$ such that
		\begin{equation}
			\label{eq_combining_nabla_dm}
			\begin{aligned}
				\|e_{\phi}^m\|_{-1,h}^2&+\|e_{\boldsymbol{u}}^m\|^2+\sum_{n=1}^{m}\|e_{\phi}^n-e_{\phi}^{n-1}\|_{-1,h}^2+\tau\sum_{n=1}^{m}\left(\|\nabla e_{\boldsymbol{u}}^n\|^2\right)\\
				&\qquad+\tau\sum_{n=1}^{m}\left(\epsilon^{-1}\|e_{\phi}^n\|_{L^4}^4+\frac{4\epsilon^4}{1-\epsilon^3}\|\nabla e_{\phi}^n\|^2+\|\nabla e_{\mu}^n\|^2\right)+d_m\\
				=&~C\tau\sum_{n=1}^{m}\left(\|e_{\boldsymbol{u}}^n\|^2+\|\nabla e_{\boldsymbol{u}}^n\|^2+\|d_{\tau}e_{\phi}^n\|_{-1,h}^2+\|\nabla e_{\mu}^n\|^2\right)\\
				&~+C\tau\sum_{n=1}^{m-1}\left(\|e_{\phi}^n\|_{-1,h}^2+\epsilon^{-(5d+16)}\|e_{\phi}^n\|_{-1,h}^6\right)\\
				&+C\left((1+\rho_0(\epsilon,d)+\rho_1(\epsilon,d))\tau^2+(1+\epsilon^{-2\alpha_3})\mathcal{E}_h^2\right).
			\end{aligned}
		\end{equation} 
		We also define $\{S_m\geq 0\}_{m\geq1}$ by
		\begin{equation}
			\begin{aligned}
				S_m:=	\|e_{\phi}^m\|_{-1,h}^2&+\|e_{\boldsymbol{u}}^m\|^2+\sum_{n=1}^{m}\|e_{\phi}^n-e_{\phi}^{n-1}\|_{-1,h}^2+\tau\sum_{n=1}^{m}\left(\|\nabla e_{\boldsymbol{u}}^n\|^2+\|\nabla e_{\mu}^n\|^2\right)\\
				&+\tau\sum_{n=1}^{m}\left(\epsilon^{-1}\|e_{\phi}^n\|_{L^4}^4+\frac{4\epsilon^4}{1-\epsilon^3}\|\nabla e_{\phi}^n\|^2\right)+d_m,
			\end{aligned}
		\end{equation}
		and \eqref{eq_combining_nabla_dm} implies that $S_1 \leq C\left((1+\rho_0(\epsilon,d)+\rho_1(\epsilon,d))\tau^2+(1+\epsilon^{-2\alpha_3})\mathcal{E}_h^2\right)$.
		Then, for all $m\geq 1$, we have 
		\begin{equation}
			\begin{aligned}
				S_{m+1}-S_m\leq C \tau S_m+C\tau\epsilon^{-(5d+16)}S_m^3.
			\end{aligned}
		\end{equation}
		As long as $\tau\to 0$ is sufficiently small or the absorption coefficient is less than 1, we obtain: 
		\begin{equation}
			S_{m}\leq C\tau\sum_{n=1}^{m-1}S_n+C\tau\sum_{n=1}^{m-1}\epsilon^{-(5d+16)}S_n^3+C\left((1+\rho_0(\epsilon,d)+\rho_1(\epsilon,d))\tau^2+(1+\epsilon^{-2\alpha_3})\mathcal{E}_h^2\right),
		\end{equation}
		Applying Lemma \ref{lemma_nonlinear_gronwall_inequalities} and the sequence $\{S_m\}_{m\geq 1}$ defined above, for all $m\geq 1$, we have 
		\begin{equation}
			S_m\leq a_m^{-1}\biggl\{S_1^{-2}-2\epsilon^{-(5d+16)}\tau\sum_{n=1}^{m-1}a_{n+1}^{-2}\biggr\}^{-\frac{1}{2}},
		\end{equation} 
		provided that 
		\begin{equation}
			S_1^{-2}-2\epsilon^{-(5d+16)}\tau\sum_{n=1}^{m-1}a_{n+1}^{-2}>0.
		\end{equation}
		Thus, we note that the terms $a_n(1\leq n\leq m)$ are all bounded as $\tau \to 0$, and under the mesh constraint stated in the theorem, hold
		\begin{equation}
			a_m=\frac{1}{(1+C\tau)^{m-1}}\geq e^{-CT},~\mathcal{E}_h^{2}+\tau^2\leq \frac{\epsilon^{\alpha_5+5d+16}}{C(T)},
		\end{equation}
		and 
		\begin{equation}
			S_m\leq Ce^{CT}S_1\leq C(T)\epsilon^{-\alpha_5}(\tau^2+\mathcal{E}_h^2),
		\end{equation}
		where $C(T):= Ce^{CT}$. Recalling the constraint on $h$ from \eqref{eq_spectral_estimate_CHNS} with $\beta = \epsilon^3$, we obtain that provided 
		\begin{equation}
			\begin{aligned}
				&h\leq C\min\biggl\{\epsilon^{\frac{17}{4}},~\epsilon^{\frac{2\alpha_3+2}{r+1}},~\epsilon^{\frac{2\alpha_3}{4-d}},~\left(\frac{\epsilon^{\alpha_5+5d+16}}{C(T)}\right)^{\frac{1}{2k_r}}\biggr\},\\
				&\tau\leq C\epsilon^{5d+16},
			\end{aligned}
		\end{equation}
		where $k_r=2$ if $r=1$, $k_r=r+2$ if $r\geq 2$, and
		the estimate \eqref{eq_eu_ephiHnev1000001_boundedness} is true.
		The proof of this lemma is completed.
	\end{proof}
	\begin{Remark}
		In inequality \eqref{eq_L3_d2_boundedness}, we employ the three-term discrete Gagliardo–Nirenberg inequality to establish a relationship between the time $\tau$ and the interface width $\epsilon$ as $\tau\leq C\epsilon^{5d+16}$. It should be noted that the resulting power exponent is lower than the corresponding value reported in the existing literature \cite{2023_ChenGang_On_the_superconvergence_of_a_hybridizable_discontinuous_Galerkin_method_for_the_Cahn_Hilliard_equation,2015_FengXiaobing_Analysis_of_symmetric_interior_penalty_discontinuous_Galerkin_methods_for_the_Allen_Cahn_equation_and_the_mean_curvature_flow}. This indicates that, compared with existing results, our approach yields a sharper (weaker) power-law dependence between time and interface width.
	\end{Remark}		\begin{Remark}\label{remark_paraters_beta}
		The choice of the scaling factor $\beta=\epsilon^3$ in the application of the 
		spectral estimate is crucial to preserve the coercivity of the gradient term 
		while controlling the nonlinear contribution. In fact, we consider a general scaling 
		factor of the form
		\begin{equation}
			\beta=\epsilon^\alpha,\qquad \alpha>0.
		\end{equation}
		We recast the continuous Cahn-Hilliard operator to get
		\begin{equation}
			\mathcal{L}(\Phi^m)
			=(1-\beta)\mathcal{L}(\Phi^m)+\beta \mathcal{L}(\Phi^m),
		\end{equation}
		where
		\begin{equation}				
			\mathcal{L}(\Phi^m)= \epsilon\|\nabla\Phi^m\|_{L^2}^2+\frac{1-\beta}{\epsilon}
			(f'(R_h\phi^m)\Phi^m,\Phi^m),
		\end{equation}
		the first part can be estimated by Lemma \ref{lemma_spectral_estimates_boundedness},
		which yields
		\begin{equation}
			(1-\beta)\mathcal{L}(\Phi^m) \geq -(1-\beta)(2C_0+1) \|\nabla\varDelta_h^{-1}\Phi^m\|_{L^2}^2 .
		\end{equation}
		The remaining part provides a positive contribution of 
		\begin{equation}
			\beta\epsilon\|\nabla\Phi^m\|_{L^2}^2 = 	\epsilon^{\alpha+1}\|\nabla\Phi^m\|_{L^2}^2,
		\end{equation}
		which is needed to absorb the higher-order nonlinear terms. Meanwhile, the
		remaining potential term satisfies
		\begin{equation}
			\frac{\beta}{\epsilon}(f'(R_h\phi^m)\Phi^m,\Phi^m)=\epsilon^{\alpha-1}(f'(R_h\phi^m)\Phi^m,\Phi^m).
		\end{equation}
		By using the $H^{-1}$-norm error estimate and Young's inequality, this term generates a
		contribution of $ \epsilon^{2\alpha-2} \|\nabla\Phi^m\|_{L^2}^2 $.
		Therefore, in order to absorb this term by the remaining gradient contribution,
		it is necessary that 
		\begin{equation}
			\epsilon^{2\alpha-2}\leq \epsilon^{\alpha+1},
		\end{equation}
		which requires $2\alpha-2\geq\alpha+1,$
		or equivalently, $\alpha\geq3.$
		Hence, any splitting weight of the form $1-\epsilon^\alpha$ with $\alpha\geq3$ is admissible. Among these choices, the smallest admissible exponent $\alpha=3$ preserves the largest possible portion of the gradient 	energy term, leading to the sharpest estimate. Consequently, throughout the above analysis, we choose 
		\begin{equation}
			\beta=\epsilon^3.
		\end{equation}
		This choice balances the spectral estimate and the nonlinear absorption argument and provides the optimal constant in the subsequent error estimates.
	\end{Remark}
	
	The exponent $\alpha=3$ in $\beta=\epsilon^{\alpha}$ is not merely convenient: it is the unique optimal choice. After applying the spectral estimate of Lemma \ref{lemma_spectral_estimates_boundedness} to the leading part of the quadratic form, the remaining potential term satisfies, by $f'(\phi)=3\phi^2-1\ge-1$,
	\begin{equation}
		\frac{\beta}{\epsilon}\big(f'(R_h\phi^m)\Phi^m,\Phi^m\big)\ \ge\ -\frac{\beta}{\epsilon}\|\Phi^m\|^2\ \ge\ -\frac{\beta\epsilon}{2}\|\nabla\Phi^m\|^2-\frac{\beta}{2\epsilon^{3}}\|\Phi^m\|_{-1,h}^{2},
	\end{equation}
	where we used $\|\Phi^m\|^2\le\|\Phi^m\|_{-1,h}\|\nabla\Phi^m\|$ and Young's inequality. The first term is absorbed by the gradient budget $\beta\epsilon\|\nabla\Phi^m\|^2$ retained in \eqref{eq_ep_geqphie_00001}; the second enters the linear part of the nonlinear Gr\"onwall inequality of Lemma \ref{lemma_nonlinear_gronwall_inequalities}, and therefore its coefficient $\beta/(2\epsilon^{3})=\epsilon^{\alpha-3}/2$ must be independent of $\epsilon$: otherwise $a_m^{-1}\le\mathrm e^{CT\epsilon^{\alpha-3}}$ and the error constants would again depend exponentially on $\epsilon^{-1}$. This forces $\alpha\ge3$. In contrast, the gradient budget $\epsilon^{\alpha+1}$, the mesh conditions of Lemma \ref{lemma_spectral_estimates_boundedness}, and the exponent of the cubic coefficient in \eqref{eq_L3_d2_boundedness} deteriorate as $\alpha$ increases. Thus, $\alpha=3$ is the optimal choice within the present error-analysis framework, and we take $\beta=\epsilon^{3}$ throughout.
	
	Next, we perform error estimates of the phase-field and chemical potential function in the $L^2$ norm. For this purpose, we proceed to error estimates in the following lemma.
	\begin{Lemma}
		\label{lemma_emuphi}
		Under the same conditions as Lemma \ref{lemma_eu_ephiHnev1000001}, the following error estimate holds for any $\tau,h>0$ and $1\leq n\leq m\leq N$,
		\begin{equation}\label{eq_L2_phimu_boundedness}
			\|e_{\phi}^m\|^2+\tau\sum_{n=1}^{m}\|e_{\mu}^n\|^2\leq \frac{C(T)}{\epsilon^{\alpha_5+1}}(\tau^2+h^{2(r+1)}),
		\end{equation}
		where $C$ is the positive constant independent of the mesh size $h$, time step $\tau$, and interface width $\epsilon$.
	\end{Lemma}
	\begin{proof}
		Taking $w_h=\epsilon e_{\phi}^n$ in \eqref{eq_error_equations_phi} and $\varphi_h=e_{\mu}^n$ in \eqref{eq_error_equations_mu}, respectively, and adding the two results, we obtain
		\begin{equation}
			\label{eq_wvarphi_boundedness0011}
			\begin{aligned}
				\frac{\epsilon}{2\tau}&(\|e_{\phi}^n\|^2-\|e_{\phi}^{n-1}\|^2+\|e_{\phi}^n-e_{\phi}^{n-1}\|^2)+\|e_{\mu}^n\|^2=\epsilon(d_{\tau}(R_h\phi^n-\phi^n),e_{\phi}^n)\\
				&~+\epsilon(R_h\phi^n\boldsymbol{u}^n-\phi_h^{n-1}\boldsymbol{u}_h^n,\nabla e_{\phi}^n)+\epsilon(R_1^n,e_{\phi}^n)+(\mathcal{G}_h\mu^n-\mu^n,e_{\mu}^n)\\
				&~+\frac{1}{\epsilon}(f(\phi^n)-f(\phi_h^n,\phi_h^{n-1}),e_{\mu}^n)
				:=\sum_{i=1}^{5}T_{ci}.
			\end{aligned}
		\end{equation}
		We obtain the error estimates for the three terms $T_{c1},~T_{c3}$ and $T_{c4}$ as follows:
		\begin{equation}
			\label{eq_error_bounded_Tc134}
			\begin{aligned}
				T_{c1}+T_{c3}+T_{c4}\leq &~\epsilon(\|d_{\tau}(R_h\phi^n-\phi^n)\|_{-1,h}+\|R_1^n\|_{-1,h})\|\nabla e_{\phi}^n\|+\|\mathcal{G}_h\mu^n-\mu^n\|_{-1,h}\|\nabla e_{\mu}^n\|\\
				\leq&~\|d_{\tau}(R_h\phi^n-\phi^n)\|_{-1,h}^2+\|R_1^n\|_{-1,h}^2+\epsilon\|\nabla e_{\phi}^n\|^2+\|\mathcal{G}_h\mu^n-\mu^n\|_{-1,h}^2+\|\nabla e_{\mu}^n\|^2.
			\end{aligned}
		\end{equation}
		The term $T_{c2}$ is bounded by
		\begin{equation}
			\label{eq_Tc3_boundedness}
			\begin{aligned}
				T_{c2}=&~\epsilon\left((\boldsymbol{u}^n-\boldsymbol{u}_h^n)\phi^n+\boldsymbol{u}_h^n(R_h\phi^n-R_h\phi^{n-1}+R_h\phi^{n-1}-\phi_h^{n-1}),\nabla e_{\phi}^n\right)\\
				\leq &~C\epsilon\left(\|\boldsymbol{u}^n-\boldsymbol{P}_h\boldsymbol{u}^n+e_{\boldsymbol{u}}^n\|_{L^4}\|\phi^n\|_{L^{4}}+\|\boldsymbol{u}_h^n\|_{L^4}\|\phi^n-\phi^{n-1}\|_{L^4}\right)\|\nabla e_{\phi}^n\|\\
				&~+C\epsilon\|\boldsymbol{u}_h^n\|_{L^6}\|e_{\phi}^n\|_{L^3}\|\nabla e_{\phi}^n\|\\
				\leq &~C\epsilon\|\nabla\phi^n\|(\|\boldsymbol{u}^n-\boldsymbol{P}_h\boldsymbol{u}^n\|_{H^1}+\|\nabla e_{\boldsymbol{u}}^n\|)\|\nabla e_{\phi}^n\|\\
				&~+C\epsilon\|\nabla\boldsymbol{u}_h^n\|\|\nabla(\phi^n-\phi^{n-1})\|\|\nabla e_{\phi}^n\|+C\epsilon\|\nabla\boldsymbol{u}_h^n\|\|\nabla e_{\phi}^n\|\|\nabla e_{\phi}^n\|\\
				\leq &~C\epsilon\|\nabla\phi^n\|(\|\boldsymbol{u}^n-\boldsymbol{P}_h\boldsymbol{u}^n\|_{H^1}+\|\nabla e_{\boldsymbol{u}}^n\|)\|\nabla e_{\phi}^n\|\\
				&~+C\tau\epsilon\|\nabla\boldsymbol{u}_h^n\|\|\nabla\phi_{t}\|\|\nabla e_{\phi}^n\|+C\epsilon\|\nabla\boldsymbol{u}_h^n\|\|\nabla e_{\phi}^n\|\|\nabla e_{\phi}^n\|.
			\end{aligned}
		\end{equation}
		By Lemma \ref{lemma_phi_Linf_boundedness} and \ref{lemma_dtau_phi_Linf_bounded}, there holds
		\begin{equation}
			|f(\phi^n)-f(\phi_h^n,\phi_h^{n-1})|\leq C(1+\epsilon^{-2\alpha_3})|\phi^n-\phi_h^n|+\tau|d_{\tau}\phi_h^n|.
		\end{equation}
		Therefore, we obtain the bound of the term $T_{c5}$.
		Next, substituting the above inequalities $T_{c1}\sim T_{c5}$ into  \eqref{eq_wvarphi_boundedness0011}, multiplying the resulting inequality by $2\tau$, taking the summation of the result from $n=1$ to $m$, and applying the Cauchy--Schwarz inequality, we obtain
		\begin{equation}
			\label{eq_summ_epsilonephi_result}
			\begin{aligned}
				\epsilon\|e_{\phi}^m\|^2+&2\epsilon\sum_{n=1}^m\|e_{\phi}^n-e_{\phi}^{n-1}\|^2+2\tau\sum_{n=1}^m\|e_{\mu}^n\|^2\\
				\leq&~C(\epsilon^{-4\alpha_3-1}+1)\tau\sum_{n=1}^m(\|\phi^n-R_h\phi^n\|^2+\|e_{\phi}^n\|^2+\|\nabla e_{\mu}^n\|^2)\\
				&~+C\tau\epsilon\sum_{n=1}^m(\|d_{\tau}R_h\phi^n-\partial_t\phi^n\|^2+\|\nabla e_{\phi}^n\|^2+\epsilon^{-\max\{8\sigma_1+7,2\sigma_3,2\sigma_1+2\sigma_2+2\}}\tau^2)\\
				&~+C\tau\sum_{n=1}^m(\|d_{\tau}(R_h\phi^n-\phi^n)\|_{-1,h}^2+\|R_1^n\|_{-1,h}^2+\|\mathcal{G}_h\mu^n-\mu^n\|_{-1,h}^2)\\
				&~+C\tau\sum_{n=1}^m(\|\boldsymbol{u}^n-\boldsymbol{P}_h\boldsymbol{u}^n\|_{H^1}^2+\|\nabla e_{\boldsymbol{u}}^n\|^2+\|\nabla e_{\mu}^n\|^2+\rho_0(\epsilon,d)\tau^2).
			\end{aligned}
		\end{equation}
		Thus, the error estimate for this term $\|d_{\tau}R_h\phi^n-\partial_t\phi^n\|^2$ in the inequality \eqref{eq_summ_epsilonephi_result} is as follows:
		\begin{equation}
			\begin{aligned}
				&\tau\sum_{n=1}^m\|d_{\tau}R_h\phi^n-\partial_t\phi^n\|^2\\
				\leq &~\tau\sum_{n=1}^m(\|d_{\tau}R_h\phi^n-d_{\tau}\phi^n\|^2+\|d_{\tau}\phi^n-\partial_t\phi^n\|^2)\\
				\leq&~\frac{1}{\tau}\sum_{n=1}^m\int_{\Omega}\left(\left|\int_{t^{n-1}}^{t^n}\partial_t(R_h\phi^n-\phi^n)dt\right|^2+\left|\int_{t^{n-1}}^{t^n}(t-t^n)\partial_{tt}\phi dt\right|^2\right)dx\\
				\leq&~\int_0^T\|\partial_t(R_h\phi-\phi)\|^2dt+\tau^2\int_0^T\|\partial_{tt}\phi\|^2dt.
			\end{aligned}
		\end{equation}
		Using the Sobolev embedding inequality,  \eqref{eq_Ritz_projections_boundedness0004} and the Lemmas \ref{lemma_stokes_projection001}, \ref{lemma_R1_Hneg_L2_boundedness0001}, \ref{lemma_nabla_emu} and \ref{lemma_eu_ephiHnev1000001}, we can see that the result \eqref{eq_L2_phimu_boundedness} of this lemma holds.
	\end{proof}

	From the above Lemmas, the triangle inequality, and error estimates of the projections, we complete the proof of the main convergence result in Theorem \ref{theorem_main_result}. For the error estimate of the pressure, we obtain it directly from the preceding error estimation analysis and the inf-sup condition as follows:
	\begin{equation}
		\tau\sum_{n=1}^m\|p^n-p_h^n\|^2\leq \frac{C(T)}{\epsilon^{\alpha_5}}(\tau^2+h^{2(r+1)}).
	\end{equation}
	From the above lemmas, we can derive all the conclusions.
	\section{Numerical results}\label{section_numerical_results}
	In this section, we carry out several numerical experiments for the following purposes: (i) to verify the optimal error estimates, both with and without smooth solutions; (ii) to demonstrate the uniform convergence order with respect to the width of the interface $\epsilon$; (iii) to confirm the unconditional stability of the energy; and (iv) to illustrate the effectiveness of the proposed method for practical problems. Unless otherwise stated, we use $P_1$ elements for the phase field, chemical potential and pressure,
	and $\boldsymbol{P}_2$ elements for velocity. Tables \ref{Table_L2_phi_mu_error_convRates_P2P2P3P2} and \ref{Table_L2_u_p_error_convRates_P2P2P3P2} additionally report the results for $P_2 \times P_2 \times \boldsymbol{P}_3 \times P_2$. The
	time discretization is first order, and the temporal error is controlled by a separate time-step refinement
	check in each spatial convergence test. Since the fully discrete scheme \eqref{eq_CHNS_variational_form_phi_discretization}-\eqref{eq_CHNS_variational_form_incompressible_u_discretization} is nonlinear, we solve the system iteratively by linearizing the nonlinear term using Newton's method.
	\begin{table}[h!]
		\centering
		\fontsize{8}{8}
		\begin{threeparttable}
			\caption{Errors and convergence rates for phase field and chemical potential with a smooth manufactured solution when $r=1$.}
			\label{Table_L2_phi_mu_error_convRates_FirstOrder}
			\begin{tabular}{c|c|c|c|c|c|c|c|c}
				\toprule
				\multirow{2.5}{*}{$h$}  & \multicolumn{2}{c}{$\left\|\phi^N - \phi_h^N\right\|_{L^2}$} & \multicolumn{2}{c}{$\left\|\mu^N - \mu_h^N\right\|_{L^2}$}  & \multicolumn{2}{c}{$\left\|\phi^N - \phi_h^N\right\|_{H^1}$} & \multicolumn{2}{c}{$\left\|\mu^N - \mu_h^N\right\|_{H^1}$} \cr
				\cmidrule(lr){2-3} \cmidrule(lr){4-5} \cmidrule(lr){6-7}  \cmidrule(lr){8-9}  
				& error & rate & error & rate & error & rate & error & rate\cr
				\midrule
				1/4   	& 2.03e$-$04 &  -     & 4.96e$-$02  & -     & 8.09e$-$03 & -     &  2.26& - \cr
				1/8   	& 4.82e$-$05 &  2.02  & 1.22e$-$02  & 2.02  & 4.33e$-$03 & 0.92  &  1.19& 0.90\cr
				1/16  	& 1.17e$-$05 &  2.00  & 3.10e$-$03  & 2.00  & 2.17e$-$03 & 1.00  &  0.59& 1.00\cr
				1/32    & 2.91e$-$06 &  2.01  & 7.77e$-$04  & 1.98  & 1.08e$-$03 & 1.00  &  0.30& 1.00\cr
				1/64  	& 7.28e$-$07 &  1.99  & 1.95e$-$04  & 2.00  & 5.45e$-$04 & 1.00  &  0.15& 1.00\cr
				\bottomrule
			\end{tabular}
		\end{threeparttable}
	\end{table}
	\begin{table}[h!]
		\centering
		\fontsize{8}{8}
		\begin{threeparttable}
			\caption{Errors and convergence rates for velocity and pressure with a smooth manufactured solution when $r=1$.}
			\label{Table_L2_u_p_error_convRates_FirstOrder}
			\begin{tabular}{c|c|c|c|c|c|c|c|c}
				\toprule
				\multirow{2.5}{*}{$h$}  & \multicolumn{2}{c}{$\left\|\boldsymbol{u}^N - \boldsymbol{u}_h^N \right\|_{\boldsymbol{L}^2}$} & \multicolumn{2}{c}{$\left\|\boldsymbol{u}^N  - \boldsymbol{u}_h^N \right\|_{\boldsymbol{H}^1}$}  & \multicolumn{2}{c}{$\left\|p^N  -p_h^N \right\|_{L^2}$} & \multicolumn{2}{c}{$\left\|p^N  - p_h^N \right\|_{H^1}$} \cr
				\cmidrule(lr){2-3} \cmidrule(lr){4-5} \cmidrule(lr){6-7}   \cmidrule(lr){8-9}  
				& error & rate & error & rate & error & rate & error & rate  \cr
				\midrule
				1/4  & 3.06e$-$04 &  -    & 1.55e$-$02   & -     & 4.88e$-$03 & -    & 2.60e$-$02  & -\cr
				1/8  & 5.37e$-$05 &  2.51 & 4.39e$-$03  & 1.82  & 7.19e$-$04 & 2.76 & 6.14e$-$03  & 2.08\cr
				1/16 & 7.59e$-$06 &  2.82 & 1.12e$-$03  & 1.97  & 1.13e$-$04 & 2.66 & 2.42e$-$03  & 1.35\cr
				1/32 & 9.80e$-$07 &  2.95 & 2.83e$-$04  & 1.99  & 2.04e$-$05 & 2.48 & 1.13e$-$03  & 1.10\cr
				1/64 & 1.23e$-$07 &  2.99 & 7.08e$-$05  & 2.00  & 4.20e$-$06 & 2.27 & 5.51e$-$04  & 1.03\cr
				\bottomrule
			\end{tabular}
		\end{threeparttable}
	\end{table}
	\begin{table}[h!]
		\centering
		\fontsize{8}{8}
		\begin{threeparttable}
			\caption{Errors and convergence rates for phase field and chemical potential with a smooth manufactured solution when $r=2$.}
			\label{Table_L2_phi_mu_error_convRates_P2P2P3P2}
			\begin{tabular}{c|c|c|c|c|c|c|c|c}
				\toprule
				\multirow{2.5}{*}{$h$}  & \multicolumn{2}{c}{$\left\|\phi^N - \phi_h^N\right\|_{L^2}$} & \multicolumn{2}{c}{$\left\|\mu^N - \mu_h^N\right\|_{L^2}$}  & \multicolumn{2}{c}{$\left\|\phi^N - \phi_h^N\right\|_{H^1}$} & \multicolumn{2}{c}{$\left\|\mu^N - \mu_h^N\right\|_{H^1}$} \cr
				\cmidrule(lr){2-3} \cmidrule(lr){4-5} \cmidrule(lr){6-7}  \cmidrule(lr){8-9}  
				& error & rate & error & rate & error & rate & error & rate\cr
				\midrule
				1/4   	& 3.53e$-$04 &  -      & 8.42e$-$01  & -     & 1.23e$-$02 & -     &  1.42e$+$01 & - \cr
				1/8   	& 5.23e$-$05 &  2.76   & 1.09e$-$01  & 2.94  & 3.31e$-$03 & 1.90  &  3.67e$+$00  & 1.95\cr
				1/16   	& 6.89e$-$06 &  2.92   & 1.39e$-$02  & 2.98  & 8.35e$-$04 & 1.99  &  9.20e$-$01  & 2.00\cr
				1/32  	& 8.75e$-$07 &  2.98   & 1.75e$-$03  & 2.99  & 2.10e$-$04 & 1.99  &  2.30e$-$01  & 2.00\cr
				1/64    & 1.10e$-$07 &  2.99   & 2.19e$-$04  & 3.00  & 5.26e$-$05 & 2.00  &  5.78e$-$02 &  2.00\cr
				\bottomrule
			\end{tabular}
		\end{threeparttable}
	\end{table}
	\begin{table}[h!]
		\centering
		\fontsize{8}{8}
		\begin{threeparttable}
			\caption{Errors and convergence rates for velocity and pressure with a smooth manufactured solution when $r=2$.}
			\label{Table_L2_u_p_error_convRates_P2P2P3P2}
			\begin{tabular}{c|c|c|c|c|c|c|c|c}
				\toprule
				\multirow{2.5}{*}{$h$}  & \multicolumn{2}{c}{$\left\|\boldsymbol{u}^N - \boldsymbol{u}_h^N \right\|_{\boldsymbol{L}^2}$} & \multicolumn{2}{c}{$\left\|\boldsymbol{u}^N  - \boldsymbol{u}_h^N \right\|_{\boldsymbol{H}^1}$}  & \multicolumn{2}{c}{$\left\|p^N  -p_h^N \right\|_{L^2}$} & \multicolumn{2}{c}{$\left\|p^N  - p_h^N \right\|_{H^1}$} \cr
				\cmidrule(lr){2-3} \cmidrule(lr){4-5} \cmidrule(lr){6-7}   \cmidrule(lr){8-9}  
				& error & rate & error & rate & error & rate & error & rate  \cr
				\midrule
				1/4   	& 1.06e$-$03 &  -     & 2.23e$-$02  & -     & 1.88e$-$01 & -    & 8.74e$-$01  & -\cr
				1/8   	& 7.09e$-$05 &  3.99  & 3.01e$-$03  & 2.89  & 2.35e$-$02 & 3.00 & 1.20e$-$01  & 2.86\cr
				1/16   	& 4.29e$-$06 &  4.01  & 3.73e$-$04  & 3.01  & 2.94e$-$03 & 3.00 & 1.75e$-$02  & 2.78\cr
				1/32  	& 2.70e$-$07 &  4.01  & 4.63e$-$05  & 3.01  & 3.68e$-$04 & 3.00 & 2.64e$-$03  & 2.73\cr
				1/64  	& 1.70e$-$08 &  3.97  & 5.83e$-$06  & 3.01  & 4.59e$-$05 & 3.01 & 4.19e$-$04  & 2.65\cr
				\bottomrule
			\end{tabular}
		\end{threeparttable}
	\end{table}
	\begin{table}[h]
		\centering
		\small
		\begin{threeparttable}
			\caption{Errors and convergence rates for phase field and chemical potential without a smooth manufactured solution.}
			\label{Table_L2_phi_mu_error_convRates_FirstOrder0002}
			\begin{tabular}{c|c|c|c|c|c|c|c|c|c}
				\toprule
				\multirow{2.5}{*}{$h$}  & \multirow{2.5}{*}{$h/2$}  & \multicolumn{2}{c}{$\left\|\phi^N_h - \phi^N_{h/2}\right\|_{L^2}$} & \multicolumn{2}{c}{$\left\|\mu^N_h - \mu^N_{h/2}\right\|_{L^2}$}  & \multicolumn{2}{c}{$\left\|\phi^N_h - \phi^N_{h/2}\right\|_{H^1}$} & \multicolumn{2}{c}{$\left\|\mu^N_h - \mu^N_{h/2}\right\|_{H^1}$} \cr
				\cmidrule(lr){3-4} \cmidrule(lr){5-6} \cmidrule(lr){7-8}  \cmidrule(lr){9-10}  
				& & error & rate & error & rate & error & rate & error & rate\cr
				\midrule
				1/4   &  1/8   & 4.75e$-$02 &  -     & 1.43e$+$00  & -      & 7.48e$-$01 & -     & 2.36e$+$01  & -      \cr
				1/8  &	 1/16  & 1.23e$-$02 & 1.95   & 3.70e$-$01  & 1.95   & 3.80e$-$01 & 0.98  & 1.32e$+$01  & 0.84   \cr
				1/16  &  1/32  & 3.14e$-$03 & 1.97   & 9.16e$-$02  & 2.02   & 1.94e$-$01 & 0.97  & 6.55e$+$00  & 1.01   \cr
				1/32  &	 1/64  & 7.89e$-$04 & 1.99   & 2.34e$-$02  & 1.97   & 9.75e$-$02 & 1.00  & 3.31e$+$00  & 0.99   \cr
				\bottomrule
			\end{tabular}
		\end{threeparttable}
	\end{table}
	
	\begin{table}[h]
		\centering
		\small
		\begin{threeparttable}
			\caption{Errors and convergence rates for velocity and pressure without a smooth manufactured solution.}
			\label{Table_L2_u_p_error_convRates_FirstOrder0002}
			\begin{tabular}{c|c|c|c|c|c|c|c|c|c}
				\toprule
				\multirow{2.5}{*}{$h$}  & \multirow{2.5}{*}{$h/2$}  & \multicolumn{2}{c}{$\left\|\boldsymbol{u}^N_h - \boldsymbol{u}^N_{h/2}\right\|_{\boldsymbol{L}^2}$} & \multicolumn{2}{c}{$\left\|p^N_h - p^N_{h/2}\right\|_{L^2}$}  & \multicolumn{2}{c}{$\left\|\boldsymbol{u}^N_h - \boldsymbol{u}^N_{h/2}\right\|_{\boldsymbol{H}^1}$} & \multicolumn{2}{c}{$\left\|p^N_h - p^N_{h/2}\right\|_{H^1}$} \cr
				\cmidrule(lr){3-4} \cmidrule(lr){5-6} \cmidrule(lr){7-8}  \cmidrule(lr){9-10}  
				& & error & rate & error & rate & error & rate & error & rate\cr
				\midrule
				1/4   &  1/8   & 8.30e$-$03 & -     & 2.93e$-$02  & -       & 6.13e$-$02 &  -     & 4.26e$-$01  & -         \cr
				1/8   &	 1/16  & 5.93e$-$04 & 3.81  & 4.35e$-$03  & 2.75    & 4.52e$-$03 & 3.76   & 1.73e$-$01  & 1.31      \cr
				1/16  &  1/32  & 7.87e$-$05 & 2.91  & 9.41e$-$04  & 2.21    & 6.21e$-$04 & 2.86   & 8.46e$-$02  & 1.03      \cr
				1/32  &	 1/64  & 1.04e$-$05 & 2.91  & 2.24e$-$04  & 2.07    & 9.24e$-$05 & 2.75   & 4.19e$-$02  & 1.01      \cr
				\bottomrule
			\end{tabular}
		\end{threeparttable}
	\end{table}
	
	\subsection{Convergence tests with or without manufactured solutions}
	The first two experiments are intended to validate the optimal convergence rates, with or without manufactured solutions, for some fixed suitable interfacial width parameters. In the first experiment, the manufactured solutions, adopted from \cite{2026_GaoHaijun_Optimal_L2_error_estimates_of_a_decoupled_FEM_for_the_Cahn_Hilliard_Navier_Stokes_equations}, take the following form:
	\begin{equation}
		\begin{cases} 
			\phi(t, x, y) = 2 + \sin(t) \cos(\pi x) \cos(\pi y), \\
			\boldsymbol{u}(t, x, y) = \left[ \pi \sin(\pi x)^2 \sin(2\pi y) \sin(t), -\pi \sin(\pi y)^2 \sin(2\pi x) \sin(t) \right]^T, \\
			p(t, x, y) = \cos(\pi x) \sin(\pi y) \sin(t).
		\end{cases}
	\end{equation}
	Simulations are performed in the domain $\Omega = [0, 1] \times [0, 1]$ with terminal time $T = 0.01$. The physical parameters are chosen as $M = 0.004$, $\epsilon = 0.04$, $\nu = 0.1$, and $\gamma = 1$. The time step is set to $\tau = 0.1h^3$ in order to balance the convergence rates in time and space. The corresponding numerical results are presented in Tables \ref{Table_L2_phi_mu_error_convRates_FirstOrder} and \ref{Table_L2_u_p_error_convRates_FirstOrder}. We further test the case where the finite element spaces are chosen as $P_2\times P_2\times\mathbf{P}_3\times P_2$ and $T=0.1$, $M = 0.004$, $\epsilon = 0.01$, $\nu = 1$, and $\gamma = 1$, and the corresponding numerical results are presented in Tables \ref{Table_L2_phi_mu_error_convRates_P2P2P3P2} and \ref{Table_L2_u_p_error_convRates_P2P2P3P2}.
	
	Furthermore, numerical experiments are conducted for problems without exact solutions, where the convergence behavior is assessed by comparing solutions on successively refined meshes. For this experiment, where no manufactured solution is available, the initial condition is prescribed as follows:
	\begin{equation}
		\begin{cases} 
			\phi(0, x, y) = \cos(\pi x) \cos(\pi y), \\
			\boldsymbol{u}(0, x, y) = \left[ \pi \sin(\pi x)^2 \sin(2\pi y), -\pi \sin(\pi y)^2 \sin(2\pi x)\right]^T, \\
		\end{cases}
	\end{equation}
	in the domain $\Omega = [0, 1] \times [0, 1]$ with terminal time $T = 0.1$. In this experiment, the Cauchy error is used to assess the convergence behavior. Specifically, for some variable $\xi$, the Cauchy error is defined by comparing the numerical solutions obtained on two successive refined meshes with spatial mesh sizes $h$ and $h/2$ as
	\begin{equation}
		\|e_\xi\|_* := \|\xi_h-\xi_{h/2}\|_*,
	\end{equation}
	where $\|\cdot\|_*$ denotes the adopted error norm. The physical parameters are chosen as $M = 4e-5$, $\epsilon = 0.04$, $\nu = 1$, and $\gamma = 1$, and the time step is set to $\tau = h^3$. Tables \ref{Table_L2_phi_mu_error_convRates_FirstOrder0002} and \ref{Table_L2_u_p_error_convRates_FirstOrder0002} present the errors and convergence rates.
	
	From Tables \ref{Table_L2_phi_mu_error_convRates_FirstOrder} and \ref{Table_L2_u_p_error_convRates_FirstOrder} (with manufactured solutions) and Tables \ref{Table_L2_phi_mu_error_convRates_FirstOrder0002} and \ref{Table_L2_u_p_error_convRates_FirstOrder0002} (without manufactured solutions), we can observe the optimal convergence rates for all variables established theoretically in Theorem \ref{theorem_main_result}. Specifically, when linear finite elements are used for the phase field, chemical potential, and pressure, optimal second-order convergence in the $L^2$-norm and first-order convergence in the $H^1$-norm are achieved. The phase field, chemical potential and pressure exhibit second-order \(L^2\) convergence and first-order \(H^1\) convergence in the \(r = 1\) tests. For velocity, the experiments exhibit approximately third-order \(L^2\) and second-order \(H^1\) convergence. The higher observed \(L^2\) velocity rate is consistent with the approximation capacity of the quadratic velocity space. For \(r = 1\), the current coupled analysis targets a second-order \(L^2\) velocity bound; the claim of polynomial dependence on \(\epsilon^{-1}\) remains conditional on closing the coupled stability estimate.
	When quadratic elements are used for the velocity, optimal third-order convergence in the $L^2$ norm and second-order convergence in the $H^1$ norm are attained.
	
	\subsection{Verification of "robustness" with respect to the interface width}
	To investigate the robustness of the proposed analysis with respect to the interfacial width parameter $\epsilon$, we perform a series of numerical experiments for CHNS system \eqref{eq_CHNS}-\eqref{eq_CHNS_initial_boundary_conditions} with different values of $\epsilon$. Except for the varying interface width parameter, the relevant parameters and the known manufactured solution are the same as in the first numerical experiment. 
	As shown by the results in Tables~\ref{Table_L2_u_p_error_convRates_FirstOrder_epsilon1}--\ref{Table_L2_u_p_error_convRates_FirstOrder_epsilon4}, the numerical errors and convergence rates remain stable for different values of \(\epsilon\), which confirms that the proposed error estimates exhibit no significant numerical variation with respect to the polynomial dependence. Therefore, this verifies that the theoretical analysis in this paper is relatively "robust" with respect to the interface width parameter \(\epsilon\).
	\begin{table}[h!]
		\centering
		\fontsize{10}{10}
		\begin{threeparttable}
			\caption{Errors and convergence rates for phase field in $\left\|\phi^N - \phi^N_h\right\|_{L^2}$ with different $\epsilon$.}
			\label{Table_L2_u_p_error_convRates_FirstOrder_epsilon1}
			\begin{tabular}{c|c|c|c|c|c|c|c|c|c|c}
				\toprule
				\multirow{2.5}{*}{$h$}  & \multicolumn{2}{c}{$\epsilon=0.1$} & \multicolumn{2}{c}{$\epsilon=0.05$}   & \multicolumn{2}{c}{$\epsilon=0.01$} & \multicolumn{2}{c}{$\epsilon=0.005$} & \multicolumn{2}{c}{$\epsilon=0.001$} \cr
				\cmidrule(lr){2-3} \cmidrule(lr){4-5} \cmidrule(lr){6-7}   \cmidrule(lr){8-9}  \cmidrule(lr){10-11}    
				& error & rate & error & rate & error & rate & error & rate & error & rate  \cr
				\midrule
				1/8  & 1.47e$-$02 &  -   & 1.44e$-$02 & -    & 1.44e$-$02 & -    & 1.43e$-$02 & -    &  1.43e$-$02 & -       \cr
				1/16 & 3.76e$-$03 & 1.96 & 3.70e$-$03 & 1.96 & 3.69e$-$03 & 1.96 & 3.69e$-$03 & 1.96 &  3.68e$-$03 & 1.96  \cr
				1/32 & 9.48e$-$04 & 1.99 & 9.34e$-$04 & 1.99 & 9.29e$-$04 & 1.99 & 9.29e$-$04 & 1.99 &  9.29e$-$04 & 1.99  \cr
				1/64 & 2.38e$-$04 & 2.00 & 2.34e$-$04 & 2.00 & 2.33e$-$04 & 2.00 & 2.33e$-$04 & 2.00 &  2.33e$-$04 & 2.00  \cr
				\bottomrule
			\end{tabular}
		\end{threeparttable}
	\end{table}
	
	\begin{table}[h!]
		\centering
		\fontsize{10}{10}
		\begin{threeparttable}
			\caption{Errors and convergence rates for chemical potential in $\left\|\mu^N - \mu^N_h\right\|_{L^2}$ with different $\epsilon$.}
			\label{Table_L2_u_p_error_convRates_FirstOrder_epsilon2}
			\begin{tabular}{c|c|c|c|c|c|c|c|c|c|c}
				\toprule
				\multirow{2.5}{*}{$h$}  & \multicolumn{2}{c}{$\epsilon=0.1$} & \multicolumn{2}{c}{$\epsilon=0.05$}   & \multicolumn{2}{c}{$\epsilon=0.01$} & \multicolumn{2}{c}{$\epsilon=0.005$} & \multicolumn{2}{c}{$\epsilon=0.001$} \cr
				\cmidrule(lr){2-3} \cmidrule(lr){4-5} \cmidrule(lr){6-7}   \cmidrule(lr){8-9}  \cmidrule(lr){10-11}    
				& error & rate & error & rate & error & rate & error & rate & error & rate  \cr
				\midrule
				1/8   	& 2.06e$+$00 &  -   & 4.07e$+$00 & -    & 2.03e$+$01 & -     & 4.05e$+$01   & -      & 2.02e$+$02    & -    \cr
				1/16   	& 5.36e$-$01 & 1.94 & 1.06e$+$00 & 1.94 & 5.28e$+$00  & 1.94 & 1.06e$+$01  &  1.94  & 5.28e$+$01    & 1.94    \cr
				1/32  	& 1.36e$-$01 & 1.98 & 2.68e$-$01 & 1.98 & 1.34e$+$00  & 1.98 & 2.67e$+$00	&  1.98  & 1.34e$+$01    & 1.98     \cr
				1/64  	& 3.40e$-$02 & 2.00 & 6.73e$-$02 & 2.00 & 3.35e$-$01  & 2.00 & 6.71e$-$01   &  2.00  & 3.35e$+$00    & 2.00     \cr
				\bottomrule
			\end{tabular}
		\end{threeparttable}
	\end{table}
	
	\begin{table}[h!]
		\centering
		\fontsize{10}{10}
		\begin{threeparttable}
			\caption{Errors and convergence rates for velocity field in $\left\|\boldsymbol{u}^N - \boldsymbol{u}^N_h\right\|_{\boldsymbol{L}^2}$ with different $\epsilon$.}
			\label{Table_L2_u_p_error_convRates_FirstOrder_epsilon3}
			\begin{tabular}{c|c|c|c|c|c|c|c|c|c|c}
				\toprule
				\multirow{2.5}{*}{$h$}  & \multicolumn{2}{c}{$\epsilon=0.1$} & \multicolumn{2}{c}{$\epsilon=0.05$}   & \multicolumn{2}{c}{$\epsilon=0.01$} & \multicolumn{2}{c}{$\epsilon=0.005$} & \multicolumn{2}{c}{$\epsilon=0.001$} \cr
				\cmidrule(lr){2-3} \cmidrule(lr){4-5} \cmidrule(lr){6-7}   \cmidrule(lr){8-9}  \cmidrule(lr){10-11}    
				& error & rate & error & rate & error & rate & error & rate & error & rate  \cr
				\midrule
				1/8   	& 5.36e$-$03 &  -     & 5.36e$-$03  & -      & 5.42e$-$03 & -     & 5.59e$-$03  & -      & 9.63e$-$03  & -    \cr
				1/16   	& 6.73e$-$04 &  2.99  & 6.74e$-$04  & 2.99   & 6.75e$-$04 & 3.00  & 6.80e$-$04  & 3.04   & 8.31e$-$04  & 3.53     \cr
				1/32  	& 8.49e$-$05 &  2.99  & 8.50e$-$05  & 2.99   & 8.50e$-$05 & 2.99  & 8.52e$-$05	& 3.00   & 9.02e$-$05  & 3.20     \cr
				1/64  	& 1.09e$-$05 &  2.97  & 1.09e$-$05  & 2.97   & 1.09e$-$05 & 2.97  & 1.09e$-$05  & 2.97   & 1.11e$-$05  & 3.03     \cr
				\bottomrule
			\end{tabular}
		\end{threeparttable}
	\end{table}
	
	\begin{table}[h!]
		\centering
		\fontsize{10}{10}
		\begin{threeparttable}
			\caption{Errors and convergence rates for pressure field in $\left\|p^N - p^N_h\right\|_{L^2}$ with different $\epsilon$.}
			\label{Table_L2_u_p_error_convRates_FirstOrder_epsilon4}
			\begin{tabular}{c|c|c|c|c|c|c|c|c|c|c}
				\toprule
				\multirow{2.5}{*}{$h$}  & \multicolumn{2}{c}{$\epsilon=0.1$} & \multicolumn{2}{c}{$\epsilon=0.05$}   & \multicolumn{2}{c}{$\epsilon=0.01$} & \multicolumn{2}{c}{$\epsilon=0.005$} & \multicolumn{2}{c}{$\epsilon=0.001$} \cr
				\cmidrule(lr){2-3} \cmidrule(lr){4-5} \cmidrule(lr){6-7}   \cmidrule(lr){8-9}  \cmidrule(lr){10-11}    
				& error & rate & error & rate & error & rate & error & rate & error & rate  \cr
				\midrule
				1/8   	& 8.39e$-$02 &  -    & 1.48e$-$01  & -      & 6.89e$-$01 & -     & 1.37e$+$00  & -      & 6.81e$+$00  & -        \cr
				1/16   	& 2.07e$-$02 & 2.02  & 3.81e$-$02  & 1.96   & 1.81e$-$01 & 1.93  & 3.60e$-$01  & 1.93   & 1.74e$+$00  & 1.92     \cr
				1/32  	& 5.20e$-$03 & 1.99  & 9.65e$-$03  & 1.98   & 4.60e$-$02 & 1.98  & 9.15e$-$02  & 1.98   & 4.56e$-$01  & 1.98     \cr
				1/64  	& 1.30e$-$03 & 2.00  & 2.42e$-$03  & 1.99   & 1.15e$-$02 & 1.99  & 2.30e$-$02  & 1.99   & 1.14e$-$01  & 1.99     \cr
				\bottomrule
			\end{tabular}
		\end{threeparttable}
	\end{table}
	\begin{Remark}
		In this numerical experiment that involves different interface widths $\epsilon$, while all other physical parameters remain unchanged, a smaller value of $\epsilon$ leads to a larger computational error in the chemical potential function $\mu$. However, this variation does not affect the convergence order of the overall discrete system. Therefore, in order to strike a balance between the error control of the various variables and the convergence order of the system, it is necessary to adjust the mobility coefficient $M$ accordingly when $\epsilon$ takes relatively small values. When mobility $M$ is set to $0.01$, $0.005$, $0.001$, $0.0005$, or $0.0001$, the interfacial thickness parameter $\epsilon$ corresponds to $0.1$, $0.05$, $0.01$, $0.005$, or $0.001$, respectively. 
	\end{Remark}
	\subsection{Shape relaxation with rotational boundary condition}\label{subsection_shape_relaxation}
\begin{figure}[h!]
		\centering
		\includegraphics[width=1\linewidth]{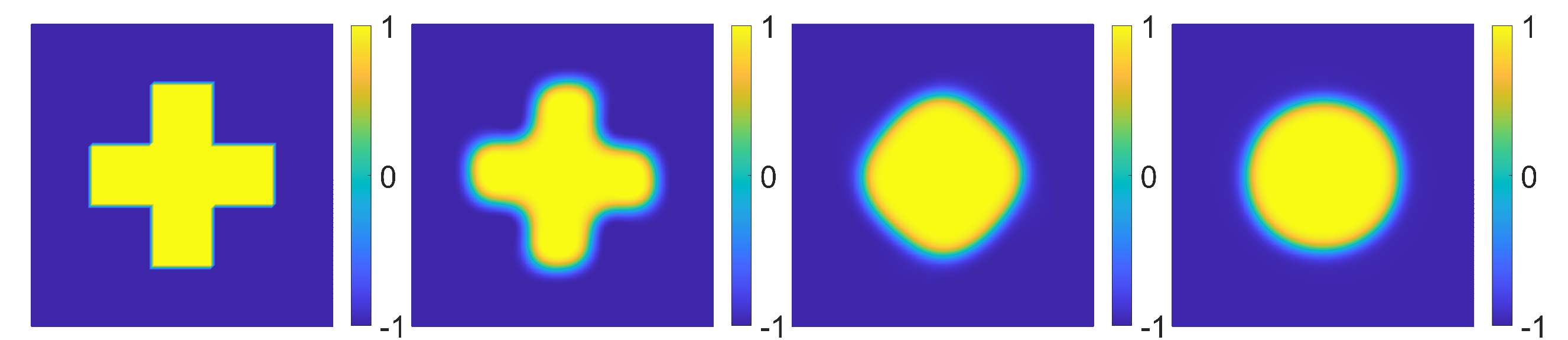}
		\caption{Snapshots of phase function at different times from left to right with $t=0, 0.1, 0.5, 5$, respectively.}
		\label{Fig_CHNS_NewEstimate_RotaBounCond_PhaseValue_dt001_20260721}
	\end{figure}

	In this experiment, we illustrate the effectiveness by simulating relaxation of a cross-shaped droplet toward a circular shape under the influence of rotational external forces. The computational domain is $\Omega=[0,1]\times[0,1]$ with final time $T=5$, and the rotational boundary conditions are prescribed as $\boldsymbol{u} = [y-0.5, \, -x+0.5]^{\textsf{T}}$ on $\partial \Omega$. The physical parameters are set as $M = 0.0001$, $\epsilon = 0.02$, $\nu = 1$, and $\gamma = 1$, and $h=1/64$, $\tau = 0.01$. The initial condition for the phase field is chosen such that $\phi = 1$ in a central polygonal cross-shaped region and $\phi = -1$ elsewhere, and for the velocity it is $\boldsymbol{u}^0 = [y-0.5, \, -x+0.5]^{\textsf{T}}$. As shown in Fig. \ref{Fig_CHNS_NewEstimate_RotaBounCond_PhaseValue_dt001_20260721}, driven by surface tension and imposed external forcing, the initial cross-shaped interface gradually evolves into a circular shape over time.
	
	\subsection{Mass conservation and unconditional energy stability}
\begin{figure}[h!]
		\centering
		\includegraphics[width=0.85\linewidth]{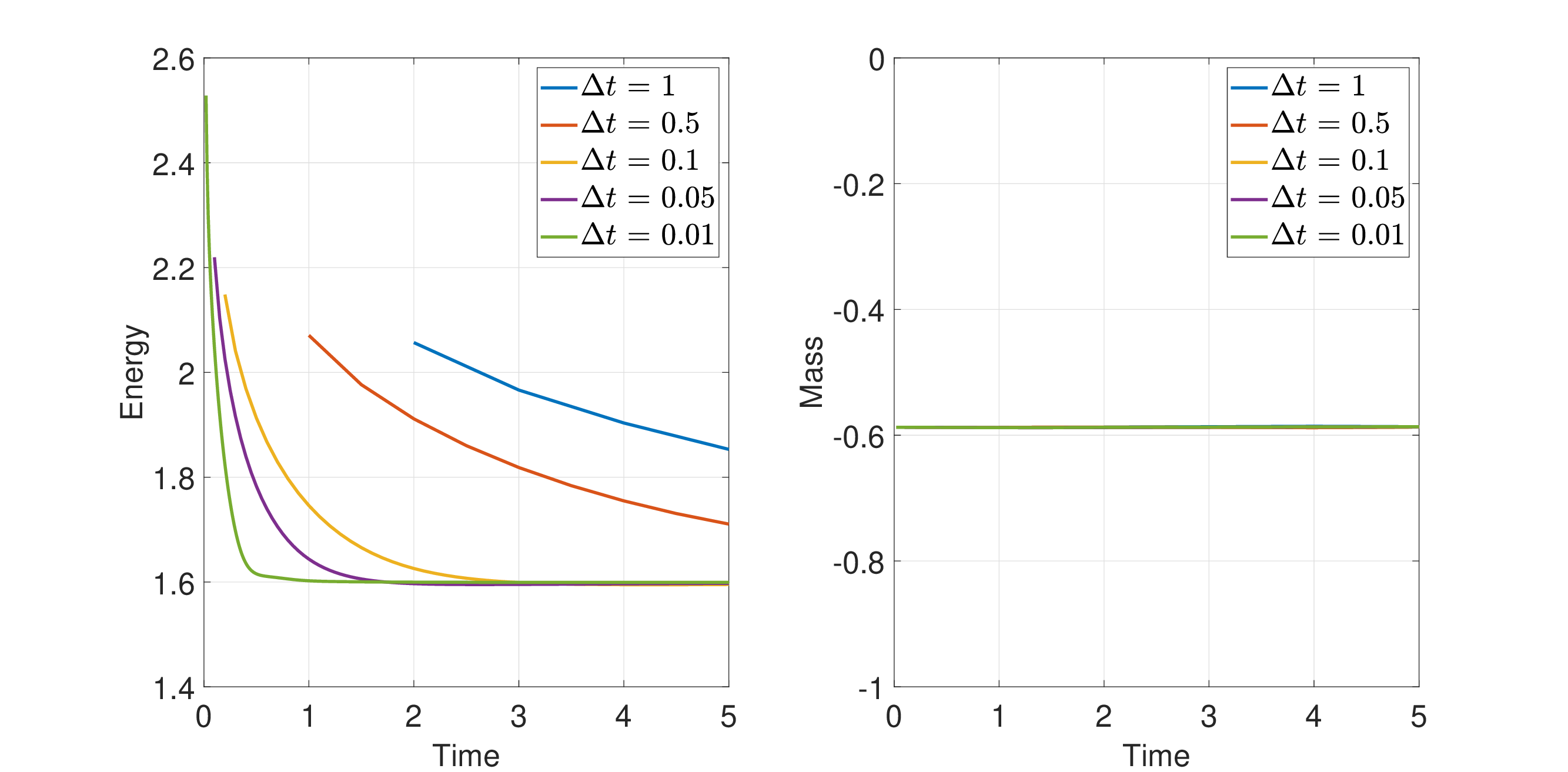}
		\caption{Evolution of the discrete mass and energy.}
		\label{Fig_CHNS_NewEstimate_RotaBounCond_Energy_Mass_20260721}
	\end{figure}

	The last numerical test is intended to validate the unconditional energy stability and mass conservation properties. To this end, we perform the aforementioned shape-relaxation example in subsection \ref{subsection_shape_relaxation} using the same parameter setup, except for the different time-step sizes $\tau$. The discrete energy $\mathcal{J}_{\epsilon}(\phi_h^n,\boldsymbol{u}_h^n)$ is defined in \eqref{eq_discrete_energy} and the discrete mass are
	\begin{equation}
		M_h^n := (\phi_h^n,1), \; n=0,1,\cdots, N.
	\end{equation}
	It can be clearly observed from Fig. \ref{Fig_CHNS_NewEstimate_RotaBounCond_Energy_Mass_20260721} that the discrete mass $M_h^n$ is well preserved throughout the simulation, while the discrete energy $\mathcal{J}_{\epsilon}(\phi_h^n,\boldsymbol{u}_h^n)$ monotonically decreases as time evolves. These results confirm that the proposed scheme preserves the mass conservation property and satisfies the discrete energy dissipation law, which is fully consistent with the theoretical analysis.
	\section{Conclusion}\label{section_conclusions}	
	In this paper, we have investigated a first-order convex-splitting fully discrete mixed finite element scheme for the CHNS system, with particular emphasis on the dependence of the error constants on the interfacial thickness parameter $\epsilon$. The analysis relies on three key ingredients. First, by constructing suitable Ritz quasi-projections and a Stokes quasi-projection, the leading error terms arising from the phase-field--flow coupling are cancelled exactly, thereby avoiding the appearance of the gradient of the phase-field error on the right-hand side of the error equations. Second, the \(H^{-1}\) approximation properties of these projections enable the phase-field error to be controlled in a negative norm, which allows the singular factor $\epsilon^{-1}$ associated with the nonlinear potential to be balanced against the discrete gradient energy. Third, the discrete spectral estimate for the perturbed Cahn--Hilliard operator avoids the crude Grönwall treatment of the nonlinear term that typically leads to exponentially growing error constants of the form $e^{c/\epsilon}$. To our knowledge, this provides the optimal $L^2$ error estimates with polynomial dependence on $\epsilon^{-1}$ for a mixed finite element discretization of the coupled CHNS system. The numerical experiments further confirm the predicted convergence rates and demonstrate that the convergence orders remain essentially unaffected by the interfacial thickness parameter $\epsilon$, while the error magnitudes exhibit only algebraic growth as $\epsilon$ decreases, in agreement with the theoretical results.
	
	Several issues remain open for future investigation. A natural direction is to extend the present framework to decoupled linear schemes and to SAV/IEQ-based formulations while preserving polynomial dependence on $\epsilon^{-1}$, which would make the analysis more suitable for large-scale computations. Higher-order temporal discretizations, such as second-order convex-splitting and BDF2 schemes, also deserve further study. Moreover, combining the present $\epsilon$-explicit analysis with \(H(\mathrm{div})\)-conforming finite element methods may lead to pressure-robust and Reynolds-semi-robust error estimates. Finally, extending the present analysis to the three-dimensional setting remains challenging due to the stronger nonlinear effects of the convective term. These issues will be considered in future work.
	
	\bibliographystyle{siamplain}
	\bibliography{references}
\end{document}